\documentclass[a4paper]{article}
\usepackage{amssymb,amsbsy,amsmath,amsfonts,amssymb,amscd}
\usepackage{latexsym}

\usepackage{latexsym}

\newcommand{\lie}[1] {\mathfrak{#1}}  
\newcommand{\bb}[1]{{\mathbb #1}}    

\newcommand{\bbR}{{\bb R}}

\newcommand{\bbZ}{{\bb Z}} 
\newcommand{\bbQ}{{\bb Q}}

\newcommand{\bbC}{{\bb C}}
\newcommand{\bbN}{{\bb N}}

\newcommand{\GL}{{\rm GL}}

\newcommand{\U}{{\rm U}}

\newcommand{\Aut}{{\rm Aut}}

\newcommand{\Out}{{\rm Out}}
\newcommand{\Inn}{{\rm Inn}}

\newcommand{\Ad}{{\rm Ad}}

\newcommand{\N}{{\rm N}}  
\newcommand{\Z}{{\rm Z}} 
\newcommand{\Fitt}{{\rm Fitt}} 

\newcommand{\C}{{\rm Z}} 

\newcommand{\A}{{\rm A}}

\newcommand{\bF}{{\mathbf{F}}}
\newcommand{\bC}{{\mathbf{C}}}
\newcommand{\bU}{{\mathbf{U}}}

\newcommand{\bH}{{\mathbf{H}}}

\newcommand{\bT}{{\mathbf{T}}}

\newcommand{\bS}{{\mathbf{S}}}

\newcommand{\cA}{{\mathcal{A}}}
\newcommand{\cB}{{\mathcal{B}}}
\newcommand{\cC}{{\mathcal{C}}}
\newcommand{\cO}{{\mathcal{O}}}
\newcommand{\cZ}{{\mathcal{Z}}}

\newcommand{\lu}{{\lie{u}}}
\renewcommand{\lg}{{\lie{g}}}

\newcommand{\id}{{\rm id}}

\newcommand{\rank}{{\rm rank }\, }
\newcommand{\ns}{{\triangleleft \,}}  

\newcommand{\ac}[1]{\overline{#1}}

{ \par  \medskip  \noindent  {\bf Definition} \hspace{0.5em} }%
{\par \medskip } 

\newtheorem{example1}{Example}[section]

{ \begin{example1} \rm }%
{ \end{example1}} 
\newenvironment {proof}%
{ \noindent {\em Proof. }}%
{\hspace*{\fill}$\Box$\par \medskip } 
\newenvironment{prf}[1]%
{ \noindent {\em #1 \hspace{0.5em}} }%
{\hspace*{\fill}$\Box$\par \medskip } 
\newenvironment {remark}%
{{\em Remark \hspace{0.5em}} }%
{\par \medskip }

\newtheorem{proposition}{Proposition}[section]
\newtheorem{definition1}[proposition]{Definition}
\newtheorem{theorem}[proposition]{Theorem}
\newtheorem{lemma}[proposition]{Lemma}
\newtheorem{corollary}[proposition]{Corollary}
\newenvironment{definition}%
{ \begin{definition1} \rm }%
{ \end{definition1}}

\input xy
\xyoption{all}

\newcommand{\bysame}{\hspace{1.5cm}}
\newcommand{\ra}{\rightarrow}
\renewcommand{\H}{{\mathrm H}}
\newcommand{\cD}{{\cal D}}
\newcommand{\cE}{{\cal E}}

\newcommand{\cG}{{\cal G}}
\newcommand{\cK}{{\cal K}}
\newcommand{\cS}{{\cal S}}
\newcommand{\cT}{{\cal T}}
\newcommand{\cN}{{\cal N}}
\newcommand{\cH}{{\cal H}}
\newcommand{\cR}{{\cal R}}

\newcommand{\tGam}{{\tilde{\Gamma}}} 
\newcommand{\tA}{{\tilde{A}}}
\newcommand{\tF}{{\tilde{F}}}

\newcommand{\Ga}{{\bb G}_a}
\newcommand{\Gm}{{\bb G}_m}

\font\Bbb=msbm10

\def\BZ{{\hbox{\Bbb Z}}}

\def\BN{\hbox{\Bbb N}}
\def\BR{{\hbox{\Bbb R}}}
\def\BC{{\hbox{\Bbb C}}}
\def\BQ{{\hbox{\Bbb Q}}}

\def\Ad{{\rm Ad}}

\def \Box{\lower .1 em
          \vbox{\hrule \hbox{\vrule \hskip .6 em \vrule height .6 em} \hrule}}
\def \Mid{\quad \vrule \quad}
\font \msbm= msbm10
\def \rfish{\mathbin {\hbox {\msbm \char'157}}}

\newcommand{\AGF}{A_{\Gamma \mid F}}
\newcommand{\tAGF}{A_{\tGam \mid \tF}}

\newcommand{\cAS}{\cA_\bS}
\newcommand{\cASid}{\cA_\bS^1}

\newcommand{\GLHC}{\GL\big(\H^*(\Gamma, \bbC)\big)} 
\newcommand{\KGone}{K(\Gamma,1)}

\title{ Automorphism groups of polycyclic-by-finite groups 
and arithmetic groups\newline}

\author{Oliver Baues \thanks{e-mail: baues@math.uni-karlsruhe.de} \\
Mathematisches Institut II\\ 
Universit\"at Karlsruhe\\ 
D-76128 Karlsruhe
\and 
Fritz Grunewald \thanks{e-mail: fritz@math.uni-duesseldorf.de}\\
Mathematisches Institut\\
Heinrich-Heine-Universit\"at\\
D-40225 D\"usseldorf}

\date{May 18, 2005}     

\bigskip
\bigskip
\bigskip
\bigskip

\begin{document}
\maketitle 
\begin{abstract} 
We show that the outer automorphism group of a 
polycyclic-by-finite group is an arithmetic group. 
This result follows from a detailed structural
analysis of the automorphism groups of such groups.
We use an extended version of the theory of the algebraic 
hull functor initiated by Mostow.
We thus make applicable refined methods from the theory of algebraic 
and arithmetic groups.
We also construct examples of polycyclic-by-finite 
groups which have an automorphism group which
does not contain an arithmetic group of finite index.
Finally we discuss applications of our results to the groups of homotopy
self-equivalences of $\KGone$-spaces and  obtain an extension 
of arithmeticity results of Sullivan in rational homotopy theory.
\newline
\hrule 

\hspace{1cm} \\

\end{abstract}

\medskip

\bigskip
\bigskip
\bigskip
\newpage
\tableofcontents 
\medskip

\bigskip
\bigskip

\newpage 
\section{Introduction} 

\subsection{The main results}

We write $\Aut(\Gamma)$ for the group of automorphisms of a group $\Gamma$.
The subgroup consisting of automorphisms induced by conjugations with elements
of $\Gamma$ is denoted by $\Inn_\Gamma$. It is normal in $\Aut(\Gamma)$ and
the quotient group 
$$ \Out(\Gamma):=\Aut(\Gamma)\big/\Inn_\Gamma$$ 
is called the
outer automorphism group of $\Gamma$. This paper is devoted to the detailed
study of the groups $\Aut(\Gamma)$ and $\Out(\Gamma)$ in case $\Gamma$ is
polycyclic-by-finite. Here we say that a group $\Gamma$ 
is $\cal E$-by-finite whenever $\cal E$ is a property of groups and when 
$\Gamma$ has a subgroup of finite index having property 
$\cal E$.  
As one of our main results we prove:

\begin{theorem}\label{teoa} 
For any polycyclic-by-finite group 
$\Gamma$, $\Out(\Gamma)$ 
is an arithmetic group.
\end{theorem} 

To explain the concept of an arithmetic group we recall that a 
{\it $\BQ$-defined linear algebraic group} $G$ is a subgroup $G \le
\GL(n,\BC)$
$(n\in\BN)$ which is also an affine algebraic set defined by polynomials with
rational coefficients in the natural coordinates of $\GL(n,\BC)$. If $R$ is a
subring of $\BC$ we put $G(R):=G\cap \GL(n,R)$.
We have $G=G(\BC)$.
Let $G$ be
a $\BQ$-defined linear algebraic group. A subgroup
$\Gamma
\le G(\BQ)$ is called an {\it arithmetic subgroup} of $G$ if $\Gamma$
is commensurable with $G(\BZ )$.
An abstract group $\Delta$ is called {\it arithmetic} if it is isomorphic to
an arithmetic subgroup of a $\BQ$-defined linear algebraic group.
Two subgroups $\Gamma_1,\, \Gamma_2$ of $\GL(n,\BC)$ are called
{\it commensurable}
if their intersection $\Gamma_1\cap \Gamma_2$ has finite index both in
$\Gamma_1$ and $\Gamma_2$. These definitions are taken from \cite{PR}. 

In case $\Gamma$ is a finitely generated nilpotent group,
Segal \cite{Segal2} observed that $\Out(\Gamma)$ contains an
arithmetic subgroup of finite index. Such a group need itself not be
arithmetic, see \cite{GP3} for examples. Our Theorem \ref{teoa} is a
strengthening of Segal's result even in this restricted case. 
The results of Sullivan \cite{Sullivan} 
imply for finitely generated nilpotent groups a
still weaker result on their outer automorphism groups. We shall come back to
this in Section \ref{applih}. 
A related result, for any polycyclic-by-finite group, is proved 
by Wehrfritz in \cite{Wehrfritz}. 
He shows that $\Out(\Gamma)$ admits a faithful representation
into $\GL(n,\bbZ)$, for some $n$,  in this general case.

A role model for these results  
is the case $\Gamma=\BZ^n$ ($n\in\BN$). Here we have
$$\Aut(\BZ^n)=\Out(\BZ^n)=\GL(n,\BZ)$$
and both $\Aut(\BZ^n)$ and $\Out(\BZ^n)$ are arithmetic groups. 
More generally, $\Aut(\Gamma)$ is arithmetic for every finitely generated
nilpotent group $\Gamma$, see \cite[Corollary 9, Chapter 6]{Segal}. This
result was obtained for torsion-free finitely generated
nilpotent groups by Auslander and Baumslag, see  \cite{AusBaumslag,Baumslag,Baumslag1}. 
We generalize it to the case of finitely generated
nilpotent-by-finite groups in Corollary \ref{coronilfi}. 
But this is not the general picture. 
We now describe a polycyclic group so that 
$\Aut(\Gamma)$ does not contain an arithmetic subgroup of finite index. To do
this we choose $d\in\BN$ not a square,  and let $K=\BQ(\sqrt d)$ be the
corresponding real quadratic number field and write $x\mapsto \bar x$ for the
non-trivial element of the Galois group  of $K$ over $\BQ$. 
Consider the subring ${\cal O}=\BZ+\BZ\sqrt d\le K$. By Dirichlet's unit
theorem we may choose a unit $\epsilon\in {\cal O}^*$ which is of infinite
order and satisfies $\epsilon \bar\epsilon=1$. Define 
\begin{equation}\label{dihe}
D_\infty:= \; \langle \, A,\,\tau \Mid \tau^2=(A\tau)^2=1 \, \rangle
\end{equation}
to be the infinite dihedral group. It is easy to see that $F={\cal O}\times
\BZ$ obtains the structure of a $\Gamma$-module by defining
\begin{equation}\label{modu}
A\cdot (m,n) =(\epsilon m,n),\quad \tau\cdot (m,n)=(\bar m,-n)\qquad
(m\in{\cal O},\, n\in\BZ).
\end{equation}   
We denote the corresponding split extension by
$$\Gamma(\epsilon):= \; F\rfish D_\infty.$$
The group $\Gamma(\epsilon)$ is polycyclic and even an arithmetic group, but
we have:
\begin{theorem}\label{teob} 
The automorphism group of\/ $\Gamma(\epsilon)$ does not
contain an arithmetic group of finite index.
\end{theorem} 

The statement of Theorem \ref{teob} is stronger than just saying that 
the automorphism group of $\Gamma(\epsilon)$ is not arithmetic, since there are
groups (see \cite{GP3}) which are not arithmetic but which contain a subgroup
of finite index which has this property.
Our examples show that the structure of the automorphism group varies
dramatically when $\Gamma$ is replaced by one of its subgroups of 
finite index. It is for example easy to see (also from results in
Subsection \ref{outline}) that the automorphism group of the subgroup of
$\Gamma(\epsilon)$ which is generated by $F$ and $A$ is an arithmetic group. 
More generally, as we will explain below, every polycyclic-by-finite group has
a finite index subgroup with an arithmetic automorphism group. 

Theorem \ref{teob} is complemented by a result of Merzljakov 
\cite{Merzljakov} who has proved,
for any polycyclic group $\Gamma$,  that 
$\Aut(\Gamma)$ has a faithful representation into $\GL(n,\bbZ)$ for some 
$n\in\BN$. This generalizes also to the case of 
polycyclic-by-finite groups, see \cite{Wehrfritz2}. 

Let $F= \Fitt(\Gamma)$ be the Fitting subgroup
of $\Gamma$. This is the maximal nilpotent normal subgroup of $\Gamma$. 
Since $F$ is characteristic in $\Gamma$ (that is normalized by every
automorphism of $\Gamma$),  we may define 
\begin{equation}\label{defiagf}
\A_{\Gamma |F} :=  \; \{\,  \phi \in \Aut(\Gamma) \Mid 
\phi|_{\, \Gamma/F} = \id_{\Gamma/F} 
\}. 
\end{equation}

As our main structural result on the automorphism group of a general 
 polycyclic-by-finite group we show:
\begin{theorem}\label{teoc} 
Let $\Gamma$ be a polycyclic-by-finite group. Then $\AGF$ is a normal 
subgroup of $\Aut(\Gamma)$ and $\AGF$ is an arithmetic group. 
There exists a finitely generated nilpotent group 
$B \leq \Aut(\Gamma)$ 
which consists of inner automorphisms 
such  that $$  \A_{\Gamma|F} \cdot B  \; \, \leq \; \Aut(\Gamma) $$
is a subgroup of finite index in the automorphism group of\/ $\Gamma$.
\end{theorem}

This theorem is a more precise version of \cite[Theorem 2, Chapter 8]{Segal}
and of certain similar theorems contained in \cite{Ausl1}. It has many
finiteness properties of $\Aut(\Gamma)$ as a consequence, 
see \cite[Chapter 8]{Segal}. The fact that $\Aut(\Gamma)$ is finitely 
presented was first proved by Auslander in \cite{Ausl1}. 


\subsection{Outline of the proofs and more results}\label{outline}


We shall describe now the general strategy in our proofs and more results on
the structure of the automorphism group $\Aut(\Gamma)$ of a general
polycyclic-by-finite group $\Gamma$.
We take the basic theory of polycyclic-by-finite groups for granted.
As a reference for our notation and results we use the book \cite{Segal}. 

A polycyclic-by-finite group $\Gamma$ has a maximal finite normal subgroup 
$\tau_\Gamma$. We say that $\Gamma$ is a {\it wfn-group} if
$\tau_\Gamma=\{1\}$. We first consider the case of 
polycyclic-by-finite groups which are wfn-groups and later reduce to this
case. 

Let now $\Gamma$ be a polycyclic-by-finite group which is a wfn-group. 
We use the construction from \cite{Baues2} of a $\BQ$-defined  solvable-by-finite 
linear algebraic group ${\bf H}_\Gamma$ which contains $\Gamma$ in its group of $\BQ$-rational 
points. The study of the functorial construction $\Gamma \mapsto \bH_\Gamma$ traces back to
work of Mostow \cite{Mostow1, Mostow2}. 
The algebraic group ${\bf H}_\Gamma$ has
special features which allow us to
identify the group of algebraic automorphisms of 
${\bf H}_\Gamma$, which we call $\Aut_a({\bf H}_\Gamma)$, with a 
$\BQ$-defined linear algebraic group $\cA_\Gamma$. In general, the group of 
algebraic automorphisms of a $\BQ$-defined algebraic group is an extension of
a linear algebraic by an arithmetic group, see \cite{BS}.
The functoriality of ${\bf H}_\Gamma$ 
leads to a natural embedding
\begin{equation}
\Aut(\Gamma) \, \le \, \Aut_a({\bf H}_\Gamma)=\cA_\Gamma.
\end{equation}
These steps are carried out in detail in Sections \ref{section:auta}, 
\ref{sechull}.

After preparations in Sections 
\ref{secthick}, \ref{secshad} we prove the following in Section \ref{secautoauto}.  
\begin{theorem}\label{teoagf}
Let $\Gamma$ be a polycyclic-by-finite 
group which is a wfn-group. 
Then $\Aut(\Gamma)$ is contained in $\cA_\Gamma(\BQ)$.
A subgroup of finite index in $\Aut(\Gamma)$ is contained in
$\cA_\Gamma(\BZ)$.
The group $\A_{\Gamma|F}$ is an arithmetic 
subgroup in its Zariski-closure in  $\cA_\Gamma$.
\end{theorem}

Theorem \ref{teoagf} is the main step in the proof of Theorem
\ref{teoc} which is given in Section \ref{generalpofi}. 

In Section \ref{subsecari} 
we develop, starting from results in \cite{GP3}, 
a theory of proving arithmeticity for linear groups which are related to
arithmetic groups in certain ways. This shows that the situation 
build up in Theorems \ref{teoagf}, \ref{teoc} 
implies that $\Out(\Gamma)$ is an arithmetic group for 
polycyclic-by-finite groups which are wfn-groups:
To construct a $\BQ$-defined algebraic group ${\cal O}_\Gamma$ which contains 
$\Out(\Gamma)$ as an arithmetic subgroup we start off by considering the
Zariski-closure ${\cal B}$ of $\A_{\Gamma|F}$ in $\cA_\Gamma$. We then take
the quotient group ${\cal B}/{\cal C}$ where ${\cal C}$ is the   
Zariski-closure of $\Inn_F$ in $\cA_\Gamma$. 
The image of $\A_{\Gamma|F}$ in ${\cal B}/{\cal C}$ is an arithmetic 
subgroup.
We show that 
$\Out(\Gamma)$ is a finite extension group of the 
quotient $\A_{\Gamma|F} \big/ \, \Inn_F$ by its central subgroup 
$D= (\A_{\Gamma|F}\cap \Inn_\Gamma) \big/\, \Inn_F$. 
We modify 
the algebraic group ${\cal B}\big/{\cal C}$ to obtain 
a $\BQ$-defined linear algebraic group ${\cal E}$ such that
the group $D$ is unipotent-by-finite, and $\cE$
still contains an isomorphic copy of $\A_{\Gamma|F} \big/\, \Inn_F$ as 
an arithmetic subgroup. The construction ensures that 
$(\A_{\Gamma|F}\cap \Inn_\Gamma) \big/ \, \Inn_F$ 
is arithmetic in its Zariski-closure $\cD$ in ${\cal E}$ and 
we can prove that     
$\Out(\Gamma)$ is arithmetic in a finite extension group 
${\cal O}_\Gamma$ of the quotient ${\cal E}/{\cD}$. The final definition of 
${\cal O}_\Gamma$ is contained in Section \ref{asec1}

To prepare for later applications we introduce some more notation. 
The group of inner automorphisms $\Inn_{{\bf H}_\Gamma}$ is 
by our constructions a $\BQ$-closed subgroup of the group of algebraic
automorphisms 
$\cA_\Gamma=\Aut_a({\bf H}_\Gamma)$. We call the quotient 
\begin{equation}
\Out_a({\bf H}_\Gamma):=\cA_\Gamma\big/\, \Inn_{{\bf H}_\Gamma}
\end{equation}
the algebraic outer automorphism group of ${\bf H}_\Gamma$.
It is again a 
$\BQ$-defined linear algebraic group, and we obtain a group homomorphism
\begin{equation}
\pi_\Gamma :\Out(\Gamma)\to \Out_a({\bf H}_\Gamma)
\end{equation}
by restricting the quotient homomorphism $\cA_\Gamma \ra  
\Out_a({\bf H}_\Gamma)$ to the subgroup $\Aut(\Gamma)\le\cA_\Gamma$.  
We then prove in Section \ref{asec1}:

\begin{theorem}\label{ftop} 
Let $\Gamma$ be a polycyclic-by-finite 
group which is a wfn-group. Then there is a $\BQ$-defined 
linear algebraic group ${\cal O}_\Gamma$ which contains $\Out(\Gamma)$ 
as an arithmetic subgroup and a $\BQ$-defined homomorphism 
$\pi_{\cO_\Gamma}:{\cal O}_\Gamma\to \Out_a({\bf H}_\Gamma)$ such that the
diagram
\begin{equation}
\xymatrix{
\Out(\Gamma) \ar[r]^{} \ar[d]_{\pi_\Gamma} 
& {\cal O}_\Gamma \ar[ld]^{\pi_{\cO_\Gamma}} \\
\Out_a({\bf H}_\Gamma) & 
}
\end{equation}
is commutative.
The kernel of $\pi_\Gamma$ is 
finitely generated, abelian-by-finite and is centralized by a finite
index subgroup of $\Out(\Gamma)$. If\/ $\Gamma$ is nilpotent-by-finite
then the kernel of $\pi_\Gamma$ is finite. 
\end{theorem}

We shall start the discussion of the general situation now. Let $\Gamma$ 
be a polycyclic-by-finite group, possibly with $\tau_\Gamma$ non-trivial.
Building on an idea of  \cite[Chapter 6, exercise 10 ]{Segal} for nilpotent
groups, we show in Section \ref{generalpofi}:
\begin{proposition}\label{equivi}  
Let $\Gamma$ be polycyclic-by-finite group and 
let $\tau_\Gamma$ denote its maximal finite normal subgroup. 
Then the groups $\Aut(\Gamma)$ and $\Aut(\Gamma /\tau_\Gamma)$ 
are commensurable. 
If $\Aut(\Gamma/\tau_\Gamma)$ is an arithmetic group 
then $\Aut(\Gamma)$ is arithmetic too. If $\Aut(\Gamma)$
is an arithmetic group then 
$\Aut(\Gamma/\tau_\Gamma)$ contains a subgroup of finite index which 
is an arithmetic group. 
\end{proposition}
Recall that two abstract groups $G$ and $G'$ are called \emph{commensurable}
if they contain finite index subgroups $G_0 \leq G$ and $G'_0 \leq G'$  
which are isomorphic. The examples in \cite{GP3} show that a 
group can contain an arithmetic subgroup of
finite index without being an arithmetic group.
We do not know whether this
phenomenon can arise in the situation of Proposition \ref{equivi}. 
As proved in \cite{GP3}, all arithmetic subgroups in algebraic groups
which do not have a quotient isomorphic to ${\rm PSL}(2)$ have only arithmetic
groups as finite extensions.

The groups $\Out(\Gamma)$ and $\Out(\Gamma/\tau_\Gamma)$ satisfy a weaker
equivalence relation (they are called \emph{S-commensurable}, see Section \ref{applih}). 
\begin{proposition}\label{finiin}
The natural homomorphism 
$\Out(\Gamma) \ra \Out(\tGam)$ has finite kernel and
maps $\Out(\Gamma)$ onto a finite index subgroup 
of\/  $\Out(\Gamma/\tau_\Gamma)$.
\end{proposition}

Clearly, $\Gamma/\tau_\Gamma$ is a wfn-group. From Theorem 
\ref{ftop} we already know that $\Out(\Gamma/\tau_\Gamma)$ is an arithmetic
group. We further know from \cite{Wehrfritz} 
that $\Out(\Gamma)$ is isomorphic to a subgroup of $\GL(n,\BZ)$, for some $n\in
\BN$, and hence is residually finite. Recall that a group $G$ is called {\it
  residually finite} if for every $g\in G$ with $g\ne 1$ there is a subgroup
of finite index in $G$ not containing $g$.
We complete the proof of Theorem \ref{teoa} in the general case by the following
result.
\begin{proposition}\label{finker} 
Let $A$ be a residually finite group and $E$ a finite normal subgroup in $A$.
If $A/E$ is arithmetic then $A$ is an arithmetic group. 
\end{proposition}
In \cite{Deligne} Deligne constructs lattices (in non-linear) Lie groups which
map onto arithmetic groups with finite kernel but which are not residually finite. 
The above gives a strong converse to this result.

To contrast the examples from Theorem \ref{teob} we formulate 
now a simple condition on a polycyclic-by-finite group which
implies that $\Aut(\Gamma)$ is an arithmetic group. We write 
$\lambda: \Gamma\to \tilde\Gamma=\Gamma/\tau_\Gamma$ for the projection
homomorphism. We define 
\begin{equation}
\tilde F(\Gamma):= \lambda^{-1}(\Fitt(\tilde\Gamma)).
\end{equation}
From general theory we infer that $\tilde F(\Gamma)$ is normal in $\Gamma$ and
that $\Gamma / \tilde F(\Gamma)$ is an abelian-by-finite group. Suppose that
$\Sigma \le \Gamma $ is a normal subgroup of finite index with 
$\tilde F(\Gamma)\le \Sigma$ and the additional property that 
$\Sigma /\tilde F(\Gamma)$ is abelian. Then the finite group 
$\mu(\Gamma,\Sigma)=\Gamma/\Sigma$ acts through conjugation on the abelian
group $\Sigma /\tilde F(\Gamma)$. We prove:

\begin{theorem}\label{teosep} 
Let $\Gamma$ be a polycyclic-by-finite group. 
Assume that there is a 
normal subgroup  $\Sigma\le\Gamma$ which is of finite index, and such that 
$\Sigma /\tilde F(\Gamma)$ is abelian and $\mu(\Gamma,\Sigma)=\Gamma/\Sigma$ 
acts trivially on $\Sigma /\tilde F(\Gamma)$.
Then
$\Aut(\Gamma)$ is an arithmetic group. 
\end{theorem}

Theorem \ref{teosep} has the following immediate corollaries.
\begin{corollary}\label{coronil} 
Let $\Gamma$ be a polycyclic-by-finite group.
If\/ $\Gamma\big/\Fitt(\Gamma)$ is nilpotent then $\Aut(\Gamma)$
is an arithmetic group.
\end{corollary}\label{coronilfi}

\begin{corollary} The automorphism group of a finitely generated and 
nilpotent-by-finite group is an arithmetic group. 
\end{corollary}

\begin{corollary}
Every polycylic-by-finite group has a finite index 
subgroup whose automorphism group is arithmetic.
\end{corollary}

Analogues of Theorem \ref{teosep} and Corollary \ref{coronil} have been proved in 
\cite{GP1} for arithmetic polycyclic-by-finite groups. The book \cite{Segal} contains in 
Chapter 8 various results which also show the
arithmeticity of $\Aut(\Gamma)$, but the hypotheses on $\Gamma$ are a lot
stronger than ours.  


\subsection{Cohomological representations}
\label{topram}


Let $\Gamma$ be a torsion-free 
polycyclic-by-finite group, and $R$ a commutative ring.  Let 
$\H^*(\Gamma,R)$ denote the cohomology
of $\Gamma$ with (trivial) $R$-coefficients. 
Since inner automorphisms
act trivially on cohomology, the outer automorphism group 
$\Out(\Gamma)$ is naturally represented on the cohomology 
space $\H^*(\Gamma,R)$. 
The finite dimensional complex vector 
space $\H^*(\Gamma,\BC)$ comes with an $\BZ$-structure
(and a resulting $\BQ$-structure)
given by the image of the base change homomorphism 
$\H^*(\Gamma,\bbZ) \ra  \H^*(\Gamma, \bbC)$.
This allows us to identify $\GLHC$ with a $\BQ$-defined
linear algebraic group.
The representation 
\begin{equation}
\rho: \Out(\Gamma) \to \GLHC
\end{equation}
is integral with respect to the above $\BZ$-structure.  
Moreover, we define in Section \ref{geometrc3}, 
building on geometric methods developed in \cite{Baues2}, a 
$\BQ$-defined homomorphism
\begin{equation}
\eta: \Out_a({\bf H}_\Gamma) \to \GLHC
\end{equation}
which extends the homomorphism $\rho$ via our homomorphism $\pi_{\cO_\Gamma}$
from Proposition \ref{ftop}. 
By composition, we obtain a $\BQ$-defined homomorphism
\begin{equation}
\rho_{\cO_\Gamma}:=\eta\circ\pi_{\cO_\Gamma}  :\  {\cal O}_\Gamma\to \GLHC. 
\end{equation}
We collect all of this together in the following theorem: 

\begin{theorem}\label{topolo1} 
Let $\Gamma$ be a torsion-free polycyclic-by-finite group.
Then there is a $\BQ$-defined homomorphism
$\rho_{\cO_\Gamma}: {\cal O}_\Gamma\to \GLHC$ such that the diagram
\begin{equation}
\xymatrix{
\Out(\Gamma) \ar[r]  \ar[d]_{\rho} 
& {\cal O}_\Gamma \ar[ld]^{\rho_{\cO_\Gamma}} \\
\GLHC&   
}
\end{equation}
is commutative.
\end{theorem}
Thus, the Zariski-closure of the image 
of $\Out(\Gamma)$ in $ \GLHC $
is a $\bbQ$-closed  
subgroup and we have: 
\begin{theorem}\label{topolo2} 
Let $\Gamma$ be a torsion-free polycyclic-by-finite group. Then
the cohomology image $\rho\big(\Out(\Gamma)\big)$
is an arithmetic subgroup of its Zariski-closure in 
$\GLHC$. The kernel of $\rho$ is a finitely generated
subgroup of $\Out(\Gamma)$.
If\/  $\Gamma$ is in addition nilpotent then the
kernel of $\rho$ is nilpotent by-finite.
\end{theorem}


\subsection{Applications to the groups of homotopy 
self-equi\-va\-lences of spaces}\label{applih}


Let $\Gamma$ be a group. The motivation to study 
the structure of $\Out(\Gamma)$
comes partially from topology since $\Out(\Gamma)$, for example, 
is naturally isomorphic to the group of homotopy classes of homotopy
self-equivalences of any $\KGone$-space (see \cite{Rutt}). 
A substantial theory 
using this tie between algebra and topology has been
developed starting from the work of Sullivan \cite{Sullivan}.  
We shall explain here the 
connections of our work with Sullivan's theory 
as described in \cite{Sullivan} and we also mention the additions 
we can make to Sullivan's theory. 

Sullivan considers spaces $X$ which have a nilpotent homotopy system
of finite type, that is, the homotopy groups 
$$(\pi_1(X);\pi_2(X),\pi_3(X),\ldots)$$
are all finitely generated, $\pi_1(X)$ is torsion-free nilpotent 
and acts nilpotently on
all $\pi_i(X)$ ($i\ge 2$). 
Sullivan associates to such a space $X$ 
the minimal model of the $\BQ$-polynomial forms of some complex representing
$X$. This is a finitely generated nilpotent differential algebra defined over
$\BQ$, it is called ${\cal X}$.
Under the further assumption that $X$ is either a finite complex or that the
homotopy groups $\pi_i(X)$ are trivial, for almost all $i\in \BN$,
he uses this construction to prove 
important results on the group $\Aut(X)$ of classes of
homotopy self-equivalences of the space $X$, cf.\ \cite[Theorem 10.3]{Sullivan}. 
 
Let us specialize, for a moment, to the case where $X$ is a $\KGone$-space,
 $\Gamma$ a torsion-free finitely generated nilpotent group. The results of
Sullivan can then be related to what we prove. 
The following table contains a dictionary
between the objects defined by Sullivan and those appearing in our theory.
\begin{center}
\begin{tabular}{|c|c|}
\hline
Sullivan  & algebraic theory \\
\hline
$\Aut(X)$ & $\Out(\Gamma)$ \\
${\cal X}$ & ${\lie g}$ \\
$\Aut({\cal X})$ & $\Aut({\lie g})$ \\
$\Aut_\BQ(X)=\Aut({\cal X})/{\rm inner\ automorphisms}$ 
& ${\cal O}_\Gamma \xrightarrow{\pi_{\cO_\Gamma}}  \Aut({\lie g})\big/\Inn_{\lie g}$ \\
\hline
\end{tabular}
\end{center}
Here ${\lie g}$ is the Lie algebra of the Malcev-completion of 
$\pi_1(X)=\Gamma$. The associated DGA $\cal X$  
coincides with the Koszul-complex of $\lie{g}$.  The $\BQ$-defined linear
algebraic group ${\cal O}_\Gamma$ is defined in the previous subsection. It
admits a $\BQ$-homomorphism $\pi_{\cO_\Gamma}$ with finite kernel onto 
the group $\Aut({\lie g})/\Inn_{\lie g}$ 
of Lie algebra automorphisms of ${\lie g}$ modulo inner automorphisms.
  
We mention the concept of 
{\it S-commensurability} of groups appearing in \cite{Sullivan}. 
This is the equivalence relation amongst groups which is generated by the 
operations of taking quotients with finite kernel and finite index subgroups. 

We shall now go through the four statements about $\Aut(X)$ which
are contained in \cite[Theorem 10.3]{Sullivan}, specialized to a 
$\KGone$ space $X$, and rephrase them in our context:

\medskip 
In Theorem 10.3 (i), it is stated that
{\it $\Aut(X)$ is S-commensurable with a full arithmetic subgroup of 
$\Aut_\BQ(X)$.} We prove (see Theorem \ref{ftop}), quite equivalently, that there is
a  homomorphism 
$\pi_{\Gamma} :\Out(\Gamma) \to  \Aut({\lie g})\big/\Inn_{\lie g}$ which has finite 
kernel and an arithmetic subgroup as image. We additionally  show that $\Out(\Gamma)$ is isomorphic
to an arithmetic subgroup in a $\BQ$-defined linear algebraic group ${\cal O}_\Gamma$ 
which is a finite extension group of $ \Aut({\lie g})\big/\Inn_{\lie g}$. 
(Examples (\cite{Deligne}, \cite{GP3}) show that the class of arithmetic
groups is not closed under any of the two operations used to define
S-commensurability. Hence
the information on $\Aut(X)=\Out(\Gamma)$ 
contained in  \cite[Theorem 10.3 i)]{Sullivan},  for our restricted case,  is
weaker then what follows from our results. 
Most dramatically there are groups which can be mapped onto an arithmetic
group with a finite kernel which are not even residually finite. As Serre
\cite{Serre} remarks, localization results in \cite{Sullivan2} imply that
the group $\Aut(X)$ is residually finite.  We wonder, whether the general theory 
can be extended similarly to our results.)

\medskip 
In Theorem 10.3 (ii) it is proved that
{\it the natural action of $\Aut(X)$ on the integral homology is compatible
  with an algebraic matrix representation of $\Aut_\BQ(X)$ on the vector
  spaces of rational homology.} We consider here cohomology instead of 
  homology and we obtain a similar result. We prove that the natural action 
  of $\Out(\Gamma)=\Aut(X)$ on the image of integral cohomology is compatible with an algebraic matrix 
representation of ${\cal O}_\Gamma$ on the vector
spaces of rational cohomology (see Theorem \ref{ftop} and Section
\ref{topram}). This has the important consequence that the stabilizers
of cohomology classes are arithmetic groups. 

\medskip 
By Theorem 10.3 (iii),
{\it the reductive part of $\Aut_\BQ(X)$ is faithfully represented on the
  natural subspace of homology generated by maps of spheres into $X$.}
We recall  (see Proposition \ref{cohom1}) that the reductive part of 
 $ \Aut({\lie g})\big/\Inn_{\lie g}$ is faithfully represented on the cohomology vector spaces.
Of course, Sullivan's result implies that the representation is already
faithful on the first homology space. We obtain a similar result for the
action of ${\cal O}_\Gamma$ on the first cohomology space. 

\medskip 
In Theorem 10.3 (iv) it is proved that
{\it as we vary $X$ through finite complexes, $\Aut(X)$ runs through every
commensurability class of groups containing arithmetic groups.}
Here our results give something different and new since
Sullivan uses simply connected spaces ($\pi_1(X)=\{1\}$)
which have non-vanishing higher homotopy groups to show his existence result:

Let \emph{N-commensurability} be the equivalence relation for groups which 
is generated by the 
operations of taking quotients 
with finitely generated nilpotent-by-finite groups as kernel 
and by taking finite index subgroups. We can prove: 
\begin{proposition}\label{exinili}
As we vary $X$ through all finite $\KGone$ complexes, $\Gamma$ a
finitely generated nilpotent group, $\Aut(X)$ runs through every
N-commensur\-ability class of groups containing an arithmetic group.
\end{proposition}
The result just mentioned follows
from our theory of cohomological
representations and from the existence theorems in \cite{BS}. We do not
provide details here. 
\medskip
 
Having discussed the relation of our results to Sullivan's we can mention
the following generalization for spaces $X$ 
which are $\KGone$ spaces for torsion-free
polycyclic-by-finite groups. In fact,  we can collect the results described in
Sections \ref{outline}, \ref{topram} to show: 

\begin{theorem}\label{teotop}
Let $\Gamma$ be a torsion-free polycyclic-by-finite group and let $X$ a
$\KGone$ space. Let $\Aut(X)$ be the group of classes of homotopy
self-equivalences of $X$ then the following hold: 
\begin{itemize}

\item[i)] $\Aut(X)$ is an arithmetic group.

\item[ii)] The action of $\Aut(X)$ on integral cohomology is 
compatible with an algebraic matrix 
representation of the $\BQ$-defined linear algebraic ${\cal O}_\Gamma$ on 
the vector spaces of cohomology.

\item[iii)] The stabilizer of an integral (rational) cohomology class is
arithmetic.  
\end{itemize} 
\end{theorem}
Again, we wonder whether our Theorem \ref{teotop} has an extension which does not
assume the vanishing of the higher homotopy groups.

As a corollary,  using well known finiteness results for arithmetic groups we
get:
\begin{corollary}\label{corotop}
Let $\Gamma$ be a torsion-free polycyclic-by-finite group and let $X$ a
$\KGone$-space. Let $\Aut(X)$ be the group of classes of homotopy
self-equivalences of $X$. Then the following hold:
\begin{itemize}
\item[i)] The group $\Aut(X)$ is residually finite.
\item[ii)] The group $\Aut(X)$ is finitely presented.
\item[iii)] The group $\Aut(X)$ contains only finitely many conjugacy classes of
finite subgroups.
\item[iv)]  Let $\mu$ be a finite group acting by group automorphisms on 
$\Aut(X)$, then the cohomology set ${\rm H}^1(\mu,\Aut(X))$ is finite.
\item[v)] The group $\Aut(X)$ is of type $WF$.
\end{itemize} 
\end{corollary}
For the definition and discussion of cohomology sets, see Section
\ref{grouplem}. The properties i),...,v) all follow from the fact that 
$\Aut(X)$ is an arithmetic group. Property i) is well known to hold for
finitely generated linear groups, 
for ii), iii) see \cite{Borel1}, iv) is proved in \cite{GP3} (for a weaker
form see \cite{BS}). 
The last item, property $WF$, means that every torsion-free subgroup $\Delta$ of 
finite index in $\Aut(\Gamma)$ appears as the fundamental group of a finite $K(\Delta,1)$.
See \cite{Serre} for discussion of this property of arithmetic groups and further references.
It is not known whether groups which are finite extensions  of  a
$WF$-group inherit the property $WF$, see \cite{Serre2, Brown}.
Note also that properties ii), iii) are compatible with
S-commensurability whereas i), iv) are not.

Let $X$ be a $\KGone$-space for a group $\Gamma$. 
The automorphism group of
$\Gamma$ can be naturally identified with the group of classes of
pointed homotopy
self-equivalences $\Aut^*(X)$ (see \cite{Rutt}). In case $\Gamma$ is
finitely generated torsion-free nilpotent it was well known that    
$\Aut^*(X)$ is an arithmetic group. Our results from Theorem \ref{teosep}
extend this to a much larger class of groups. 
Our examples constructed in Section \ref{exnona}, Theorem
\ref{teob} give rise to $\KGone$-spaces $X$, where $\pi_1(X)$ is nice and
 arithmetic whereas
$\Aut^*(X)$ is far from being an arithmetic group.

\medskip
{\em Acknowledgement } \ The authors wish to thank Wilhem Singhof
for helpful conversations concerning the contents of this article.


\section{Prerequisites on linear algebraic groups and arithmetic groups}
\label{prereq}


In the sequel we shall use a certain amount of the theory of 
linear algebraic groups and also of the theory of arithmetic groups. We 
briefly review here what we need and also add certain consequences of the
general theory. 
Our basic references are \cite{Borel}, \cite{BT}, \cite{PR}. 


\subsection{The general theory}


We use the usual terminology of Zariski-topology. Thus a 
linear algebraic group $\cA$ is a Zariski-closed subgroup of 
$\GL(n,\BC)$, for some $n\in \BN$. It is $\BQ$-defined if it is $\BQ$-closed.
We use the shorter term {\it $\BQ$-closed}
for
what should be called closed in the Zariski-topology with closed subsets
being those affine algebraic sets defined by polynomials with coefficients in
$\BQ$. If $R$ is a subring of $\bbC$, we put $\cA(R) = \cA \cap \GL(n,R)$.
We denote the connected component of the identity in 
a linear algebraic group $\cA$ by $\cA^\circ$. 
The connected component $\cA^\circ$ always has finite index in $\cA$ and 
$\cA$ is called \emph{connected} if $\cA=\cA^\circ$. If $\cA$ is defined
over $\bbQ$, its group of $\bbQ$-points $\cA(\bbQ)$ is Zariski-dense in $\cA$.

A {\it homomorphism of algebraic groups} 
is a morphism of the underlying affine algebraic varieties which is also
a group homomorphism. A homomorphism is {\it defined over $\BQ$} 
(or {\it $\BQ$-defined})
if the corresponding morphism of algebraic varieties is defined over $\BQ$.
It is called a {\it $\BQ$-defined isomorphism} if its inverse
exists and is also a homomorphism defined over $\BQ$.
An {\it automorphism} is an isomorphism of a linear algebraic group to itself.
We also use the more abstract concept of a $\BQ$-defined
linear algebraic group. 
It is well known that a $\BQ$-defined affine variety equipped with a
group structure given by $\BQ$-defined morphisms is $\BQ$-isomorphic with a 
$\BQ$-closed subgroup of 
$\GL(n,\BC)$, for some $n\in \BN$. As usually we write $\Ga=\BC$ 
for the additive
and $\Gm=\BC^*$ for the multiplicative group.
Quotients exist in the category of $\BQ$-defined linear algebraic groups. That
is given a $\BQ$-closed normal subgroup ${\cal N}$ of a 
$\BQ$-defined linear algebraic group $\cA$ the quotient group $\cA/{\cal N}$
is a $\BQ$-defined linear algebraic group and the natural map 
$\cA\to\cA/{\cal N}$ is a $\BQ$-defined homomorphism. We say that 
$\cA$ is the {\it almost direct product} of two Zariski-closed subgroups
${\cal B},\,{\cal C}$ if ${\cal B}$ and ${\cal C}$ centralize each other, 
their intersection ${\cal B}\cap{\cal C}$ is finite and if  
$\cA={\cal B}\cdot{\cal C}$ holds. 

Let now $\cA$ be a $\BQ$-defined linear algebraic group, we write 
$\bU_{\cA}$ for its {\it unipotent radical}. This is the largest Zariski-closed
unipotent normal subgroup in $\cA$, it is $\BQ$-closed. 
The algebraic group $\cA$ is called {\it reductive} if $\bU_{\cA}=\{1\}$.
In particular, the quotient ${\cA}^{\rm red}  = \cA \big/ \bU_{\cA}$ is reductive,
for all linear algebraic groups $\cA$.
We write 
${\cA}^{\rm sol}$ for the solvable radical of $\cA$. This is the largest 
connected Zariski-closed
solvable normal subgroup in $\cA$, it is $\BQ$-closed. 
A connected linear algebraic group $\cA$ is called \emph{semisimple} 
if ${\cA}^{\rm sol}  =\{1\}$. 

A linear algebraic group is called a {\it d-group} if
it consists of semisimple elements only. 
A Zariski-closed subgroup which is a
d-group is called a {\it d-subgroup}.
A d-group is reductive and
abelian-by-finite. A {\it torus} is a linear algebraic group isomorphic to 
$\Gm^n$,  for some $n\in\BN$. If ${\cal S}$ is a d-group then 
${\cal S}^\circ$ is a torus.
 
Let $\cA$ be a $\BQ$-defined linear algebraic group. Then the following hold
(see \cite{Borel,BS,BT, PR}).
\begin{itemize}
\item[{\bf AG1:}] 
There exists a 
$\BQ$-closed reductive complement $\cR \leq \cA$, that is, 
$\cA=\bU_{\cA}\cdot \cR$ is a direct product of subgroups,
and $\cR  \cong \cA^{\rm red}$. In particular,  we also have 
$\cA \cong \bU_{\cA}\rfish \cA^{\rm red}$. Here the symbol $\cong$ indicates $\bbQ$-isomorphism
of linear algebraic groups,  and $\rfish$ indicates a semi-direct product.
\item[{\bf AG2:}]  All $\BQ$-closed reductive complements are conjugate by
elements of $\bU_{\cA}(\BQ)$.
\item[{\bf AG3:}] Let $\cA$ be a $\BQ$-defined connected and reductive 
linear algebraic group.  Then $\cA$ is the almost direct product of its 
center, which is a $\BQ$-closed d-subgroup, 
and the commutator subgroup $[\cA,\cA]$, which is
$\BQ$-closed, connected and semisimple. 
\item[{\bf AG4:}] Let $\cA$ be a $\BQ$-defined connected and semisimple 
linear algebraic group. Then its group $\Aut_a(\cA)$ of algebraic
automorphisms is  a $\BQ$-defined linear algebraic group and $\Inn_\cA$ is its
connected component.
\item[{\bf AG5:}] Let $\cA$ be a $\BQ$-defined connected and reductive 
linear algebraic group. Let ${\cal B}$ be $\BQ$-closed normal subgroup  of
$\cA$. Then there is a $\BQ$-closed  subgroup ${\cal C}$ of
$\cA$ which centralizes ${\cal B}$, satisfies $\cA={\cal B}\cdot {\cal C}$
such that ${\cal B}\cap {\cal C}$ is finite. That is, $\cA$ is the almost
direct product of ${\cal B}$ and ${\cal C}$. 
\item[{\bf AG6:}] Let ${\cal S}$ be a commutative d-group and $G$ 
a finite group of $\BQ$-defined
automorphisms of ${\cal S}$. Let ${\cal S}_1\le {\cal S}$ be a $G$-invariant
$\BQ$-closed subgroup. Then there is a $G$-invariant $\BQ$-closed subgroup 
 ${\cal S}_2\le {\cal S}$ such that ${\cal S}$ is the almost direct product of
 ${\cal S}_1$ and ${\cal S}_2$.
 \item[{\bf AG7:}]  Let $\cT \leq \cA$ be a torus. Then the centralizer
 $\Z_\cA(\cT)$ of $\cT$ in $\cA$ has finite index in the normalizer $\N_\cA(\cT)$.
\end{itemize}

Statement {AG6} is proved by using 
the category equivalence between the category of
$\BQ$-defined commutative 
d-groups and the category of continuous, $\BZ$-finitely
generated modules for the absolute Galois group of $\BQ$ 
(see \cite[\S 8]{Borel}) together with Maschke's theorem. 
The last fact AG7 is called the \emph{rigidity of tori}, see
\cite[Corollary 2 of III.8]{Borel} for a proof.

We shall add some remarks concerning the structure of 
solvable-by-finite groups. 
A linear algebraic group $\bH$ is solvable-by-finite if
its identity component $\bH^\circ$ is solvable.
A Cartan subgroup $\bC$ of  
$\bH$ is by definition 
the normalizer of a maximal torus $\bT$ in $\bH$.

For a $\BQ$-defined, solvable-by-finite linear algebraic group $\bH$ the
following hold:

\begin{itemize}
\item[{\bf SG1:}] There are maximal tori $\bT\le \bH$ which are $\BQ$-closed.
\item[{\bf SG2:}] Let $\bT$ be a maximal torus of 
$\bH$ and $\bC=\N_{\bH}(\bT)$. Then 
$\bC^\circ =  \N_{\bH^\circ}(\bT)$ 
equals the centralizer of $\bT$ in
$\bH^\circ$ and is a connected nilpotent group.    
Moreover, $\bC$ contains a maximal d-subgroup $\bS$ 
with $\bS^\circ = \bT$. 
\item[{\bf SG3:}] Let $\bT\le \bH$ be a maximal torus 
which is $\BQ$-closed. Then 
$\bC:=\N_{\bH}(\bT)$ is $\BQ$-closed and contains a maximal d-subgroup 
$\bS$ of
$\bH$ which is $\BQ$-closed and which satisfies $\bS^\circ = \bT$
\item[{\bf SG4:}] Let $\bS$ be a maximal d-subgroup of $\bH$. Then 
$\bH={\bf U}_{\bH}\cdot \bS={\bf U}_{\bH}\rfish \bS$.
\item[{\bf SG5:}] All maximal 
tori and also all maximal d-subgroups in $\bH$ are conjugate by elements 
of $[\bH^\circ, \bH^\circ] \leq \bU_{\bH}$. 
\item[{\bf SG6:}] All $\BQ$-closed maximal tori 
and all $\BQ$-closed maximal d-subgroups 
are conjugate by elements 
of $[\bH^\circ, \bH^\circ](\bbQ)$.
\item[{\bf SG7:}]
Let $\bC$ be a Cartan subgroup of $\bH$ and $\bF\le \bU_\bH$ a
normal subgroup of $\bH$ which contains 
$[\bH^\circ, \bH^\circ]$ then we have $\bH = \bF\cdot \bC$. 
\end{itemize}

For almost all of this see \cite[Chapter III]{Borel}. 
Statement SG6 is contained  in \cite{BS}, for a detailed 
proof see \cite[Section 2]{GP1}.

We have defined the notion of an arithmetic subgroup of a $\BQ$-defined linear
algebraic group and that of an arithmetic group in the beginning of the
introduction. We shall often have to use the behaviour of arithmetic subgroups
under $\BQ$-defined homomorphisms. To describe this, let 
$\rho :\cA_1\to\cA_2$ be a $\BQ$-defined homomorphism between 
$\BQ$-defined linear algebraic groups $\cA_1,\, \cA_2$ and $A\le \cA_1$ a
subgroup. Furthermore we suppose that $\rho$ is surjective. We have:  
\begin{itemize}
\item[{\bf AR1}:] If $A$ is an arithmetic subgroup of $\cA_1$ then     
$\rho(A)$ is an arithmetic subgroup of $\cA_2$.
\item[{\bf AR2}:] Suppose that $\rho$ is injective, then $A$ 
is an arithmetic subgroup of $\cA_1$ if and only if $\rho(A)$ is commensurable
with $\cA_2(\BZ)$.
\item[{\bf AR3}:] Every abelian subgroup of an arithmetic group is finitely
generated. 
\item[{\bf AR4}:] Let $\cA$ be a $\bbQ$-defined 
group and  $\cA_1$, $\cA_2$ $\bbQ$-defined subgroups
such that $\cA = \cA_1 \rtimes \cA_2$ is a semi-direct product. 
Let $A  \leq \cA(\BQ)$ be a subgroup. Assume there exists subgroups 
$A_1 $, $A_2$ of $A$ such that $A = A_1 \rtimes \A_2$ is a semi-direct
product. If  $A_1$ and $A_2$ are arithmetic subgroups in $\cA_1$ and $\cA_2$, respectively, then $A$ is an arithmetic subgroup of $\cA$. 
\end{itemize}
Statement AR1 is proved in \cite{Borel2} and AR2 is a consequence of it. 
Statement AR2 shows that the notion of an arithmetic subgroup does not
depend on the particular embedding of the ambient algebraic group into
$\GL(n,\BC)$. Statement AR3 follows by consideration 
of the Zariski-closure of
the given abelian subgroup.  Statement AR4 is a consequence of AR1.


\subsection{Algebraic groups of automorphisms}


Let $G$ be a group acting by automorphisms on a group $A$ and let $B\le A$ be
a subgroup of $A$. We write
$$\N_G(B):=\{\, g\in G\ |\ g(B)=B\,\},\quad  
{\rm Z}_G(B):=\{\, g\in G\ |\ g(b)=b\ {\rm  for\ all}\ b\in B\,\}$$
for the normalizer and centralizer of $B$ in $G$.

\begin{definition}\label{algauto} 
Let ${\cal A}$ be a $\BQ$-defined linear algebraic group 
and let $G$ be a group which acts by $\BQ$-defined automorphisms on $\cA$. 
Then $G$ normalizes the unipotent radical $\bU_\cA$ and we
say that $G$ acts as an {\it algebraic group of automorphisms} on $\cA$ if 
${\rm Z}_G({\cal T})$ has finite index in $\N_G({\cal T})$ for every
$\bbQ$-defined torus  ${\cal T}$ of ${\cA}^{\rm red} = \cA\big/\bU_\cA$.   
\end{definition}

Note that $G=\GL(n,\BZ)$ acts by $\BQ$-defined automorphisms on the torus
$\cA=\Gm^n$. But for $n\ge 2$ this is not an algebraic group of automorphisms
on $\cA$. Positive examples are given in the next lemma. 

\begin{lemma}\label{algautoex}
Let ${\cal B}$ be a $\BQ$-defined linear algebraic group and 
$\cA\le{\cal B}$ a 
$\BQ$-closed subgroup. Let further $G\le {\cal B}(\BQ)$ be a subgroup which
normalizes $\cA$. Then $G$ acts by conjugation on $\cA$ as an algebraic group
of automorphisms.
\end{lemma}
\begin{proof}
Replacing ${\cal B}$ by the Zariski-closure of $\cA \cdot G$, 
we may assume that $\cA$ is
normal in ${\cal B}$. Then $\bU_\cA$ is also normal in ${\cal B}$ and the image
of $G$ in ${\cal B}\big/\bU_\cA$ acts by conjugation on $\cA\big/\bU_\cA$. 
Thus  the lemma follows from the rigidity of tori  (AG7). 
\end{proof}

\begin{lemma}\label{algformal1} 
Let ${\cal A}$ be a $\BQ$-defined linear algebraic group. 
Let $G$ be an algebraic  group of automorphisms on $\cA$.
Then 
${\rm Z}_G({\cal S})$ has finite index in $\N_G({\cal S})$ for every
$\bbQ$-closed commutative d-subgroup $\cS$ of $\cA\big/\bU_\cA$. 
\end{lemma}
\begin{proof}
To prove the lemma, we may replace $\cA$ by $\cA\big/\bU_\cA$. 
We consider a $\BQ$-closed commutative
d-subgroup ${\cal S}$ of $\cA$. 
Its connected component ${\cal S}^\circ$ is a
$\BQ$-defined torus in $\cA^\circ$.  
We have 
$\N_G({\cal S})\le \N_G({\cal S}^\circ)$.
Our hypothesis implies that
$Z_1=\N_G({\cal S})\cap{\rm Z}_G({\cal S}^\circ)$ has finite index in 
$\N_G({\cal S})$. 

There is a finite Zariski-closed
subgroup ${\cal E}\le{\cal S}$ such that ${\cal S}={\cal S}^\circ\cdot  
{\cal E}$. Let  ${\cal E}_0$ be the kernel of the homomorphism from ${\cal S}$ to  
${\cal S}$ which sends $x\in {\cal S}$ to $x^{|{\cal E}|}$.
It is finite, contains ${\cal E}$ and is normalized by 
$\N_G({\cal S})$. 
Trivially $Z_2=\N_G({\cal S})\cap{\rm Z}_G({\cal E}_0)$
has finite index in $\N_G({\cal S})$. Since $Z_1\cap Z_2$ is contained in 
${\rm Z}_G({\cal S})$, this  proves that ${\Z}_G(\cS)$ is of finite index
in ${\N}_G({\cal S})$.
\end{proof}

The following lemmata describe natural operations which can be perfomed with
algebraic groups of automorphisms.

\begin{lemma}\label{algformal2} 
Let ${\cal A}$ be a $\BQ$-defined linear algebraic group and ${\cal B}$ a
$\BQ$-closed normal subgroup in $\cA$. 
Let $G$ be an algebraic group of automorphisms of $\cA$ which normalizes
${\cal B}$.
Then $G$ acts as a group of algebraic automorphisms on ${\cal B}$ and also on 
$\cA\big/{\cal B}$. 
\end{lemma}
\begin{proof} This is clear for $\cB$, since $\cB$ is normal in $\cA$.
To prove that $G$ acts as an algebraic group of
automorphisms on $\cA\big/\cB$, we use a 
$G$-invariant decomposition AG5.  We leave the details to
the reader. 
\end{proof}

Remark that, in general,  the converse of the statements of Lemma \ref{algformal2} 
is not true. Nevertheless, we have:

\begin{lemma}\label{algformalproducts} 
Let ${\cal A}$ and $\cB$ be $\BQ$-defined linear algebraic groups, and let 
$G$ act as an algebraic group of automorphisms on $\cA$ and on $\cB$. Then
$G$ acts as an algebraic group of automorphisms on the product 
$\cC = \cA \times \cB$.  \end{lemma}
\begin{proof} Let $G$ act on $\cC$ by the product action. 
Let $\cT$ be a torus in $\cC^{\rm red} =  {\cA}^{\rm red} \times {\cB}^{\rm red} $ 
and put  $G' = N_G(\cT)$. The torus $\cT$ is contained in the direct product  
of factors  $\cT_\cA \leq {\cA}^{\rm red} $, $\cT_\cB \leq {\cB}^{\rm red} $ which 
are stabilized by $G'$. Hence,  both are centralized by a finite index subgroup of $G'$.
Therefore, $\cT$ is centralized by this finite index subgroup of $G'$.
\end{proof}

\begin{lemma}\label{algformal3} 
Let ${\cal A}$ be a $\BQ$-defined linear algebraic group and $G$ a group of
$\BQ$-defined automorphisms of $\cA$. Suppose $L,\, H\le G$ 
are subgroups such that $L$ is normal in $G$ and $G=L\cdot H$ holds. If
both $L$ and $H$ are $\cA$-algebraic groups of automorphisms of $\cA$ then also $G$
has this property.
\end{lemma}
\begin{proof} We consider $G$ as a group of
$\BQ$-defined automorphisms of\ ${\cal B}=(\cA\big/\bU_\cA)^\circ$.
We have the almost direct decomposition 
${\cal B}={\cal Z}\cdot [{\cal B},{\cal B}]$, where ${\cal Z}$ is the center of
${\cal B}$. The group $G$ normalizes both ${\cal Z}$ and 
$[{\cal B},{\cal B}]$. Let ${\cal T}$ be a $\bbQ$-defined torus
of ${\cal B}$. If ${\cal T}\le [{\cal B},{\cal B}]$ then 
${\rm Z}_G({\cal T})$ has finite index in   
${\rm N}_G({\cal T})$ since $G$ contains a subgroup of finite index which acts
by inner automorphisms on ${[\cB,\cB]}$ (AG4). 
Now assume $\cT \leq \cZ$.
By our assumption on $L$ and $H$ it follows that 
${\rm Z}_G({\cal Z}^0)$ has finite index in   
$G$  (since $G$ normalizes ${\cal Z}^0$), and, a fortiori,
${\rm Z}_G({\cal T})$ has finite index in   
${\rm N}_G({\cal T})$. 
The case of a general ${\cal T}$ can be treated 
by projecting $\cT$ on the factors of the almost direct
product decomposition AG3.
\end{proof}

In Section \ref{subsecari} the following fact plays a crucial role.

\begin{proposition}\label{decozen}
Let ${\cal A}$ be a connected\/  $\BQ$-defined linear algebraic group,  equipped 
with an algebraic group of automorphisms $G$.
Then there are: \begin{itemize}

\item[i)] a $\BQ$-closed, $G$-invariant central d-subgroup 
${\cal Z}_1$ in $\cA$ and 

\item[ii)] a $\BQ$-closed, $G$-invariant 
subgroup ${\cal A}_1$ of $\cA$ with unipotent-by-finite
center 
\end{itemize}
such that
$\cA$ is decomposed as the almost direct product of $\cZ_1$ and $\cA_1$. 
\end{proposition}
\begin{proof}
Let $\cA^{\rm red}$ be a maximal reductive subgroup of $\cA$ which
is defined over $\bbQ$. Let ${\cal Z}$ be the center of $\cA^{\rm red}$
and  let $\cZ_1$ be the maximal torus in the center of $\cA$. Thus $\cZ_1$
is $\bbQ$-closed in $\cA$, it is 
contained in $\cZ$, and it is normalized by $G$.

By Lemma \ref{algformal1}, the group $G_1$ of automorphisms of $\cZ$ 
which is induced by the quotient action of $G$ on the isomorphic image 
of $\cZ$ in $\cA\big/\bU_\cA$ is finite.  Hence, by AG6, we may choose a $\BQ$-closed and 
$G_1$-invariant subgroup ${\cal Z}_2\le {\cal Z}$ such that
${\cal Z}$ is the almost
direct product of ${\cal Z}_1$ and ${\cal Z}_2$.
This shows that  ${\cal Z}_1$ and $\cA_1 ={\bf U}_\cA\cdot{\cal Z}_2\cdot [\cA^{\rm red},\cA^{\rm red}]$ satisfy the requirements of
the lemma.
\end{proof}


\section{The group of automorphisms of a solvable-by-finite 
linear algebraic group} \label{section:auta}


Let $\cA$ be a linear algebraic group
defined over $\bbQ$.     
We let 
\begin{equation} \Aut_{a}(\cA) \; , \; \Aut_{a,\bbQ}(\cA)
\end{equation} 
denote 
the group of all automorphisms of $\cA$, respectively 
the group of all $\bbQ$-defined automorphisms of $\cA$. 
We also need the following concept. 
\begin{definition}
We say that a linear algebraic group $\cA$ has a {\it strong unipotent
  radical} if the centralizer $\C_{\cA}(\bU_\cA)$ of its unipotent radical
$\U_\cA$ is contained in $\U_\cA$. 
\end{definition}
Given a $\BQ$-defined solvable-by-finite linear algebraic group $\bH$  
with a strong unipotent radical we will identify $\Aut_{a}(\bH)$ in a natural 
way with a $\bbQ$-defined linear algebraic group. We obtain special features 
of this identification which we will need later. 

For a general linear algebraic group $\cA$,
the group $\Aut_{a}(\cA)$ is a group of type ALA, that is, 
$\Aut_{a}(\cA)$ is an extension 
of an affine algebraic group by 
an arithmetic group, see \cite{BS}. 
Our approach is a variation
of corresponding results in \cite{BS}. A construction similar to ours, but in a more
restricted situation, is contained in \cite{GP1}.


\subsection{The algebraic structure of $ \Aut_{a}(\bH)$} 
\label{subsect:algebraicstructure}


We will assume here that $\bH$ is a $\BQ$-defined solvable-by-finite
linear algebraic group which has a 
strong unipotent radical $\bU:=\bU_\bH$.
Let $\lu$ denote the Lie-algebra of $\bU$. 
The Lie-algebra $\lu$ is
defined over $\bbQ$ and $\Aut(\lu) \leq \GL(\lu)$
is a $\bbQ$-defined linear algebraic group. The exponential map
\begin{equation} 
\exp:  \lu \rightarrow  \bU
\end{equation} 
is a $\bbQ$-defined isomorphism of varieties. Thus, via the 
bijective homomorphism  
\begin{equation} \label{eq:exp}
  \Aut_a(\bU) \ni \, \Phi \; \mapsto \; \exp^{-1} \circ \,
\Phi \circ \exp  \; \, \in \Aut(\lu)\; 
\end{equation}
the group $\Aut_a(\bU)$ attains a natural structure of a linear 
algebraic group which is defined over $\bbQ$.  

We let $\bS\le \bH$ be a maximal $\BQ$-closed d-subgroup and obtain the
decomposition $\bH=\bU\cdot \bS$ described in SG4 of Section \ref{prereq}.
We define
\begin{equation}
\Aut_{a}(\bH)_\bS \, := \, 
\{\, \Phi \in  \Aut_{a}(\bH)  \mid \Phi(\bS) = \bS \} \; \leq 
\Aut_a(\bH) . 
\end{equation}

\begin{lemma}
The restriction map $\Phi \mapsto \Phi|_{\bU}$
identifies $ \Aut_{a}(\bH)_\bS$ with a 
$\bbQ$-defined closed subgroup of $\Aut_a(\bU)$.
\end{lemma} 
\begin{proof}
Since $\bH$ has a strong unipotent radical the restriction map is
injective on  $\Aut_{a}(\bH)_\bS$. 
Let $\Ad(\bS)$ denote the image of $\bS$ 
under the restriction of the adjoint representation to $\lu$.
Then  $\Ad(\bS)$ is a $\bbQ$-defined 
subgroup of $\Aut(\lu)$. Note that the obvious homomorphism
$\bS\to \Ad(\bS)$ is a $\BQ$-defined isomorphism since $\bH$ has a strong
unipotent radical (see also \cite{GP1}, Lemma 2.3).
It is clear that, via the 
exponential map, the restriction of $ \Aut_{a}(\bH)_\bS$ to $\bU$ 
is isomorphic to the normalizer 
of $\Ad(\bS)$ in $\Aut(\lu)$. Since $\Ad(\bS)$ is a $\bbQ$-defined
group, its normalizer and centralizer are Zariski-closed subgroups
of $\Aut(\lu)$. By the Galois criterion for 
rationality (see \cite[AG 14]{Borel})
the normalizer and centralizer of $\Ad(\bS)$ are defined over $\bbQ$.
In particular, $ \Aut_{a}(\bH)_\bS$ restricts to a $\bbQ$-defined subgroup of
$\Aut_a(\bU)$. 
\end{proof}

Thus the lemma furnishes a natural structure of $\bbQ$-defined 
linear algebraic group on $ \Aut_{a}(\bH)_\bS$. 
Since $ \Aut_{a}(\bH)_\bS$ acts on $\bU$ by
$\bbQ$-defined morphisms, 
the semi-direct product $\bU \rtimes \Aut_{a}(\bH)_\bS$ is
an affine algebraic group over $\bbQ$, and thus is
also equipped with the natural structure of 
a $\bbQ$-defined linear algebraic group, see \cite{Borel}.
Let $\Inn_u \in \Aut_a(\bH)$ denote the inner automorphism 
corresponding to $u\in \bU$, that is,  $\Inn_u(h) = u h u^{-1}$, for
all $h \in \bH$. We consider the  homomorphism 
\begin{equation} \label{eq:Theta}
 \Theta: \, \bU \rtimes   \Aut_{a}(\bH)_\bS  \longrightarrow  \Aut_{a}(\bH) 
\,\; , \; \;  (u, \Phi) \mapsto \Inn_u \circ \Phi \; . 
\end{equation}

\begin{lemma} The homomorphism $\Theta$ is surjective. 
\end{lemma}
\begin{proof}
The group $\bH$ 
decomposes as a semi-direct
product of $\bbQ$-defined algebraic groups
$ \bH = \bU \cdot \bS  \cong  \bU \rtimes \bS$.  
All maximal d-subgroups of $\bH$ are conjugate by elements of $\bU$, 
and $\bbQ$-defined
d-subgroups are conjugate by elements of $\bU(\bbQ)$ (see Section
\ref{prereq}). 
Therefore, the homomorphism $\Theta$ 
is onto $\Aut_a(\bH)$.
\end{proof}

\begin{lemma}
The kernel of the homomorphism $\Theta$ is a $\bbQ$-defined 
unipotent subgroup of $\bU \rtimes  \Aut_{a}(\bH)_\bS$.
\end{lemma}
\begin{proof}
In fact, $\ker \Theta = \{ (u, \Inn_u^{-1} ) \mid 
\Inn_u \in  \Aut_{a}(\bH)_\bS \}$ is a unipotent subgroup 
of  $\bU \rtimes  \Aut_{a}(\bH)_\bS$. Since the homomorphism 
$\Theta$ corresponds to
a $\bbQ$-defined action of  $\bU \rtimes  \Aut_{a}(\bH)_\bS$
on the variety $\bU$, its kernel is a $\bbQ$-defined
subgroup. 
\end{proof}
Let 
\begin{equation}
 \cA_\bH :=  
\bU \rtimes  \Aut_{a}(\bH)_\bS \big/ \, \ker \Theta
\end{equation}  
denote the quotient group 
of $ \bU \rtimes   \Aut_{a}(\bH)_\bS$ by the kernel of $\Theta$. 
The group $\cA_\bH$ has a natural structure of
a $\bbQ$-defined linear algebraic group, 
such that
\begin{equation}
 \cA_\bH(\bbQ) = \bU \rtimes  \Aut_{a}(\bH)_\bS  \left(\bbQ\right) \big/ 
\left(\ker \Theta\right)(\bbQ) \;. 
\end{equation}
Therefore, we have:

\begin{theorem}
The homomorphism 
$   \cA_\bH  \rightarrow  \Aut_{a}(\bH) $
which is induced by \eqref{eq:Theta} naturally identifies $\Aut_a(\bH)$ 
with the complex points 
of the $\bbQ$-defined linear algebraic group $\cA_\bH$.  
\end{theorem}
By arguments using the above setup, 
the rational points of the corresponding affine algebraic group 
are naturally interpreted in the following way: 

\begin{proposition} \label{rationalisq}
Under the homomorphism $   \cA_\bH  \rightarrow  \Aut_{a}(\bH) $
which is induced by \eqref{eq:Theta} 
the group of rational points $\cA_\bH(\bbQ)$ of $\cA_\bH$ 
corresponds to the group $\Aut_{a,\bbQ}(\bH)$ of\/
$\bbQ$-defined automorphisms of\/ $\bH$. 
\end{proposition}
\begin{proof} Let $\Phi \in \cA_\bH(\bbQ)$. It is clear from
our construction that $\Phi$ preserves the group of rational 
points $\bH(\bbQ) \leq \bH$ which is Zariski-dense in $\bH$
(cf. \cite[\S18.2]{Borel}). 
By the Galois-criterion for rationality (compare \cite[AG 14]{Borel}), 
$\Phi$ is defined over $\bbQ$.

Conversely, assume $\Phi \in \Aut_a(\bH)_\bbQ$ is a $\bbQ$-defined
automorphism of $\bH$. Then, by AG2, there exists $v \in \bU(\bbQ)$ such that 
$\Psi = \Inn_v \circ \Phi \in \Aut_a(\bH)_\bS$, and $\Psi$ is
defined over $\bbQ$ as well.  The exponential correspondence shows
that $\Psi \in \Aut_a(\bH)_\bS(\bbQ)$. It follows that 
$\Phi \in   \Aut_a(\bH)(\bbQ)$. 
\end{proof}
Henceforth, we will identify $\Aut_a(\bH)$ with (the complex points of)
the linear algebraic group $\cA_\bH$, and
$\Theta: \bU \rtimes   \Aut_{a}(\bH)_\bS \rightarrow \cA_\bH = \Aut_a(\bH)$
becomes a morphism of algebraic groups which is defined over $\bbQ$.


\subsection{Arithmetic subgroups of $\Aut_a(\bH)$} \label{section:asubgroups}


We keep the notation of the previous subsection. In particular, 
$\bH$ denotes a $\BQ$-defined solvable-by-finite linear algebraic group with 
a strong unipotent radical  $\bU=\bU_\bH$. 
Let further $\bS$ be a maximal $\BQ$-defined d-subgroup
of $\bH$.

Let $\theta \leq \bU({\bbQ})$ be a finitely generated 
subgroup which is Zariski-dense in $\bU$. Then  $\theta$ is an
arithmetic subgroup of the $\bbQ$-defined group $\bU$. 
We explain now how the choice of $\theta$ gives rise
to an arithmetic subgroup $A_\theta$ of  $\Aut_a(\bH)=\cA_\bH$.  

\begin{definition}\label{LThedef}
We  define the following subgroups of $\Aut_a(\bH)$: 
\begin{equation}
\cA_\theta :=  \{\, \Phi \in \Aut_a(\bH)  \mid \Phi(\theta) = 
\theta \, \} \; ,\qquad 
 \Lambda_\theta \, := \,  \Aut_{a}(\bH)_\bS  \cap  \cA_\theta \; ,  
\end{equation} 
\begin{equation}
 \A_\theta :=  \Theta\left(\theta \rtimes 
\Lambda_\theta\right) \,  \leq \cA_\bH{(\bbQ)} \; .
\end{equation}
\end{definition}

\begin{lemma}   \label{lemma:Lambdatheta}
Let $\theta \leq \bU({\bbQ})$ be a finitely generated 
subgroup which is Zariski-dense in $\bU$. Then 
the group $\Lambda_\theta$ is an arithmetic subgroup 
of $ \Aut_{a}(\bH)_\bS$. 
\end{lemma} 
\begin{proof} As explained in Subsection \ref{subsect:algebraicstructure}, 
the group $ \Aut_{a}(\bH)_\bS$ is $\bbQ$-isomorphic to the normalizer
$\N_{\Aut(\lu)}(\Ad(\bS))$ via the isomorphism \eqref{eq:exp}.
Embedding $\Aut(\theta)$ as a subgroup of  $\Aut_a(\bU)$, it
corresponds to 
$$ N_{\Aut(\lu)}\big( \log \theta\big) = 
\{ \Psi \in \Aut(\lu) \mid \Psi (\log 
\theta) \subseteq \log \theta \}\, \leq \Aut(\lu) \; . $$
The group $N_{\Aut(\lu)}( \log \theta)$
stabilizes the  lattice $L$ in $\lu$ 
which is generated by the set $\log \theta$. 
Moreover, $N_{\Aut(\lu)}( \log \theta)$ is arithmetic in $\Aut(\lu)$, 
see  \cite[Chapter 6]{Segal}. 
It follows that $\Lambda_\theta \leq  \Aut_{a}(\bH)_\bS$  
corresponds to 
$$ \N_{\Aut(\lu)}\big(\Ad(\bS)\big)  
\cap N_{\Aut(\lu)}\big( \log \theta(\Gamma)\big) \; . $$ 
Thus,
under the natural linear representation of $ \Aut_{a}(\bH)_\bS$ in
$\GL(\lu)$,  $\Lambda_\theta$ is commensurable to
the stabilizer of a lattice $L  \subset \lu({\bbQ})$.    
Therefore,  $\Lambda_\theta$ is arithmetic in $ \Aut_{a}(\bH)_\bS$.
\end{proof}

Using AR4 in Section \ref{prereq} we deduce that $\theta \rtimes \Lambda_\theta $
is arithmetic in $\bU \rtimes  \Aut_{a}(\bH)_\bS$.
Since $A_\theta$ is the image of the arithmetic group 
$\theta \rtimes \Lambda_\theta$ under the $\bbQ$-defined 
surjective homomorphism
$\Theta$, $A_\theta$ is arithmetic in $ \Aut_{a}(\bH)$. This proves:

\begin{proposition} \label{prop:Atheta}
Let $\theta \leq \bU({\bbQ})$ be a finitely generated 
subgroup which is Zariski-dense in $\bU$.  Then 
the group $\A_\theta$ is an arithmetic subgroup of $ \Aut_{a}(\bH)$. 
\end{proposition}


\subsection{Some closed subgroups of $\Aut_a(\bH)$}\label{section:subgroups}


We stick to the conventions about $\bH$, $\bU$, $\bS$.
In addition,  we introduce $\bF \leq \bU$ to be a 
$\bbQ$-closed normal subgroup of $\bH$ 
which contains the commutator group $[\bH^\circ, \bH^\circ]$. 
By SG7 of Section \ref{prereq},  we have $\bH=\bF\cdot\bC$, for every
$\BQ$-closed Cartan subgroup of $\bH$.
As before
$\bS$ denotes a maximal d-subgroup in $\bH$. 

Let $\N_\cA(\bF)$ denote  the subgroup of elements in 
$\Aut_a(\bH)$ which preserve $\bF$. Additionally, we define: 
$$   \cA_{\bH | \bF}   :=  \,    \{ \Phi \in \N_\cA(\bF)  \mid 
  \Phi|_{{\bH /\bF}} =\id_{\bH / \bF} \}  \; , $$
$$    \cAS   :=  \,  \{ \Phi \in \N_\cA(\bF)  \mid 
  \Phi|_{{\bH^\circ /\bF}} = \id_{\bH^\circ  / \bF} \, , \, \Phi(\bS) = \bS \} \; , $$
$$  \cASid   :=  \, \{ \Phi \in \cA_{\bS} 
\mid \Phi|_\bS  = \id_\bS \} \; . $$
It is easy to see that these groups are $\bbQ$-closed 
subgroups of $\Aut_a(\bH)$. 

\begin{lemma}  \label{lemma:subgroups1}
The following hold: \begin{itemize}

\item[i)] $\cASid \leq \cA_{\bH | \bF} $.

\item[ii)] Define
$ \cA_{\bH | \bU} =  \{ \Phi \in \cA \mid \Phi|_{{\bH/ \bU }} 
= \id_{\bH/\bU} \}  $, 
then $ \cA_{\bH | \bU} \cap \cAS = \cASid$.

\item[iii)]  $\cA_{\bH | \bF} \cap \cAS = \cASid$. 
\end{itemize}
\end{lemma}
\begin{proof} Let $\Phi \in \cASid$ and $u \in \bU$. 
Since  $\Phi|_{{\bH^\circ /\bF}} = \id_{\bH^\circ / \bF}$,
we have
$\Phi(u) = u f_u$, where $f_u \in \bF$. Now let $h \in \bH$, 
and write $h = s u$, where $s \in \bS$, $u \in \bU$. It follows
that $\Phi(h) = h f_u$. Hence, $\Phi \in \cA_{\bH | \bF}$. This 
proves i). 

Let $\Phi \in  \cA_{\bH | \bU} \cap \cAS$ and $s \in \bS$.
Then $\Phi(s) = s u_s$, where $u_s  \in \bU$. Since $\Phi \in \cA_S$, 
$\Phi(s) \in \bS$. It follows that
$u_s = 1$. Therefore, $\Phi \in \cASid$. This proves
ii). 

By i), $\cASid \subset  \cA_{\bH | \bF}$. Conversely, 
$\cA_{\bH | \bF} \cap \cA_S \leq  \cA_{\bH | \bU}  \cap \cA_S$,
and hence by ii),  $\cA_{\bH | \bF} \cap \cA_S \leq \cASid$.
This proves iii). 
\end{proof}

Before we proceed let us introduce some notation concerning inner
automorphisms.

\begin{definition}\label{innn} 
Let $A$ be a group and $B\le A$ a subgroup. Let $a \in A$. 
We write $\Inn_a\in \Aut(A)$ for the inner
automorphism of $A$ defined by $g \mapsto a g a^{-1}$, for all $g \in A$. 
Given an element $b\in B$, 
we set $\Inn_b^A\in \Aut(A)$ for the corresponding inner 
automorphism of $A$ (to distinguish it from the induced inner 
automorphism of $B$). We write $\Inn_B^A$ for the subgroup of $\Aut(A)$
consisting of all elements  $\Inn_b^A$, $b \in B$.
\end{definition}
Let $\bC$ be Cartan subgroup of $\bH$ which contains $\bS$. 
Thus $\bC = \bU_\bC\cdot \bS$. We define
\begin{equation}  
\cA_\bC :=   \{ \Phi \in \N_\cA(\bF)  \mid   \Phi|_{{\bH^\circ /\bF}} = \id_{\bH^\circ  / \bF} \, , \,  
\Phi(\bC) \subset \bC \} \; .
\end{equation}

If $\bC$ is defined over $\bbQ$ then $\cA_\bC$
is a $\bbQ$-defined subgroup of $\Aut_a(\bH)$.

\begin{lemma} \label{lemma:subgroups2}
The following hold:  \begin{itemize}

\item[i)] Let $u \in \bU_\bC$ such that $\Inn^\bH_u \in \cA_{\bH | \bF}$.
Then there exists $v \in \bC \cap \bF$ 
such that, for all $s \in \bS$,  
$\Inn^\bH_u\,  (s) \;  = \;  \Inn^\bH_v \,  (s) $.

\item[ii)]  $ \cA_{\bH | \bF} \cap \cA_\bC = \Inn^\bH_{\bF \cap \bC}\cdot 
\cASid$.   

\end{itemize}
\end{lemma} 
\begin{proof} 
Let $s \in \bS$. Then $ \Inn^\bH_u (s) =  
s f_s$, where $f_s  \in \bF \cap \bC$.  
Moreover, $f_s =1$, for $s \in \bS^\circ$ since $u$ normalizes $\bS^\circ$. 
By a standard argument (compare for example \cite{BS})
the cocycle $s \mapsto f_s$ is of the form $f_s = s^{-1} v s v^{-1}$, for
some $v \in \bF \cap \bC$. Thus, $\Inn^\bH_v (s) =   \Inn^\bH_u  (s)$, 
for all $s \in \bS$.   

Let $\Phi \in  \cA_{\bH | \bF} \cap \cA_\bC$. Since 
$\Phi(\bC) \leq \bC$, 
there exists $v \in \U_\bC$ 
such that $ v \bS v^{-1} = \Phi(\bS)$. Then 
$\Psi = \Inn^\bH_{v^{-1}}\circ \Phi \in \cAS \cap \cA_{\bH|\bU}$ holds. 
By Lemma \ref{lemma:subgroups1} (ii), $\Psi \in \cASid$ follows.
In particular, $\Psi(s) = s$, for all $s \in \bS$. 
Hence, $\Phi(s) = \Inn^\bH_v(s)$, for all $s \in \bS$. 
Since we have $\Phi \in 
\cA_{\bH | \bF}$,  this implies that $\Inn^\bH_v \in  \cA_{\bH | \bF}$ holds.
By the first part,  there exists $w \in \bC \cap \bF$ such that 
$\Inn^\bH_w (s) = \Inn^\bH_v (s)$ holds for all $s \in \bS$. 
In particular we find 
$\Inn^\bH_{w^{-1}}\circ \Phi \in \cAS \cap \cA_{\bH | \bF} = 
\cASid$. 
Hence, $\cA_{\bH | \bF} \cap \cA_\bC \subset \Inn^\bH_{\bF \cap \bC}
\cdot \cASid$ holds.
The lemma follows.
\end{proof}

\begin{proposition} \label{cAHF}
Let $\bH$ be a
$\BQ$-defined solvable-by-finite linear algebraic group and 
$\bF \leq \bU_\bH$ be a 
$\bbQ$-closed subgroup  
which contains $[\bH^\circ, \bH^\circ]$, then
$$ \cA_{\bH | \bF} \, = \, \Inn^\bH_\bF  \cdot \cASid$$ holds.
\end{proposition} 
\begin{proof}  
It is clear that $\Inn^\bH_\bF\cdot \cASid 
\leq  \cA_{\bH | \bF}$.  
Now let $\Phi \in  \cA_{\bH | \bF}$. Since $\bF$ contains 
$[\bH^\circ, \bH^\circ]$, there exists $v \in \bF$ such
that $ v \bC v^{-1} = \Phi(\bC)$. Therefore, $\Inn^\bH_{v^{-1}}\circ \Phi 
\in \cA_\bC \cap  \cA_{\bH | \bF}$ holds. The proposition follows from 
the previous lemma, part ii).  
\end{proof}


\section{The algebraic hull of a polycyclic-by-finite
wfn-group}\label{sechull}


Let $\Gamma$ be a polycyclic-by-finite group. Its maximal nilpotent
normal subgroup $\Fitt(\Gamma)$ is called the Fitting subgroup of
$\Gamma$. We assume that $\Fitt(\Gamma)$ is 
torsion-free and $\C_\Gamma\left(\Fitt(\Gamma)\right) \leq \Fitt(\Gamma)$.
These two conditions are equivalent
to the requirement that $\Gamma$ has no non-trivial finite normal
subgroups. We call a group with this property a {\it wfn-group}.
Proofs of the following results may be found in \cite[Appendix A]{Baues2}.
 
\begin{theorem} \label{alghull} 
Let $\Gamma$ be a polycyclic-by-finite wfn-group.  
Then there exists a $\bbQ$-defined linear algebraic group $\bH_\Gamma$ and
an injective group homomorphism $\psi: \Gamma \to 
\bH_\Gamma({\bbQ})$ 
such that: \begin{itemize}

\item[i)] $\psi(\Gamma)$ is Zariski-dense in $\bH_\Gamma$,
 
\item[ii)] $\bH_\Gamma$ has a strong unipotent radical $\bU=\bU_{\bH_\Gamma}$,

\item[iii)] $\dim \bU =  \rank \Gamma$. 

 \end{itemize}
 
Moreover, $\psi(\Gamma) \cap \bH_\Gamma({\bbZ})$ is of finite index in 
$\Gamma$. 
\end{theorem} 
Here,  ${\rank}\Gamma$ denotes the number of
infinite cyclic factors in a composition series of $\Gamma$.
(This invariant is sometimes also called the Hirsch-rank of $\Gamma$.)

We remark that the group $\bH_\Gamma$ is determined by the conditions
i)-iii) up to $\bbQ$-isomorphism of algebraic groups:
 
\begin{proposition} \label{prop:rigidity1}
Let $\Gamma$ be a polycyclic-by-finite wfn-group. 
Let $\bH'$ be a $\bbQ$-defined 
linear algebraic group and $\psi': \Gamma \longrightarrow \bH'({\bbQ})$ 
an injective homomorphism  which satisfies 
i) to  iii) from above. Then there 
exists a $\bbQ$-defined
isomorphism $\Phi: \bH_\Gamma \rightarrow \bH'$ such that 
$\psi' = \Phi \circ \psi$.   
\end{proposition}

\begin{corollary} \label{cor:rigidity2}
Let $\Gamma$ be a polycyclic-by-finite wfn-group. 
The algebraic hull\/ $\bH_{\Gamma}$ 
of\/ $\Gamma$ is unique up to $\bbQ$-isomorphism
of algebraic groups. 
In particular,  every automorphism $\phi$ of\/ $\Gamma$ 
extends uniquely to a
$\bbQ$-defined automorphism 
$\Phi$ of\/ $\bH_{\Gamma}$. 
\end{corollary}
We call the $\bbQ$-defined linear algebraic group $\bH_\Gamma$ 
the {\em algebraic hull for $\Gamma$}. We shall identify $\Gamma$ with the
corresponding subgroup of its algebraic hull $\bH_\Gamma$. 

If $\Gamma$ is finitely generated torsion-free nilpotent then
$\bH_\Gamma$ is unipotent and Theorem \ref{alghull} and Proposition
\ref{prop:rigidity1} are essentially
due to Malcev \cite{Malcev}. If $\Gamma$ is torsion-free polycyclic,
Theorem \ref{alghull} is due to Mostow \cite{Mostow1} 
(see also \cite[\S IV, p.74]{Raghunathan} for a different
proof). 

\begin{proposition} \label{Fittu} 
Let $\Gamma$ be a polycyclic-by-finite wfn-group.
Let $\bH_\Gamma$ be the algebraic hull for $\Gamma$.
Then $\Gamma\cap \bU_{\bH_\Gamma}={\rm Fitt}(\Gamma)$ holds.
\end{proposition} 

\begin{definition} We define $\bF =\bF_\Gamma := 
\ac{\Fitt(\Gamma)} \leq \bH_\Gamma$ as
the Zariski-closure of the Fitting subgroup of $\Gamma$.
\end{definition}

Thus, in particular, $\bF$ is a connected unipotent normal
subgroup of $\bH_\Gamma$, and $\bF$ is defined over $\bbQ$.
Moreover: 

\begin{proposition} The commutator subgroup 
$[\bH_\Gamma^\circ, \bH_\Gamma^\circ]$ is 
contained in $\bF$. Let $\bC$ be a Cartan subgroup 
of\/  $\bH_\Gamma$.  Then
there is a decomposition 
\begin{equation} \label{eq:fcdecomp}
 \bH_\Gamma \, = \, \bF \cdot \bC \;  
\end{equation}
\end{proposition}
\begin{proof} Since $\Gamma$ is Zariski-dense, we see that 
$[\bH_\Gamma^\circ, \bH_\Gamma^\circ] = 
\ac{[\Gamma \cap \bH_\Gamma^\circ , \Gamma \cap \bH_\Gamma^\circ]} 
\leq \ac{\Fitt(\Gamma)}= \bF$. The decomposition of $\bH_\Gamma$
follows (see SG7 of Section \ref{prereq}).
\end{proof}


\subsection{A faithful rational representation of $\Aut(\Gamma)$}
\label{auto111}


Let $\Aut_a\!\left(\bH_\Gamma\right)$ 
be the group of algebraic automorphisms of $\bH_\Gamma$ with 
its natural structure of $\bbQ$-defined linear algebraic group. 
We view $\Aut_a(\bH_\Gamma)$ as a $\bbQ$-defined closed subgroup
of $\GL(n,\bbC)$, for some $n \geq 0$. 
Then Proposition \ref{rationalisq} shows that the extension 
\begin{equation} \label{eq:extensionofphi}
 \Aut(\Gamma) \ni \; \; \phi \, \mapsto \, \Phi  \; \; \in \,\Aut_{a,\bbQ}\left(\bH_\Gamma\right) = \Aut_{a}\! \left(\bH_\Gamma\right) \left(\bbQ\right)
\end{equation}
which is  defined in Corollary \ref{cor:rigidity2} gives rise  
to  a faithful representation of $\Aut(\Gamma)$ 
into $\GL(n,\bbQ)$:  

\begin{corollary} The extension homomorphism \eqref{eq:extensionofphi} 
is  a faithful homomorphism of $\Aut(\Gamma)$ into the 
group of\/ $\bbQ$-points of the linear algebraic group $\Aut_a(\bH_\Gamma)$.  
\end{corollary} 

Using this fact, if an embedding $\Gamma \leq \bH_\Gamma$ is
fixed we identify $\Aut(\Gamma)$ with
a subgroup of $\Aut_a(\bH_\Gamma)$.
Remark though, that, in general, $\Aut(\Gamma)$ is not
Zariski-dense in $\Aut_a(\bH_\Gamma)$ because the
elements of $\Aut(\Gamma)$ preserve the Fitting subgroup
$F$, and hence also the Zariski-closure $\bF$ of $F$. Thus, 
with the conventions from Section \ref{section:subgroups}, 
$\Aut(\Gamma)$ is contained in the subgroup ${\rm N}_\cA(\bF)$
of $\Aut_a(\bH_\Gamma)$.

Let $\AGF\le \Aut(\Gamma)$ be the subgroup defined in (\ref{defiagf}).
Since $\Gamma$ is Zariski-dense in $\bH_\Gamma$, the following is clear.
\begin{lemma}\label{agfinAGF}
We have $\AGF\le \cA_{\bH | \bF}(\BQ)$ under the extension 
homomorphism \eqref{eq:extensionofphi}.
\end{lemma}

 
\subsection{Thickenings of $\Gamma$ in $\bH_\Gamma$}


Before introducing the thickening,
we recall some results about finitely generated subgroups in
unipotent algebraic groups.
Let $\bU$ be a unipotent $\bbQ$-defined linear algebraic group, and
let $ F \leq \bU({\bbQ})$ be a finitely generated subgroup. Then 
$F$ is a torsion-free nilpotent group.
Let $\bF:= \ac{F} \leq \bU$ be the Zariski-closure of $F$. 
Then
$\bF$ is $\bbQ$-defined, and $\dim \bF = \rank F$. The group 
of $\bbQ$-points of $\bF$ is isomorphic to the Malcev 
radicable hull of $F$, i.e., $\bF({\bbQ})$ is radicable, 
and for every $x \in \bF({\bbQ})$ there
exists $k \in \bbN$ such that $x^k \in F$. For $m \in \bbN$,
we define 
$$   F^{\frac{1}{m}} \, :=   
\; \langle  x \in \bF({\bbQ}) \mid x^m \in F \rangle \; . $$   
Then $F^{\frac{1}{m}} \leq \bF({\bbQ})$ is finitely generated,
$F \leq F^{\frac{1}{m}}$ and 
$|F^{\frac{1}{m}} : F| < \infty$.    
Every finitely generated subgroup $G \leq \bF({\bbQ})$
is contained in $F^{\frac{1}{m}}$, for some $m \in \bbN$.
See \cite{Segal} for more details on all of this. 

Now let $\Gamma$ be a polycyclic-by-finite wfn-group,  
let $F$ denote the Fitting subgroup of $\Gamma$,
and let $\bF \leq \bH_\Gamma$ be 
the Zariski-closure of $F$ in the algebraic hull 
of $\Gamma$. Since $\bF$ is unipotent,  
$F^{\frac{1}{m}}$ is defined as a subgroup of $\bF(\bbQ)$.  
\begin{definition} \label{def:thickening}
A subgroup $\tilde{\Gamma}$ of $\bH_\Gamma$
which is of the form $\tilde{\Gamma} = F^{\frac{1}{m}}\cdot \Gamma$
is called a \emph{thickening} of $\Gamma$. 
\end{definition}
Clearly thickenings $\tilde{\Gamma}$ exist, for every $m \in \bbN$.
A thickening $\tilde{\Gamma}$ of $\Gamma$ 
is a finitely generated subgroup of $\bH_{\Gamma}(\BQ)$ which
is of finite index over $\Gamma$. We further remark 
that $\Fitt(\tilde{\Gamma}) = F^{\frac{1}{m}}$.  
The inclusion of $\tilde{\Gamma}$ into $\bH_{\Gamma}$ 
shows that $\bH_{\Gamma}$ is an algebraic hull also for  
the thickening $\tilde{\Gamma}$. 


\subsection{The automorphism group of the thickening}


Let $\Gamma$ be a polycyclic-by-finite wfn-group and let $\tilde{\Gamma} \leq \bH_\Gamma(\bbQ)$ be a thickening of $\Gamma$.
Let $\phi \in \Aut(\Gamma)$, and let  
$\Phi: \bH_\Gamma  \rightarrow \bH_\Gamma$ 
denote the extension of $\phi$ in $\Aut_a(\bH_\Gamma)$. 
Since the automorphism $\Phi$ preserves  $\bH_\Gamma(\bbQ)$, it is
clear that  $\Phi$ preserves $\tilde{\Gamma} \leq \bH_\Gamma(\bbQ)$
as well. Restricting $\Phi$ to $\tilde{\Gamma}$ we thus obtain a 
natural inclusion $\Aut(\Gamma) \hookrightarrow
\Aut(\tilde{\Gamma})$. This shows that we may identify $\Aut(\Gamma)$ with 
a finite index subgroup of $\Aut(\tilde{\Gamma})$ in a natural way:

\begin{proposition} \label{findex}
Let $\Gamma$ be a polycyclic-by-finite wfn-group, and let $\tilde{\Gamma}$ be a
thickening of $\Gamma$.  Then
the group $\Aut(\Gamma) = 
\{ \psi \in \Aut(\tilde{\Gamma}) \mid \psi(\Gamma) = \Gamma \, \} $ is a subgroup of finite
index in $\Aut(\tilde{\Gamma})$.
\end{proposition} 
\begin{proof}
Put $d = [ \tilde{\Gamma} : \Gamma ]$  for the index of $\Gamma$
in $\tilde{\Gamma}$. Remark that there are only finitely many
subgroups of $\tilde{\Gamma}$
with index $d$, since $\tilde{\Gamma}$ is a finitely generated group. 
The automorphism group  $\Aut(\tilde{\Gamma})$ acts on the set 
of such subgroups 
and the group $\Aut(\Gamma)$ 
is the stabilizer of the subgroup $\Gamma$. Hence, we have
$ [\Aut(\tilde{\Gamma}) : \Aut(\Gamma) ]  \, \leq \, \ell $, where $\ell$ 
is the number of subgroups of index $d$.
\end{proof}


\section{Thickenings of $\Gamma$ admit a supplement}\label{secthick}


We give here a short account of the construction
of nilpotent-by-finite supplements in polycyclic-by-finite groups. 
Similar results are contained 
in the book \cite{Segal} where nilpotent supplements
in polycyclic groups are considered. 

\begin{definition}
Let $\Gamma$ be a polycyclic-by-finite group and
let $C \leq \Gamma$ be a nilpotent-by-finite subgroup. 
We call $C$ a {\em nilpotent-by-finite supplement in 
$\Gamma$} if $\Gamma = \Fitt(\Gamma)\cdot C$.
\end{definition}

Nilpotent-by-finite supplements do not exist for general 
groups $\Gamma$.
We will show below that a polycyclic-by-finite wfn-group 
$\Gamma$ admits a thickening 
which has a nilpotent-by-finite supplement.

As standing assumption for this section we have that $\Gamma$ is a 
polycyclic-by-finite wfn-group, and we put $F = \Fitt(\Gamma)$. 
We fix an inclusion $\Gamma \leq \bH_\Gamma(\bbQ)$ 
of $\Gamma$ into its algebraic hull $\bH_\Gamma$. 
We put $\bF$ for the Zariski-closure of $F$ in
$\bH_\Gamma$, and  put $N= \bF(\bbQ)$.

\begin{lemma}
Let $\Gamma$ be a polycyclic-by-finite wfn-group and 
let $\bC \leq \bH_\Gamma$ be a $\bbQ$-defined Cartan-subgroup.
Then $\hat{C}= \Gamma N \cap \bC$ is a nilpotent-by-finite
subgroup of $\Gamma \cdot N$ such that  $\Gamma \cdot N = \hat{C}\cdot N$ holds.  
\end{lemma} 
\begin{proof}
The decomposition \eqref{eq:fcdecomp} induces a corresponding
decomposition for the group of $\bbQ$-points of $\bH_\Gamma$,
that is, $ \bH_{\Gamma}(\bbQ) \, = \, N \cdot  \bC(\bbQ) $.

Let $\gamma \in \Gamma$. 
Since $\Gamma \leq \bH_{\Gamma}(\bbQ)$ holds, it follows that 
$\gamma = n_\gamma c_\gamma$, where $n_\gamma \in N$,
$c_\gamma \in \bC(\bbQ) \cap \Gamma N$. Hence, 
$\Gamma \cdot N = \hat{C} \cdot N$ holds.   
\end{proof}

We prove now the existence of supplements in
a thickening $\tilde{\Gamma} = F^{\frac{1}{m}}\cdot\Gamma$.

\begin{proposition}\label{supsplits1}
Let $\Gamma$ be a polycyclic-by-finite wfn-group and
let $\bC \leq \bH_\Gamma$ be a $\bbQ$-defined Cartan subgroup.
Then there exists $m \in \bbN$ such that  
$$ F^{\frac{1}{m}}\cdot\Gamma = F^{1 \over m} \cdot C \, , \, \;
\text{ where } C = (F^{\frac{1}{m}}\cdot \Gamma) \cap \bC \; .$$
\end{proposition} 
\begin{proof} 
By the previous lemma, $\Gamma N = \hat{C} N$, 
where $\hat{C} = \Gamma N \cap \bC$. 
It follows that the natural map 
$\hat{C} \rightarrow \Gamma N/ N$ is surjective.

Since $\Gamma N/ N$ is finitely generated, there exists  a
finitely generated group $C \leq  \hat{C}$ so that 
$C \, \rightarrow \, \Gamma N/ N $ 
is surjective. 
Let $c_1, \ldots, c_k$ be generators for $C$, $c_i= \gamma_i n_i$, where
$\gamma_i \in \Gamma$ and $n_i \in N$. Choose $m \in \bbN$ such that 
$n_i \in F^{\frac{1}{m}}$, $i= 1 \ldots k$. Then $C \leq  
\tilde{\Gamma}= \Gamma F^{\frac{1}{m}}$, and, in particular,
$ F^{\frac{1}{m}} C \leq \tilde{\Gamma}$. 

The surjectivity of the map $C \rightarrow \Gamma N/ N$ 
shows that every $\gamma \in \Gamma$
is of the form $\gamma= c n$, where $c \in C$, $n \in N \cap \tilde{\Gamma} = F^{\frac{1}{m}}$. 
This shows that 
$\tilde{\Gamma} = F^{\frac{1}{m}}\cdot C$. 
\end{proof}
A nilpotent-by-finite supplement is called \emph{maximal} if it 
is a maximal element of the set of all nilpotent-by-finite supplements
with respect to inclusion of subgroups. We show that the maximal
supplements are those which arise by the construction of Proposition \ref{supsplits1}.   

\begin{proposition}\label{supsplits2}
Let $\Gamma$ be a polycyclic-by-finite wfn-group. 
Let $\tilde{\Gamma}= 
F^{\frac{1}{m}}\cdot\Gamma$ be a thickening 
which admits a maximal nilpotent-by-finite supplement $C$.
Then there exists is a $\BQ$-defined 
Cartan subgroup $\bC$
of\/  $\bH_\Gamma$ such that $C=\tilde{\Gamma} \cap \bC$
and  $C$ is Zariski-dense in $\bC$. 
\end{proposition}
\begin{proof}
Let $C$ be any nilpotent-by-finite supplement in $\tilde{\Gamma}$. 
Put $\cC = \ac{C}$ for the Zariski-closure of $C$. Since $C \leq 
\bH_\Gamma(\bbQ)$,
$\cC$ is defined over $\bbQ$. 
Since  $\tilde{\Gamma}= F^{\frac{1}{m}}\cdot C$
is Zariski-dense in $\bH_\Gamma$, we have $\bH_\Gamma = \bF\cdot \cC$. 
Let $\bS$ be a maximal d-subgroup of $\cC$. 
Then $\bS$ is also maximal in $\bH_\Gamma$. 
In particular,
$\cC$ contains a maximal $\bbQ$-defined torus
$\bT$ of $\bH_\Gamma$. Since $\cC$ is nilpotent by-finite,
$\bT$ is unique and normal in $\cC$. We let $\bC$ denote 
the Cartan-subgroup corresponding to $\bT$. Then 
$\cC \leq \bC$. It follows that $C \leq \tilde{\Gamma} \cap \bC$.
In particular, if $C$ is maximal, then $C =  \tilde{\Gamma} \cap \bC$.

Let us prove now  that every
maximal nilpotent by-finite supplement $C$ is Zariski-dense
in the Cartan subgroup $\bC$ which contains $\cC$.  
In fact, since $\bH \leq \bF\cdot \cC$,  $\bC = \U_\bC\cdot \bS$, 
where $\bS \leq \bC$ is a maximal d-subgroup of $\cC$, and
$ \U_\bC$ is a unipotent normal subgroup. Hence, $\bS $ is
a maximal d-subgroup of $\bC$ as well. 
Furthermore, we can write $u \in \U_\bC$
as $u = f u_1$, where $u_1 \in \U_\cC$, $f \in \bF \cap \bC$.
By the maximality of $C$, $F \cap \bC \leq C$.  Therefore,
$C \cap \bF = F \cap \bC$ is Zariski-dense in $\bF \cap \bC$.
This means, $\bF \cap \bC \leq \cC$. Hence, $f \in \cC$ and
$u \in \cC$. It follows that $\cC = \bC$.
\end{proof}

\begin{proposition}\label{subsplits3}
Let $\Gamma$ be a polycyclic-by-finite wfn-group, and
$\tilde{\Gamma}= 
F^{\frac{1}{m}}\cdot\Gamma$ a thickening of $\Gamma$. 
Then there  are at most finitely many $F^{\frac{1}{m}}$-conjugacy classes of  
maximal nilpotent by-finite supplements in $\tilde{\Gamma}$.
\end{proposition}
\begin{proof}
Let $C$ be a maximal nilpotent-by-finite
supplement in $\tilde{\Gamma}$. By 
Proposition \ref{supsplits2}, $C = \Gamma \cap \bC$, 
where $\bC$ is a Cartan subgroup of $\bH$.  
We consider 
$ \tilde{\Gamma}_0 =\tilde{\Gamma} \cap \bH^\circ$,
$C_0 = C \cap \tilde{\Gamma}_0 $. 
Then $\tilde{\Gamma}_0$ is a polycyclic normal
subgroup of $\tilde{\Gamma}$
and $C_0 \, \ns C$.  Also 
$C_0 = C \cap \tilde{\Gamma}_0 = \bC \cap \tilde{\Gamma}_0 
= \bC^\circ \cap \tilde{\Gamma}_0$. Since $\bC^\circ$ is a Cartan subgroup  
in $\bH^\circ$, $C_0$ is a maximal nilpotent supplement in $\tilde{\Gamma}_0$
(see Proposition \ref{supsplits2}).
Note further that $C =  \tilde{\Gamma} \cap \N_{\Gamma N}(C_0)$ is 
uniquely determined by $C_0$.    
By \cite[Chapter 3, Theorem 4]{Segal}, there are only 
finitely many $F^{\frac{1}{m}}$-conjugacy classes of 
maximal nilpotent supplements
$C_0 \leq \tilde{\Gamma}_0$. This also implies that there are only 
finitely many $F^{\frac{1}{m}}$-conjugacy classes of 
maximal nilpotent by-finite supplements in $\tilde{\Gamma}$. 
\end{proof}


\section{Lemmas from group theory}\label{grouplem}


We provide here some simple facts which shall be needed later.

We start off with a few remarks on group cohomology with
non-abelian coefficients. Let $\mu$ be a group, and let $L$ be a group
on which $\mu$ acts by automorphisms. If $s \in \mu$ 
we write $v \mapsto v^s$, $v \in L$ , for the action of $s$ on $L$. 

The set ${\rm Z}^{1}(\mu, L) = \{ z: \mu \to L \mid
z(s_1 s_2) = z(s_1) \, z(s_2)^{s_1} \} $
is called the set of 1-cocycles. 
Two 1-cocycles $z_1$ and $z_2$ are cohomologous if and
only if there exists $v \in L$ such that 
$z_1(s) = v^{-1} z_2(s)  v^s$.   
Let ${\H}^{1}(\mu, F)$
denote the set of equivalence classes of cocycles. It is 
called the first cohomology set for $\mu$ with  
coefficients in $L$. The following lemma is well
known. 
\begin{lemma} \label{lemma:cohomology}
Let $\mu$ be a finite group, and $L$ a 
finitely generated nilpotent
group on which $\mu$ acts by automorphisms. 
Then the cohomology set ${\H}^1(\mu,L)$
is finite.  
\end{lemma} 
\begin{proof} Let $G=L\rfish\mu$ be the split extension corresponding to the
 given action of $\mu$ on $L$. We consider $L$ as a normal subgroup in $G$.
A 1-cocycle $z:\mu\to L$ gives rise to the
 finite subgroup $\mu_z:=\{\, (z(s),s) \mid s\in \mu\,\}$ of $G$. 
Two 1-cocycles $z_1$ and $z_2$ are cohomologous if and
only if the corresponding subgroups $\mu_{z_1}$ and $\mu_{z_2}$ are conjugate
by an element of $L$. Since $G$ is a finitely 
finitely generated nilpotent-by-finite group we know 
(see \cite{Segal}, Chapter 8, Theorem 5) that $G$ has only finitely many
conjugacy classes of finite subgroups. Since $L$ has finite index in $G$, the
lemma follows.  
\end{proof}

We also need the following lemma.
\begin{lemma}\label{InnLem} 
Let $N$ be a group and $M \leq N$ a 
finitely generated torsion-free abelian normal subgroup of finite index. 
Define
$$\Aut(N,M) = \{\, \phi \in \Aut(N) \mid \phi(M)=M,\ \phi|_{M} = \id_M \}$$
then the inner automorphisms $\Inn_M^N$ form a subgroup of finite index in $\Aut(N,M)$. 
\end{lemma} 
\begin{proof} We briefly sketch the argument.  
Remark first that is suffient to prove the 
lemma in the case that the extension $M \leq N$
is effective, that is, ${\rm Z}_N(M) \leq M$. We set $\mu = N/M$.     
Assuming effectiveness, there exists a finite extension 
group $N \leq L$ of $N$ which splits, that is, $L$ is
a semi-direct product  
$L = M_L \rtimes \mu$, where 
$M_L\ge M$ is a torsion-free abelian group which contains $M$ as a
subgroup of finite index. 
Every automorphism of $N$ extends uniquely 
to an automorphism of $L$ which preserves $M_L$. 
Therefore, it is enough to show the lemma for $\Aut(L,M_L)$. 
Now let $\phi \in \Aut(L,M_L)$, that is, $\Phi|_{M_L} = \id_{M_L}$,
and assume additionally  that $\phi$ is the identity on the finite
quotient $L/M_L$. The group of all such $\phi$ is isomorphic to 
the group of 1-cocycles in $Z^1(\mu, M_L)$ with the inner 
automorphisms corresponding to 1-coboundaries.  
Since $\H^1(\mu, M_L)$ (see Lemma \ref{lemma:cohomology}) is finite, 
$\Inn_{M_L}^L$ is of finite index in  $\Aut(L,M_L)$.
\end{proof}
The following can be deduced from \cite{Segal}, Section 6, we skip the proof.
\begin{lemma}\label{uniple}
Let $\bU$ be a unipotent $\BQ$-defined linear algebraic group. 
The following hold: \begin{itemize}

\item[i)] Let $U_1\le U_2\le \bU(\BQ)$ be two finitely 
generated subgroups and suppose that
$U_1$ is Zariski-dense in $\bU$ then the index $[U_2:U_1]$ is finite.

\item[ii)]  Let $U \le \bU(\BQ)$ be a Zariki-dense finitely generated subgroup and let
$d\in\BN$. Let $V \le \bU(\BQ)$ be a subgroup which contains $U$ and satisfies 
$[V:U]\le d$. Then $V$ is contained in $U^{1 \over {d!}}$. In particular, 
the set of all such subgroups $V$ is finite.
\end{itemize}
\end{lemma}


\section{Unipotent shadows of $\Gamma$}\label{secshad}


Let $\Gamma$ be a polycyclic-by-finite wfn-group. We set $F = \Fitt(\Gamma)$
and write $\bF$ for its Zariski-closure in $\bH_\Gamma$. 
Furthermore, we choose a thickening $\tilde\Gamma=F^{\frac{1}{m}}\cdot \Gamma$
which has a (maximal) nilpotent-by-finite supplement $C\le\tilde\Gamma$.  
We use this setup to 
construct (in a controlled way, depending on $\Gamma$)
a finitely generated nilpotent group $\theta \leq \bU(\bbQ)$
which is Zariski-dense in $\bU$. We shall later use $\theta$ to find
arithmetic subgroups in $\Aut(\Gamma)$.

Using the above data we start our construction. We shall use the results of
Section \ref{secthick}.
Let $\bC= \ac{C} \leq \bH_\Gamma$ denote the 
Zariski-closure of $C$. Then $\bC$ is
a $\bbQ$-defined Cartan-subgroup of $\bH_\Gamma$.
We have $C=\tilde\Gamma\cap \bC$, by Proposition \ref{supsplits2}.
Let $\bS \leq \bC$ be a maximal $\bbQ$-defined 
d-subgroup in $\bC$. Then $\bS$ is
a finite extension of the maximal 
torus $\bS^\circ$. The torus $\bS^\circ$ 
is central in $\bC^\circ$, and $\bC= \N_\bH(\bS^\circ)$.
We set $C_0= C \cap \bC^\circ$. Then $C_0$ is a nilpotent
finite index normal subgroup of $C$.

We consider the split decompositions  
$\bC= \bU_\bC \cdot \bS$, $\bC^\circ= \bU_\bC \cdot \bS^\circ$,
where $\bU_\bC$ is the
unipotent radical of $\bC$. Every $c \in \bC({\bbQ})$
can be  written uniquely as
\begin{equation}\label{eledec}     
c = u_c \cdot s_c  \; \; \text{ with } u_c \in \bU_\bC(\bbQ) ,  s_c \in \bS(\bbQ) \; .
\end{equation}
If $c \in \bC^\circ({\bbQ})$ holds, then $s_c\in \bS^\circ(\bbQ)$ 
follows. We define
\begin{equation} 
U_{C,\bS} \, =  \, \langle\, u_c \mid c \in C\,\rangle \, ,  \quad 
U_{C_0} \, =  \, \langle\, u_c \mid c \in C_0\, \rangle \; . 
\end{equation}

\begin{lemma}\label{unipolem}
The groups $U_{C_0}\le U_{C,\bS}$ are finitely generated,
Zariski-dense subgroups of\/ $\bU_\bC$, and $U_{C_0}$ is of finite index
in $U_{C,\bS}$. 
Moreover, $U_{C,\bS}$ is normalized by $C$, and $U_{C_0}$
is normalized by $C_0$.
\end{lemma}
\begin{proof} 
Since $\bS^\circ \leq \bC^\circ$ is central in $\bC^\circ$, 
the map $C_0 \ni c \, \mapsto u_c \in \bU_\bC$ is a homomorphism.
Therefore, the group $U_{C_0}$
is finitely generated. Moreover, since 
$\bU_\bC = \bU_{\ac{C_0}}$, the group $U_{C_0}$ is 
Zariski-dense in $\bU_\bC$. 
Let $S \leq \bS$ denote the image of the homomorphism 
$ C \rightarrow \bS \, , \; c \mapsto s_c$,  
and $S_0 \leq S$ the corresponding image of
$C_0$. The group $C$ acts on $C_0$ and on $\bU_\bC$ by conjugation.
Since $\bS^\circ$ is central in $\bC^\circ$, this action factors over the finite group 
$\mu= S/S_0$. 
For all $c,d \in C$, we have the formula 
\begin{equation}\label{compati}
u_{c d} \, = \, u_c \, s_c u_d s_c^{-1} \; .
\end{equation}
This shows that the group 
$U_{C,\bS}$ is generated by $U_{C_0}$ and a
finite set $u_{c_1}, \ldots u_{c_l}$, where
the $c_i \in C$
represent generators for $C/C_0$.
Therefore, $U_{C,\bS}$ is finitely generated. 
The statement about the finite index follows from Lemma
\ref{uniple}.  

Equation (\ref{compati}) also shows that the action of the finite
group $\mu$ on $\bU_{\bC}$ preserves the subgroup
$U_{C,\bS}$. This implies that $U_{C,\bS}$ is normalized by $C$. The second
statement follows by similar reasoning. 
\end{proof}

\begin{definition}\label{unipodefi}
We define  
\begin{equation}
\theta_{C,\bS} \, 
:= \; \langle\,  F^{\frac{1}{m}},\, U_{C,\bS}\, \rangle \, , \quad 
\theta_{C_0} :=  \langle\,  F^{\frac{1}{m}},\, U_{C_0}\, \rangle \, .
\end{equation} 
The groups $\theta_{C,\bS}$ are called \emph{unipotent shadows} of $\Gamma$. 
\end{definition}
Since we have $\bU=\ac{F}\cdot \bC_\bU$ (see SG7 of Section \ref{prereq}),
Lemma \ref{unipolem} shows that each unipotent shadow
$\theta_{C,\bS}$ 
is a finitely generated subgroup of $\bU(\bbQ)$ which is
Zariski-dense in $\bU$, and it contains the
group $\theta_{C_0}$ as a normal subgroup of finite index. 

\begin{definition}\label{godsha}
We call $\theta_{C,\bS}$ a \emph{good unipotent shadow} 
if the conditions
\begin{equation} \label{eq:FTheta}
\theta_{C_0} \cap \bF =\theta_{C,\bS} \cap \bF = 
F^{\frac{1}{m}}=\Fitt(\tilde\Gamma) 
\end{equation}
are satisfied. 
\end{definition}
Good shadows may be obtained by further
thickening of the Fitting-\-subgroup.

\begin{proposition}\label{supsplits4}
Let $\Gamma$ be a polycyclic-by-finite wfn-group. Then there is a thickening 
$\tilde\Gamma=F^{\frac{1}{m}}\cdot\Gamma$ with a nilpotent-by-finite
supplement $C\le\tilde\Gamma$, such that, for every maximal $\BQ$-defined d-subgroup 
$\bS\le\ac{C}$, $\theta_{C,\bS}$ is a good unipotent shadow. 
\end{proposition} 
\begin{proof} 
Choose $\ell\in\BN$ such that the thickening
$F^{\frac{1}{\ell}}\cdot\Gamma$ admits a nilpotent-by-finite supplement $C$.
Let  $\bS$ be a maximal $\BQ$-defined d-subgroup in the Zariski-closure
$\bC$ of $C$. Let $\theta_{C,\bS}$ be defined as in Definition \ref{unipodefi}. 

Since $\theta_{C,\bS}$ is finitely generated, 
we may choose $m \in \bbN$, divisible by $\ell$, 
such that $\theta_{C,\bS} \cap \bF 
\leq F^{\frac{1}{m}}$.
Now put $\tilde{\Gamma} = F^{\frac{1}{m}}\cdot\Gamma$, and remark that  
$C$ is a nilpotent-by-finite supplement in $\tilde\Gamma$. 
By Proposition \ref{supsplits2},
$C_1= \bC \cap \tilde\Gamma$ is a maximal 
nilpotent-by-finite supplement in $\tilde\Gamma$,
and contains $C$. Since every element $c_1\in C_1$ can be expressed as 
$c_1=fc$ with $c\in C$ and $f\in F^{\frac{1}{m}}$, we find,  going through the
definitions of $\theta_{C,\bS}$ and $\theta_{C_1,\bS}$, 
that $\theta_{C_1,\bS}\cap \bF\le 
F^{\frac{1}{m}}\cdot\big(\theta_{C,\bS}\cap \bF\big)$. 
This implies that $\theta_{C_1,\bS}\cap
\bF=F^{\frac{1}{m}}=\Fitt(\tilde\Gamma)$. 
Hence the requirements of Definition \ref{godsha} follow.
\end{proof}

The following compatibility results are very important for our future
constructions.  
\begin{proposition}\label{supsplits5}
Let $\Gamma$ be a polycyclic-by-finite wfn-group.  Let\/ 
$\tilde\Gamma=F^{\frac{1}{m}}\cdot\Gamma$ be a thickening
of\/ $\Gamma$ with a nilpotent-by-finite
supplement $C\le\tilde\Gamma$. Let
$\theta_{C,\bS}$ be a corresponding unipotent shadow.
Then the following hold: \begin{itemize}

\item[i)] Let $\phi\in \Aut(\tilde\Gamma)$ be an automorphism which satisfies
$\phi(C)=C$, and let $\Phi$ be its extension to an automorphism of\/
$\bH_\Gamma$. Then we have $\Phi(\theta_{C_0})=\theta_{C_0}$.

\item[ii)] For a finite finite index subgroup of the group of all automorphisms $\phi \in \Aut(\tilde\Gamma)$ with 
$\phi(C)=C$, the extension $\Phi$ satsifies
$\Phi(\theta_{C,\bS})=\theta_{C,\bS}$. 

\item[iii)] The group $\tilde\Gamma$ normalizes $\theta_{C,\bS}$.
\end{itemize}
\end{proposition}
\begin{proof} Since $\phi(C) = C$, we have
$\Phi(\bC)=\bC$, hence also
$\Phi(\bC^\circ)=\bC^\circ$ and $\Phi(\bS^\circ)=\bS^\circ$. The definition of 
$U_{C_0}$ shows that $\Phi(U_{C_0})=U_{C_0}$. Since
$\Phi$ also stabilizes $\Fitt(\tilde\Gamma)=F^{\frac{1}{m}}$, it stabilizes
$\theta_{C_0}$. This proves i).

Since $\theta_{C_0}$ is of finite index in 
$\theta_{C,\bS}$,  we can use i) together with part ii) of Lemma \ref{uniple} to
prove ii).

By Lemma \ref{unipolem},  $U_{C,\bS}$ is normalized by $C$. Since $C$ also
normalizes $\Fitt(\tilde\Gamma)=F^{\frac{1}{m}}$ it normalizes
$\theta_{C,\bS}=F^{\frac{1}{m}}\cdot U_{C,\bS}$. The Fitting subgroup
$F^{\frac{1}{m}}$ normalizes $\theta_{C,\bS}$ because it is contained in $\theta_{C,\bS}$.
Hence, $\tilde\Gamma=F^{\frac{1}{m}}\cdot C$ normalizes $\theta_{C,\bS}$.
Hence, iii) holds.
\end{proof}


\section{Arithmetic subgroups of $\Aut(\Gamma)$}\label{secarisag}


Let $\Gamma$ be a polycyclic-by-finite wfn-group. As usually,  
the group $\Gamma$ is considered as embedded in the $\bbQ$-points of its 
algebraic hull $\bH_\Gamma$.  This also fixes an embedding of $\Aut(\Gamma)$ 
in the $\bbQ$-points of the $\BQ$-defined linear algebraic group
$\Aut_a(\bH_\Gamma)$.

We set $F = \Fitt(\Gamma)$
and write $\bF$ for its Zariski-closure.
We assume for this section that $\Gamma$ admits a nilpotent-by-finite supplement.
Thus we may choose a nilpotent-by-finite subgroup $C$ of $\Gamma$ 
such that $\Gamma=F\cdot C$ holds. We choose $C$ maximal with these
properties. 
We write $\bC=\ac{C}$ for its Zariski-closure. Then $C = \Gamma \cap \bC$. 
We further choose a 
$\BQ$-defined d-subgroup $\bS\le\ac{C}$. Associated with these
data comes a unipotent shadow $\theta=\theta_{C,\bS}$,  as constructed in
the previous section. We make the additional assumption that $\theta$ 
is a  good unipotent shadow (see Definition \ref{godsha}).
Our general philosophy is that we always can replace a general wfn-group 
$\Gamma$ by one of its thickenings to enforce these assumptions.

We define $\bU$ to be the unipotent radical of $\bH_\Gamma$. The unipotent
shadow $\theta\le \bU(\BQ)$ provides us with arithmetic subgroups of suitable 
$\BQ$-closed subgroups of $\Aut_a(\bH_\Gamma)$ (compare 
Section \ref{section:asubgroups}). 
We then will find the
position of $\Aut(\Gamma)\le \Aut_a(\bH_\Gamma)$ relative to them. 
Given a subgroup $B\le \Aut_a(\bH_\Gamma)$ we define
\begin{equation}
B[\theta]:=\{\, \Phi\in B\mid \Phi(\theta)=\theta\,\}
\end{equation}
to be the stabilizer of $\theta$ in $B$. We have:  
\begin{lemma}\label{construction} 
Let $\cB\le \Aut_a(\bH_\Gamma)$ be a $\BQ$-closed subgroup which acts faithfully
on $\bU$.  Then $\cB[\theta]$ is an arithmetic subgroup of $\cB$.
\end{lemma}
The lemma follows along the principles used in Section 
\ref{section:asubgroups}, that is, by
linearizing the action on $\bU$ via the exponential function to a linear
action on the Lie algebra of $\bU$. 

As a first application of Lemma \ref{construction}, we obtain that
$\cASid[\theta]$ is arithmetic in $\cASid$, (see
Section \ref{section:subgroups} for the definition of $\cASid$).
We deduce: 
\begin{proposition}\label{construction1}
Given the data $(\Gamma,C,\bS)$ as described above. Then 
$$ \Inn_F^{\bH_\Gamma}\cdot \cASid[\theta]  \, \leq \, 
\cA_{\bH_\Gamma | \bF} (\bbQ)$$   
is an arithmetic subgroup of $\cA_{\bH_\Gamma | \bF}$.
\end{proposition}
\begin{proof}
By our definitions, 
$\cASid[\theta]$ normalizes both $\bF$ and $\theta$, and hence also
$\bF\cap\theta$. We have $\bF\cap\theta=F$ since $\theta$ is a good shadow.
It follows that 
$$F\rfish\cASid[\theta] \, \le \, \bF\rfish\cASid$$
is an arithmetic subgroup. We consider the natural $\BQ$-defined homomorphism
$ \bF \rfish \cASid\to 
\cA_{\bH_\Gamma |\bF}$ which is induced by \eqref{eq:Theta}.  By Proposition
\ref{cAHF},  it is surjective. This implies the result. Of course we have also used that the image
of an arithmetic group under a $\BQ$-defined homomorphism is arithmetic (see 
AR1 of Section \ref{prereq}).
\end{proof}

We turn now to the task of comparing $\Aut(\Gamma)$ to the above 
arithmetic groups.
We define, as in the introduction,
$$   A_{\Gamma |F} :=  \{ \phi \in \Aut(\Gamma) 
\mid \phi|_{\, \Gamma/F} = \id_{\Gamma/F} 
\} \; .$$
Clearly, $A_{\Gamma |F}$ is a characteristic subgroup of $\Aut(\Gamma)$.
Set also 
$$ A_{\Gamma |F}^C:=\{ \phi \in  A_{\Gamma |F} 
\mid \phi(C)=C \}  \; .
$$
We obtain from Proposition \ref{subsplits3}:
\begin{lemma}\label{ter1}
Given the data $(\Gamma,C,\bS)$ described above, 
$\Inn_F^\Gamma\cdot A_{\Gamma |F}^C$ has finite index in 
$A_{\Gamma |F}$.
\end{lemma}

Next we analyze the group $A_{\Gamma |F}^C[\theta]$. We obtain from 
Proposition \ref{supsplits5} ii):
\begin{lemma}\label{ter2}
Given the data $(\Gamma,C,\bS)$ described above,
$A_{\Gamma |F}^C[\theta]$ has finite index in $A_{\Gamma |F}$.
\end{lemma} 

We have the semi-direct product decomposition $\bC=\bU_\bC\cdot \bS$. Relative to
this decomposition we can consider the quotient homomorphism 
$\pi_\bS :\bC\to \bS$,  and define
\begin{equation}\label{esss}
S:=\pi_\bS(C), \qquad S_0:=\pi_\bS(C\cap \bC^\circ)
\end{equation}
Notice that, by the
constructions in Section \ref{secshad}, we have 
$c\cdot \pi_\bS(c)^{-1}\in \theta$,  for every $c\in C$.
We need the following technical observations:
\begin{lemma}\label{ter3}
Given the data $(\Gamma,C,\bS)$, described above and $S,\, S_0$ as
defined in (\ref{esss}), we have \begin{itemize}

\item[i)] The group $S$ normalizes $F\cap C$ and $S_0$ centralizes it:

\item[ii)] Let $\Phi$ be the extension of the automorphism $\phi\in
A_{\Gamma |F}^C[\theta]$ to an automorphism of\/ $\bH_\Gamma$. Then 
$\Phi(s)\cdot s^{-1}\in F\cap C$ holds, for every $s\in S$. If $s\in S_0$ 
then $\Phi(s)\cdot s^{-1}=1$.
\end{itemize}
\end{lemma} 
\begin{proof}
i): Let $s$ be in $S$, choose $c\in C$ with $\pi_\bS(c)=s$ and define 
$v=c\cdot\pi_\bS(c)^{-1}$. As remarked above, we have $v\in \theta$.
Since $c$ normalizes $F\cap C$, the following holds 
$$s \left(F\cap C\right) s^{-1}=v^{-1}  \left(F\cap C\right) v\, .$$
The right hand  side is in $\theta$ and in $\bF$, hence in $\theta\cap \bF=F$.
This implies that the right hand  side is in $D=F\cap \bC$. The subgroup $D\le
\Gamma$ is nilpotent-by-finite and normalized by $C$. Both $C$ and $D$ are
contained in $\bC$. Hence $\langle \, C,\, D\,\rangle$ is a
nilpotent-by-finite supplement in $\Gamma$. Since $C$ is maximal, we have 
$D=F\cap \bC\le C$ and the first part of i) follows. For the second part,
notice that $F\cap C\le \bC^\circ$.  

ii): Let $s$ be in $S$, choose $c\in C$ with $\pi_\bS(c)=s$ and define
$v=c\cdot\pi_\bS(c)^{-1}$. Then
$$\phi(c)c^{-1}=\Phi(v)\Phi(s)s^{-1} v^{-1}.$$  
Since $\Phi$ is the identity modulo $\bF$, the right hand side is in $\bU$,
whereas the left hand side is in $C$. Hence, the right hand side is in $C\cap
\bU\le\theta$. Our assumptions imply $\Phi(s)s^{-1}\in\theta$ and then  
$\Phi(s)s^{-1}\in\theta\cap \bF=F$. Here we have used that $\theta$ is a good
unipotent shadow. Furthermore the above equation shows that 
$\Phi(s)s^{-1}\in\bC$ holds. As under i), we finish the proof 
of the first part of ii) by remarking
that $F\cap \bC$ is contained in $C$. For the second notice that if $s\in
\bS^\circ$ holds then  
$\Phi(s)s^{-1}\in\bF\cap\bS^\circ=\{1\}$ follows.
\end{proof}
We proceed with the construction of subgroups in $\Aut(\Gamma)$. We define
\begin{equation}
A_{\Gamma |F}^C[\theta]^1:=\{ \, \phi\in A_{\Gamma |F}^C[\theta]
\mid \Phi(\bS)=\bS,\ \Phi|_\bS={\rm id}_\bS\,\}.
\end{equation}
Here $\Phi$ is as always the extension of the automorphism
$\phi\in\Aut(\Gamma)$ to an automorphism of $\bH_\Gamma$.
We have:

\begin{lemma}\label{ter4}
Given the data $(\Gamma,C,\bS)$, as above. Then
$\Inn_{F\cap C}\cdot A_{\Gamma |F}^C[\theta]^1$ is of finite index in
$A_{\Gamma |F}^C[\theta]$.
\end{lemma}
\begin{proof}
Notice first that $\Inn_{F\cap C}$ is contained in $A_{\Gamma |F}^C[\theta]$, 
as follows from the definitions. 
Let $S,\, S_0\subset  \bS$ be the subgroups
defined in (\ref{esss}). The quotient group $\mu=S/S_0$ is finite and acts by
conjugation on $F\cap C$ (see Lemma \ref{ter3} i)). We let 
${\rm Z}^{1}(\mu, F \cap C)$ be the corresponding set of 1-cocycles and 
${\H}^1(\mu, F\cap C)$ the cohomology set (see Section \ref{grouplem} for
definitions). This cohomology set is finite,  by Lemma \ref{lemma:cohomology}.

Given $\phi\in A_{\Gamma |F}^C[\theta]$ with extension $\Phi$,  we obtain,
using Lemma \ref{ter3}, a map $D_\phi :S \to F\cap C$ by setting
$D_\phi(s)=\Phi(s)\cdot s^{-1}$.
The verification of the following is straightforward:
\begin{itemize}
\item the above $D$ induces a (well defined) map $D: A_{\Gamma|F}^C[\theta]\to 
{\rm Z}^{1}(\mu, F \cap C)$,
\item the map $D$ from the previous item induces a (well defined) map
$$ \hat D: \; A_{\Gamma|F}^C[\theta] \Big/
\Inn_{F\cap C}^\Gamma \cdot A_{\Gamma|F}^C[\theta]^1  \, 
\longrightarrow \, {\H}^1(\mu, F\cap C) ,$$
\item the map $\hat D$ is injective.
\end{itemize}
As remarked before, ${\H}^1(\mu, F\cap C)$ is finite and the lemma is proved. 
\end{proof}

We put now Lemmas \ref{ter1}, \ref{ter2}, \ref{ter4} together and obtain:
\begin{lemma}\label{ter5}
Given the data $(\Gamma,C,\bS)$,  as above. Then
the group
$\Inn_{F}^\Gamma\cdot A_{\Gamma |F}^C[\theta]^1$ is of finite index in   
$A_{\Gamma |F}$.
\end{lemma}
The link between Proposition \ref{construction1} and Lemma \ref{ter5} is given
by:

\begin{proposition}\label{terfin}
Given the data $(\Gamma,C,\bS)$,  as above. Then
$A_{\Gamma |F}^C[\theta]^1=\cASid[\theta]$.
\end{proposition}
\begin{proof}
By the definitions, 
$A_{\Gamma |F}^C[\theta]^1\le\cASid[\theta]$ holds.
Now let $\Phi$ be an element of $\cASid[\theta]$.
We show that $\Phi$ is contained in $\Aut(\Gamma)$.
First of all, we have $\Phi(F) = \Phi(\theta \cap \bF)
= \Phi(\theta) \cap \Phi(\bF) = F$. Here we used that $\theta$ is a
good unipotent shadow. 
Since $\Phi|_{\bH^\circ/\bF} = \id_{\bH^\circ/\bF}$, it follows
that 
$ \Phi|_{{\theta/F}} \, = \, \id_{\theta/F}.$ 
Let $c \in C$, $c= u s$ with $u \in \theta$ and $s \in \bS$. 
By the above,  $\Phi(u) = f u$, for some $f \in F$.  Therefore, we get 
$$  \Phi(c) \, =     \Phi(u) \Phi(s) = 
\Phi(u) s = f u s= f c \, \in \Gamma \;  .$$
Since $\Phi(\bC) = \bC$, it follows that $\Phi(c) \in \Gamma \cap \bC = C$. 
This shows that $\Phi$ stabilizes $C$. 
Hence, $ \Phi(\Gamma)= \Phi(F C) = F C =  \Gamma.$ 
Thus,  $\Phi \in \Aut(\Gamma)$ holds. 
The lemma follows.
\end{proof}
Putting together Lemma \ref{ter5}, Proposition \ref{terfin} and 
Lemma \ref{agfinAGF}, we obtain:
\begin{corollary} \label{AGFisa1} 
The group $A_{\, \Gamma | F}$ is an 
arithmetic subgroup of $\cA_{\bH_\Gamma|\bF}$. 
\end{corollary}

Finally, let us consider the arithmetic subgroup $\A_\theta$ 
(see Definition \ref{LThedef}) of the group $\Aut_a(\bH_\Gamma)$.
The group $A_{\tilde\Gamma |\tilde F}[\theta]$ is defined as 
the stabilizer of $\theta$ in $\AGF$. We note: 

\begin{proposition}  \label{propcontat}
A finite index subgroup of\/ 
$\Inn_{\Gamma}\cdot \AGF[\theta]$ 
is contained in $\A_\theta$.   
\end{proposition}
\begin{proof} 
By construction, we have 
$\AGF^C[\theta]^1 \le \A_{\theta}$. The group 
$\Inn_{\Gamma}$ stabilizes $\theta$,  by Proposition \ref{supsplits5} iii).
Obviously  $\Inn_{\tilde F}$ is 
contained in $A_\theta$. We have $\Gamma=\tilde F\cdot C$.  
Let $c$ be an element of $C$,  we write $c=v\cdot s$ with $v\in \theta$, $s\in\bS$. Then 
$\Inn_v \in \A_\theta$,  and hence $\Inn_s$ stabilizes $\theta$. Clearly,
$\Inn_s \in \Aut_a(\bH_\Gamma)_\bS$, therefore 
$\Inn_s$ is also in $\A_\theta$. This shows that
$\Inn_{\Gamma}$ is contained in $\A_\theta$.  
Now $\Inn_{\Gamma}\cdot A_{\Gamma |\tilde F}[\theta]^1$ 
is a finite index subgroup of 
$\Inn_{\Gamma}\cdot \AGF[\theta]$, by Lemma \ref{ter5}. 
\end{proof}


\section{The automorphism group of $\Gamma$ as a subgroup of 
$\Aut_a(\bH_{\Gamma})$}\label{secautoauto}


This section contains the final proof of Theorem \ref{teoagf}. We
also provide the input for the proofs of Theorem  \ref{teoc} and 
Theorem \ref{teoa}. (These proofs will
be given in Section \ref{secari}.)

Let $\Gamma$ be a polycyclic-by-finite wfn-group. We stick to our usual conventions. 
Namely, $\Gamma$ is embedded in the $\bbQ$-points of its 
algebraic hull $\bH_\Gamma$,  $\Aut(\Gamma) \leq \Aut_a(\bH_\Gamma)(\bbQ)$, 
and  $\bU$ is the unipotent radical of $\bH_\Gamma$.
We also fix,  as a reference, a thickening 
$\tilde\Gamma=F^{\frac{1}{m}}\cdot\Gamma$ of $\Gamma$ in 
$\bH_\Gamma$. We choose a thickening which satisfies the assumptions of Section \ref{secarisag}
on the data $(\tilde{\Gamma},C, \bS)$. In particular, $C \leq \tilde{\Gamma}$ is a
maximal nilpotent-by-finite supplement,  $\bC=\ac{C}$ its Zariski-closure, $\bS \leq \bC$ 
a maximal $\BQ$-defined $d$-subgroup,  and
$\theta = \theta_{\tilde{\Gamma}}$ is a good unipotent shadow for $\tilde{\Gamma}$.  
Such a thickening exists,  by 
Proposition \ref{supsplits4}.

\begin{proposition} \label{InnGAGFfi} 
Let $\Gamma$ be 
a polycyclic-by-finite wfn-group.
Then the subgroup $\Inn_\Gamma\cdot A_{\Gamma |F}$ 
has finite index in $\Aut(\Gamma)$.   
\end{proposition}

\begin{proof} Recall that $A_{\Gamma |F} = \{ \phi \in \Aut(\Gamma) \mid \phi|_{\Gamma/F} = 
\id_{\Gamma/F} \} $. We shall use that 
\begin{equation}
\AGF \, = \, \Aut(\Gamma) \cap \cA_{\bH_\Gamma|\bF} 
 \, = \, \Aut(\Gamma) \cap \cA_{\bH_\Gamma|\bU}.
\end{equation}
This is  a straightforward consequence of Proposition \ref{Fittu}.

Let us put  $N =  \pi_{\bU}(\Gamma)$, and $M= \pi_{\bU}(\Gamma_0)$,
where $\Gamma_0 = \Gamma \cap \bH^\circ$. Define 
$\hat{\bS} = \bH_\Gamma \big/ \bU$ and note that $\hat{\bS}$ is a
$\BQ$-defined d-group.

Let $\phi \in \Aut(\Gamma)$ and let $\Phi\in\Aut_a(\bH_\Gamma)$ be its
extension to $\bH_\Gamma$. The $\BQ$-defined automorphism $\Phi$
induces a $\BQ$-isomorphism 
$\Phi_{\hat{\bS}}$ of $\hat{\bS}$, 
which preserves $N$ and $M$. 
The restriction
of $\Phi_{\hat{\bS}}$ to $N$ will be denoted by
$\phi_N$. By the rigidity of tori (AG7), 
for all $\phi$ in a finite index subgroup of $\Aut(\Gamma)$,
$\phi_N$ is the identity on $M$, that is $\phi_N  \in \Aut(N,M)$. 
Thus, by Lemma \ref{InnLem}, $\phi_N \in \Inn^N_M$ holds
in a finite index subgroup of $\Aut(\Gamma)$. 
If $\phi_N  \in \Inn_M^N$,  there exists 
$c \in \Gamma_0$ 
such that $\left(\Inn_c^\Gamma \circ \phi\right)_N = \id_N$.
Since $N$ is Zariski-dense in
$\hat{\bS}$, this implies 
$\Inn_c^\Gamma\circ \phi \in \Aut(\Gamma) \cap \cA_{\bH_\Gamma| \bU} 
= \A_{\Gamma|F}$.
Therefore, $\Inn_{\Gamma_0}^\Gamma \cdot \A_{\Gamma|F}$ is of finite
index in $\Aut(\Gamma)$. 
\end{proof}

\begin{prf}{Proof of Theorem \ref{teoagf}.}
By the results  proved in Section \ref{auto111},  $\Aut(\Gamma)$ is
contained in the $\BQ$-points of $\Aut_a(\bH_\Gamma)$. 
We come now to the statement about $\AGF$. We use that 
$\Aut(\Gamma)$ is naturally contained in 
$\Aut(\tilde\Gamma)$ as a subgroup of finite index, see Proposition \ref{findex}. 
The 
obvious fact that
$\tilde F=\Fitt(\tilde\Gamma)$ satisfies $\tilde F\cap\Gamma=F$ implies that
$\AGF = A_{\tilde\Gamma |\tilde F} \cap \Aut(\Gamma)$ 
is of finite index in $A_{\tilde\Gamma |\tilde F}$. Hence, by Corollary
\ref{AGFisa1},  $\AGF$ is an arithmetic subgroup of
$\cA_{\bH_\Gamma |\bF}$.

Now we prove that $\Aut(\Gamma)$ is,  up to finite index, contained in the
arithmetic group $\A_\theta$  defined relative
to the good unipotent shadow $\theta = \theta_{\tilde{\Gamma}}$. 
(For the construction of $\A_\theta$ refer to  Definition \ref{LThedef}.)  
Since $\Aut(\Gamma)$ is with finite index naturally contained in 
$\Aut(\tilde\Gamma)$ (see Proposition \ref{findex}) it is enough to prove 
that a finite index subgroup of  $\Aut(\tilde\Gamma)$ is contained 
in $\A_\theta$. The latter is implied  by 
Proposition \ref{propcontat} and Proposition \ref{InnGAGFfi} .
\end{prf}

\begin{prf}{Proof of Theorem \ref{teoc} for wfn-groups.}
By Proposition \ref{InnGAGFfi},  we have that
$\Inn_\Gamma \cdot  \AGF$ is of finite index in $\Aut(\Gamma)$.
Moreover,  $\tilde\Gamma =\tilde F\cdot C$ contains 
$\Gamma$ as a subgroup of finite index. Hence, 
$ F \cdot (\Gamma \cap C)$ is of finite index in $\Gamma$. 
We have $\Inn_{ F}$ is in 
$\AGF$. Therefore, $\Inn_C^\Gamma \cdot \AGF$ is of finite
index in $\Aut(\Gamma)$.  Now choose a finite index invariant 
nilpotent subgroup $B$ of $\Inn_C^\Gamma$ to obtain
the result. 
\end{prf}

The usual induction procedure gives the following 
immediate corollary of Theorem \ref{teoagf}:
\begin{corollary}\label{imBett} 
Let $\Gamma$ be 
a polycyclic-by-finite wfn-group.
Then there exists a faithful 
representation of $\Aut(\Gamma)$ into 
$\GL(n, \bbZ)$, for some $n \in\BN$. 
\end{corollary}


\section{Extensions and quotients of arithmetic groups}\label{subsecari}


In the following we shall accumulate some results about extensions and
quotients of arithmetic groups which will be necessary for the proofs in the next
subsection. To formulate these results we need the following 
concepts. 

\begin{definition} Let $\cA$ be a $\BQ$-defined 
linear algebraic group and $A\le \cA$ a
subgroup. An automorphism $\phi$ of $A$ is said 
to be \emph{$\cA$-rational} if there is a
$\BQ$-defined automorphism of $\cA$ which normalizes $A$ and coincides with
$\phi$ on $A$. Moreover, a homomorphism $\rho: A \ra \cG(\bbQ)$ 
into a $\bbQ$-defined linear algebraic group $\cG$ is called 
\emph{$\cA$-rational} if $\rho$ extends to a $\bbQ$-homomorphism 
$\rho_\cA: \cA \ra \cG$.   
\end{definition}

\begin{definition}
Let $G$ be a group of automorphisms of $A$.
If the action of $G$ on $A$ extends to an algebraic group of automorphisms on $\cA$  (see Definition \ref{algauto}) then $G$ is said to be an  \emph{$\cA$-algebraic group of automorphisms} of $A$. 
\end{definition}

As explained in the introduction, a finite extension group of an arithmetic
group need not be arithmetic.  We shall give now a slight generalization of a criterion
from \cite{GP3} which allows to show that certain finite 
extension groups of arithmetic groups are again arithmetic. 

\begin{lemma}\label{crit}
Let $\cA$ be a $\BQ$-defined linear algebraic group and $A\le \cA$ a
Zariski-dense
arithmetic subgroup. Let $B\ge A$ be a group containing $A$ as a 
normal subgroup of finite index. Suppose that conjugation of $B$ on
$A$ is  $\cA$-rational. Then the following hold:  

\begin{itemize}
\item[i)] The inclusion of $A$ into $\cA$ can be extended to an embedding of
the group $B$ as an arithmetic subgroup into a 
$\BQ$-defined linear algebraic group ${\cal B}$ which contains $\cA$ as a
subgroup of finite index. 

\item[ii)] Every automorphism of $B$ which
normalizes $A$ and which induces an $\cA$-rational automorphism of $A$   
is ${\cal B}$-rational.

\item[iii)] Let $G$ be a group which acts by automorphisms 
on $B$ which normalize $A$. If $G$ acts
as an $\cA$-algebraic group of automorphisms on $A$
then $G$ is a ${\cal B}$-algebraic group of automorphisms of $B$.  
\item[iv)] Let $\rho:  B \ra \cG(\bbQ)$ be a representation of $B$ which
restricts to an $\cA$-rational representation of $A$. Then $\rho$ is
$\cB$-rational. 
\end{itemize}
\end{lemma}
\begin{proof}  
Statement i) is an application  of \cite[Proposition 2.2]{GP3}. 
The $\BQ$-defined linear algebraic group ${\cal B}$ is
constructed through the usual induction procedure, and 
if $R = \{ r_1,\ldots,r_n \} \subset B$ 
is a complete set of coset representatives for $A$ 
in $B$ then $R$ also forms a complete set of 
coset representatives for $\cA$ in ${\cal B}$. 
There exists thus $r_{ij} \in R$ and $a_{ij} \in A$ such that
$ r_i r_j =  r_{ij} a_{ij}$.

We show iv).  Let $\rho:B\to \cG$ be a homomorphism of $B$ whose
restriction to  $A$ is $\cA$-rational. That
is, there is a $\BQ$-defined homomorphism $\rho_\cA :\cA\to \cG$ 
with $\rho_\cA(a)=\rho(a)$
for all $a\in A$. We shall show now that $\rho$ is ${\cal B}$-rational.
Note first that, for all $b\in B$ and all $a\in\cA$, we have
$$ \rho_\cA(b^{-1} a b)=\rho(b)^{-1} \rho_\cA(a)\rho(b) \; , $$
since this identity is valid on the
Zariski-dense subgroup $A$ of $\cA$. 
Define now a map $f:{\cal B}\to {\cal G}$ by 
$$f(r_ia):= \rho(r_i) \rho_\cA(a)\qquad (i=1,\ldots,n,\ a\in\cA).$$
Clearly, $f$ is a $\BQ$-defined morphism of varieties.  A straightforward 
computation using the above mentioned identity shows that $f$ is a
homomorphism of groups. It follows that $f$ is a $\bbQ$-defined 
homomorphism of linear algebraic groups. 
It is clear that $f$ coincides with $\rho$ on $B$.
This proves iv).

Note that  ii) is an immediate consequence of  iv). 

To prove iii) use ii), and note that the condition of being an algebraic group of automorphisms on $\cB$ depends only on the connected component  
${\cal B}^\circ = \cA^\circ$. 
\end{proof}

The following remark is evident from the definition of an
arithmetic group.
\begin{lemma}\label{trivilem} 
Let $\cB$ be a $\BQ$-defined linear algebraic group and $A\le B\le \cB(\BQ)$
be subgroups. Assume that $A$ is of finite index in $B$ and is an
arithmetic subgroup of its Zariski-closure. Then $B$ is an 
arithmetic subgroup of its Zariski-closure in $\cB$.
\end{lemma}
We apply the lemmas just proved to show Proposition \ref{finker}.  

\begin{prf} {Proof of Proposition \ref{finker}.}
Since $A$ is residually finite we can choose a subgroup $C$ of finite index in
$A$ which is normal in $A$ and satisfies $ E\cap C=\{1\}$.
To obtain a normal such $C$, take a finite index subgroup of $A$ with this 
property, and then intersect over (the finite number of) all
subgroups of the same index.

Choose a $\BQ$-defined linear algebraic group $\cA$ and a group
homomorphism $\rho:A\to \cA$ with kernel $E$ such that the image is a
Zariski-dense arithmetic subgroup of $\cA$. The homomorphism $\rho$ is injective 
on $C$ and $\rho(C)$ is an arithmetic subgroup of $\cA$.
Clearly, the  conjugations by elements of $A$ induce 
$\cA$-rational automorphisms of $\rho(C)$. We may 
now finish the proof by using Lemma \ref{crit}.
\end{prf}

Next we study certain arithmetic quotients of arithmetic
groups.

\begin{proposition}\label{propquot}
Let $\cA$ be a $\bbQ$-defined linear algebraic group and 
$A\le \cA$ a Zariski-dense arithmetic subgroup. Let
$N \leq A$ be a normal subgroup of $A$ and let ${\cal N}$ denote the
Zariski-closure of $N$ in $\cA$.  Assume furthermore that $N$ has finite 
index in ${\cal N}\cap A$. Then:
\begin{itemize} \item[i)] 
The group $A/N$ embeds  an 
arithmetic subgroup into a 
$\bbQ$-defined linear algebraic group ${\cal D}$.
 \item[ii)] 
Every $\cA$-rational automorphism of $A$ which normalizes $N$ induces a
${\cal D}$-rational automorphism of $A/N$. 
\item[iii)]  Every $\cA$-algebraic group of  automorphisms on $A$ which normalizes
 $N$ induces a ${\cal D}$-algebraic group of 
automorphisms of $A/N$.
\item[iv)] Let $\rho$ be an $\cA$-rational representation of $A$ with $N \leq \ker \rho$.
Then the induced quotient representation of $A/N$ is $\cD$-rational. 
\end{itemize}
\end{proposition}

For the proof of Proposition \ref{propquot} we need the following
lemma.
\begin{lemma} \label{lemmapropquot}
Under the assumptions of Proposition \ref{propquot} there exists
a finite index subgroup $C \leq A$ such that $\cN \cap C \leq N$. 
\end{lemma}
\begin{proof}
To prove the lemma, we may clearly assume that the
group $\cA$ is connected. 
Let $\bU_\cA$ be the unipotent radical of $\cA$.
The unipotent radical 
$\bU_\cN$ of $\cN$ is contained as a normal subgroup in $\bU_\cA$. 
We may now choose a reductive complement 
$\cA^{{\rm red}}$ for $\bU_\cA$, such that 
${\cal N}^{\rm red} =\cN \cap  \cA^{{\rm red}}$ is
a reductive complement for $\bU_\cN$ 
in ${\cal N}$. In particular, $ {\cal N}^{\rm red}$
is normal in $\cA^{{\rm red}}$. Thus,  by AG5,
we may choose an almost direct 
complement ${\cal H}$ for ${\cal N}^{\rm red}$ in 
$ \cA^{{\rm red}}$.  That is, ${\cal H}$ is a $\bbQ$-defined subgroup of   
$ \cA^{{\rm red}}$ which centralizes  ${\cal N}^{\rm red}$, satisfies 
$\cA^{{\rm red}}= {\cal N}^{\rm red}\cdot {\cal H}$ and has finite
intersection ${\cal N}^{\rm red}\cap {\cal H}$. 
 
We put $N_{\bU}={ \bU_\cN }\cap N$, and
$N_1={\cal N}^{\rm red}\cap N$. 
Since $N$ is arithmetic in $\cN$,  $N_{{\bU}} \cdot N_1$ has finite
index in $N$. Since $N_{\bU}$ is 
arithmetic in the unipotent group ${\bU_\cN }$, there is a congruence subgroup 
$G$ of $A$ with the property ${\bU_\cN } \cap G \le N_{\bU}$ (see 
\cite{Segal}, Chapter 4, Theorem 5). This congruence subgroup
may be chosen torsion-free as well. We now set 
$$G_{\bU} =\bU_\cA \cap G,\quad G_{1} =
{\cal N}^{\rm red} \cap G,\quad G_\cH ={\cal H} \cap G\; .$$
Since $G$ is an arithmetic subgroup of $\cA $, the product 
$G_{\bU} \cdot  G_{1} \cdot G_\cH$ is of finite index in $G$. 
Both $N_1$ and $G_1$ are arithmetic subgroups of 
${\cal N}^{\rm red}$.  Hence, $C_1= N_1\cap G_1$ has finite index in  
$N_1$ and $G_1$. Therefore, 
 $$ C=G_\bU\cdot C_1\cdot  G_\cH $$ has finite index in $G$ and $A$.
Since $G$ is torsion-free, $\cN \cap C$ is contained in 
$G_{\bU} \cdot C_1$. Now we find that 
$${\cal N}\cap C  \le  G_{\bU} \cdot C_1
\le \big( {\bU_\cA} \cap C \big) \cdot \left({\cal N}^{\rm red}\cap G \right)
\le  N_{\bU} \cdot N_1 \leq N \; .$$  
This finishes the proof of the lemma. 
\end{proof}

\begin{prf}{Proof of Proposition \ref{propquot}.}
Replacing the subgroup $C$  constructed in Lemma  \ref{lemmapropquot} 
above by one of its subgroups of
finite index and also possibly by $C\cdot N$, the following can be arranged:
\begin{itemize}
\item $C$ is a normal subgroup of finite index in $A$,
\item $C\cap{\cal N}= N$,
\item $C$ is normalized by every automorphism of $A$ which normalizes $N$.
\end{itemize}
Let ${\cal C}$ denote the Zariski-closure of $C$. This is a $\BQ$-closed subgroup
(of finite index) in $\cA$ which contains ${\cal N}$. Also $C/N$ is contained as a
Zariski-dense subgroup in the 
$\BQ$-defined linear algebraic group ${\cal C}/{\cal N}$. Moreover, $C/N$ is an
arithmetic subgroup of ${\cal C}/{\cal N}$, by AR1. 

Now $B= A/N$ is a finite extension group of $C/N$. Clearly, conjugations by
elements of $A$ give rise to ${\cal C}\big/{\cal N}$-rational automorphisms of
$C/N$ and so does every $\cA$-rational 
automorphism of $A$ which normalizes $N$. Therefore,  the group ${\cal D}$ may be 
constructed by application of Lemma \ref{crit}, thus proving i), ii), iii). 

To prove iv), let $\rho_\cA: \cA \ra \cG$ denote the algebraic extension of $\rho$, 
and let $\rho_{A/N}:  A/N \ra \cG(\bbQ)$ be the quotient 
representation induced by $\rho$. Clearly, its restriction to $C/N$ is $\cC\big/ \cN$-rational,
since $\rho_\cA$ factors over $\cC\big/ \cN$.  Thus, part iv) of Lemma
\ref{crit} shows that  $\rho_{A/N}$ is  $\cD$-rational.
\end{prf}

In the following we deal with the fact that a group which is isomorphic to an
arithmetic group may admit essentially different arithmetic embeddings  
 into linear algebraic groups. This phenomenon plays a role in our 
arithmeticity proofs. At this point we will
also need the full strength of the assumption that $G$ acts as an 
algebraic group of automorphisms on $A \leq \cA$, in order to
extend the action of $G$  to a modification of the ambient group $\cA$. 

\begin{proposition}\label{changerep}
Let  $A$ be an arithmetic subgroup of  a $\bbQ$-defined linear algebraic group 
$\cA$ and $G$ an $\cA$-algebraic group of automorphisms of $A$.
Let $D$ be a normal subgroup of $A$ which is contained in the
center of $A$ and which is normalized by $G$. 
Then the group $A$ can be embedded as an arithmetic subgroup
into a $\bbQ$-defined linear algebraic group ${\cE}$ such that:
\begin{itemize}
\item[i)] The subgroup $D$ of $A$ is unipotent-by-finite in $\cE$.
\item[ii)] The group $G$ acts as an $\cE$-algebraic group of automorphisms of $A$.
\item[iii)] If $\rho: A \ra \cG(\bbQ)$ is an $\cA$-rational representation 
which satisfies $D \leq \ker \rho$ then $\rho$ is  $\cE$-rational. 
\end{itemize}

\end{proposition}
\begin{proof}
Clearly, there is no harm in assuming that $A$ is Zariski-dense in $\cA$.
We choose a subgroup $C\leq A$, subject to the following conditions:
\begin{itemize}
\item  $C$ is torsion-free and a normal subgroup of finite index in $A$,
\item  $C$ is contained in the connected component $\cA^\circ$ of $\cA$, 
\item  $C$ is normalized by  $G$.
\end{itemize}
Note that $C$ is Zariski-dense in $\cA^\circ$.  
We define $G_1 = \Inn_A\cdot G$ 
to be the subgroup of the automorphism group of $A$ 
which is generated by $\Inn_A$ and $G$. Then $G_1$ acts by $\bbQ$-defined 
automorphisms on $\cA$.  Now Lemma \ref{algautoex} and Lemma \ref{algformal3} show 
that $G_1$ is an $\cA$-algebraic group of automorphisms of $A$ (and also of $C$).  
Hence, according to Proposition \ref{decozen},  there exists a  
$G_1$-invariant almost direct product decomposition 
$$\cA^\circ={\cal Z}_1 \cdot \cA_1 \; \; , $$ 
where ${\cal Z}_1$ is the $\BQ$-closed central d-subgroup consisting
precisely of the semisimple elements contained in the center of $\cA^\circ$,
and  $\cA_1$ is a $\BQ$-closed normal subgroup of $\cA$ 
with unipotent by-finite center. 

Next we define 
$$Z_1 =C\cap {\cal Z}_1 \; , \; C_1 =C\cap \cA_1 \; , \;  C_2 = Z_1  \cdot  C_1 \; .$$
Observe that $Z_1$ is an arithmetic subgroup of ${\cal Z}_1$ and that 
$C_1$ is arithmetic in $\cA_1$. It follows that $C_2$ is
arithmetic in $\cA$ and of finite index in $C$.  
Since $C$ is torsion-free, we have $Z_1 \cap C_1=\{1\}$,  and 
therefore $C_2$ is isomorphic to the direct product   
$Z_1 \times C_1$. The action of $G_1$ as an $\cA$-algebraic group 
of automorphisms of $C$
stabilizes the factors $Z_1$ and $C_1$, and hence also $C_2$. 

Now put $\cD$ for the Zariski-closure of $D$.  Since $D$ is
central in $A$, $\cD \leq \Z(\cA)$. It follows that the maximal
$d$-subgroup $\cS_\cD$ of $\cD$ is contained in $\cZ_1$.
Since $\cS_\cD$ is invariant in $\cD$, there exists, by virtue of AG6,
an almost direct product decomposition $\cZ_1 = \cS_\cD \cdot \cS_2$
which is respected by $G_1$. We define $Z_\cD = Z_1 \cap \cS_\cD $
and $Z_2 = Z_1 \cap \cS_2$. By the arithmeticity of the factors
$Z_1$ and $Z_2$, the product $Z_\cD\cdot  Z_2$ is of finite index 
in $Z_1$. Also this decomposition is preserved by $G_1$. 

Now define $C_3= Z_\cD \cdot Z_2 \cdot C_1 \leq \cS_\cD \cdot \cS_2 \cdot \cA_1$.
This group is an arithmetic subgroup, $G_1$-invariant and of finite index in $A$. 
Since it is torsion-free
it is also a direct product of its factors. Let
us put $\cA_2 =  \cS_2 \cdot \cA_1$  and
 $\bar{\cA_2} = \cA_2 \big/  (\cA_2 \cap \cS_\cD)$.  
Then we have an induced arithmetic and Zariski-dense embedding  
$$    C_3 =  Z_\cD \times  \left(Z_2 \times  C_1\right) \, \leq \, \cS_\cD \times \bar{\cA_2}  \; . $$
Since $Z_\cD$ is isomorphic 
to $\BZ^n$, for some $n \geq 0$,  we may embed this group as an
arithmetic subgroup into $\Ga^n$. This gives rise to an arithmetic and
Zariski-dense embedding
$$  C_3 =   Z_\cD \times  \left(Z_2 \times  C_1\right)  \, \leq \, \Ga^n \times \bar{\cA_2}  \;  $$
which has the property that  $D \cap C_3$ is unipotent by-finite in 
$\Ga^n\times \bar{\cA_2}$.
 
Note that $G_1$ induces 
an $\bar{\cA_2}$-algebraic group of automorphisms of $C_1$.
We consider now the arithmetic embedding of $Z_\cD$ 
into the unipotent group $\Ga^n$.  The action of $G_1$ 
on $Z_\cD$ extends to an action by $\bbQ$-defined 
automorphisms of $\Ga^n$. This turns $G_1$ into
a $\Ga^n$-algebraic group of automorphisms of $Z_\cD$.
We infer from Lemma \ref{algformalproducts} that the
product action of $G_1$ on $\Ga^n\times \bar{\cA}_2$ turns
$G_1$ into an algebraic group of automorphisms of $\Ga^n\times \bar{\cA}_2$.
This, in particular, turns $G_1$ into an $\Ga^n\times \bar{\cA}_2$-algebraic group 
of automorphisms of $C_2$. 

Now let $\rho$ be an $\cA$-rational representation of $A$ which contains
$D$ in its kernel. In particular, it satisfies $\cD \leq \ker \rho_\cA$, where
$\rho_\cA$ denotes the extension of $\rho$ to $\cA$. The restriction of 
$\rho_\cA$ to $\cA_2$ gives rise to a $\bbQ$-defined homomorphism
$\bar{\rho}: \Ga^n \times \bar{\cA_2} \ra \cG$ which has the subgroup $\Ga^n$ in its kernel. 
We contend that $\bar{\rho}$ extends the representation $\rho$ on $C_3$.
This is easily verified. Thus $\rho: C_3 \ra \cG$ is  $\Ga^n \times \bar{\cA_2}$-rational. 

Since $C_3$ is normal and of finite index in $A$, this 
allows, by application of Lemma \ref{crit},  to embed the
group $A$ as an arithmetic subgroup in a finite extension ${\cE}$ of 
$\Ga^n\times \bar{\cA}_2$ such that $G$ acts as a ${\cal E}$-algebraic group
of automorphisms of $A$.  This embedding has the property that 
the finite index subgroup $C_3 \cap D \leq D$ is unipotent. Hence, $D$ is
unipotent-by-finite under the embedding of $A$ into ${\cal E}$. 
This proves i) and ii). The last statement of Lemma \ref{crit}
asserts that $\rho$ is $\cE$-rational, since the restriction of $\rho$
to $C_3$ is $\Ga^n\times \bar{\cA}_2$-rational.
 \end{proof}

Part of Proposition \ref{changerep} is reminiscent of Corollary 3.5 from \cite{GP1}
and of Proposition 3.3 from \cite{GP3}, but it is stronger since no passages
to subgroups of finite index are required. Our ultimate arithmeticity result
is contained in the next proposition.

\begin{proposition}\label{endprop}
Let $\cA$ be a  $\bbQ$-defined linear algebraic group and let 
$A\le \cA(\bbQ)$ be a Zariski-dense subgroup. 
Assume $N$,$B$,$C$ are normal subgroups of $A$ such that
the following hold: \begin{itemize} 

\item[i)] $N\cdot B$ has finite index in $A$,

\item[ii)] $B$ is an arithmetic subgroup in its Zariski-closure ${\cal B}$, 

\item[iii)] $C$ is an arithmetic subgroup in its Zariski-closure ${\cal C}$,

\item[iv)] $C\le N\cap B$ and $D=  (N\cap B)\big/C$ is in the center of $B/C$.
\end{itemize}
Then $A/N$ is an arithmetic group.

Moreover, there exists an arithmetic embedding of  $A/N$ into a $\bbQ$-defined 
linear algebraic group $\cA_N$ which has the following property: For any
$\cA$-rational representation of $A$ with $N \leq \ker \rho$, the induced
quotient representation of $A/N$ is  an  $\cA_N$-rational representation. 
\end{proposition}
\begin{proof} The group $A$ induces by conjugation an ${\cB}$-algebraic
group of automorphisms of $B$ which we call $G$. 
Note that ${\cal B}$ and ${\cal C}$ are normal in $\cA$, and preserved 
by $G$ as well.  Since $B$ and $C$ are 
arithmetic subgroups of their respective Zariski-closures, we find
that $C$ has finite index in ${\cal C}\cap B$. We may hence use Proposition
\ref{propquot} to embed the group $B/C$ as an arithmetic subgroup into a 
$\bbQ$-defined linear algebraic group ${\cal D}$. This embedding has the
property that the group $G$ of automorphisms of $B/C$ is ${\cal D}$-algebraic.

We consider now the subgroup $D= (N\cap B)\big/C$ in $B/C$. 
By our assumption iv),  $D$ is central in $B/C$. Since $G$ is ${\cal D}$-algebraic,
we may, by Proposition \ref{changerep}, change the arithmetic 
embedding of $B/C$ in $\cD$ to an arithmetic embedding of $B/C$ into  
a $\bbQ$-defined linear algebraic group ${\cal E}$ such that 
$D$ is unipotent-by-finite in ${\cal E}$. Moreover, 
the group $G$ acts as an ${\cal E}$-algebraic group of
automorphisms. 

Let us  consider now the Zariski-closure ${\cal D}_1$ of $D$ in ${\cal E}$.
Since $D$ has a unipotent finite index subgroup, 
$D$ is an arithmetic subgroup of the 
$\bbQ$-defined algebraic group ${\cal D}_1$. (For a proof consult 
\cite[Chapter 8]{Segal}).
Since $D$ is an arithmetic subgroup it has finite index in ${\cal D}_1 \cap B\big/C$.
Therefore, we may apply Proposition
\ref{propquot} to embed the quotient $(B/C)\big/D$ as an arithmetic
subgroup into a $\bbQ$-defined linear algebraic group ${\cal B}_1$ such that $G$ acts by
${\cal B}_1$-rational automorphisms on $(B/C)\big/D$.

By assumption i), $A/N$ is isomorphic to a finite extension group of 
$$(B/C)\big/D\cong B\big/B\cap N\cong \left(N\cdot B\right) \big/N.$$ 
The elements of $A/N$ act on $(B/C)\big/D$ as ${\cal B}_1$-rational automorphisms
since the elements of $G$ have this property.   
Finally,  we apply Lemma
\ref{crit} to find that $A/N$ is arithmetic 
in a $\bbQ$-defined linear algebraic group .  

To prove that the restriction of $\rho$ to $A/N$ is 
$\cA_N$-rational, we have to carry over the rationality of $\rho$ in
each of the construction steps above. The details are easily verified. 
\end{proof}


\section{The arithmeticity of $\Out(\Gamma)$} \label{secari}


This section contains the complete proof of Theorem
\ref{teoa}  which proceeds in two steps. These steps are carried out
in Section \ref{asec1} and in Section \ref{generalpofi}. On our way,
we provide (respectively, finish) the proof of  Theorem \ref{ftop} in 
Section  \ref{asec1}, as well as the proofs of
Theorem \ref{teoc} and of Theorem \ref{teosep} 
in Section \ref{generalpofi}.

Throughout this section, $\Gamma$ denotes a polycyclic-by-finite
group. We also stick to the notation introduced in Sections 2 to
6.  In particular,  $F\le \Gamma$ denotes the Fitting subgroup of $\Gamma$. 
If in addition $\Gamma$ is a wfn-group,  $\bH_\Gamma$ denotes the
algebraic hull of $\Gamma$,
and $\bF$ the Zariski-closure of $F$. 


\subsection{The case of polycyclic-by-finite wfn-groups} \label{asec1}


The purpose
of this subsection is to prove Theorem \ref{ftop} of the introduction. 
Let us therefore assume here that $\Gamma$ is a wfn-group. 
The arithmeticity of  $\Out(\Gamma)$, 
in the case that $\Gamma$ is a wfn-group,  
is an immediate consequence of the
structural properties of the embedding $\Aut(\Gamma) \leq \Aut_a({\bf H_\Gamma})(\bbQ)$ 
together with Proposition \ref{endprop}. To see this let us put now 
\begin{equation}\label{festsetz}
A=\Aut(\Gamma) \, ,\quad N=\Inn_\Gamma\, ,\quad B=\AGF \, , \quad C=\Inn_F \, .
\end{equation}
Then we have: 

\begin{proposition}  Let $\Gamma$ be a wfn-group. 
Then the subgroups 
$N,B,C\le A$ defined in \eqref{festsetz} 
satisfy the hypotheses of Proposition \ref{endprop}
with respect to the Zariski-closure $\cA$
of $\Aut(\Gamma)$ in $\Aut_a({\bf H_\Gamma})$.
\end{proposition} 

\begin{proof}
Condition i) requires that  $N\cdot B = \Inn_\Gamma \cdot \AGF$
has finite index in $A = \Aut(\Gamma)$. This is contained in  
Proposition \ref{InnGAGFfi}.  
Theorem \ref{teoagf} says that $B=\AGF$ 
is arithmetic in its Zariski-closure ${\cal B}$ in $\cA$.
This implies condition ii). 

The construction of the algebraic structure on $\Aut_a({\bf H_\Gamma})$ 
(see Subsection \ref{subsect:algebraicstructure}) shows that 
the group $\Inn_\bF$ is a Zariski-closed subgroup 
of the unipotent radical of 
$\Aut_a({\bf H_\Gamma})$. Moreover, $\Inn_\bF$
contains the finitely generated group $\Inn_F$ as a Zariski-dense 
subgroup of rational points. 
In particular,  $C=\Inn_F$ is arithmetic in its Zariski-closure.
This implies condition iii).
 
We shall finally verify the conditions iv) of Proposition \ref{endprop}.
Clearly we have $\Inn_F \leq \Inn_\Gamma \cap \AGF$, that is, 
$C\le N\cap B$.  
Now let $\Phi \in \AGF$ and $\gamma\in\Gamma$ with $\Inn_\gamma \in \AGF$. Then
$\Phi(\gamma) = \gamma f$, where $f \in F$. It follows that 
$$\Phi\circ \Inn_\gamma\circ \Phi^{-1} = 
\Inn_{\Phi(\gamma)}= \Inn_\gamma \, \Inn_f.$$
This shows that $\Inn_\Gamma \cap \AGF = N\cap B$ projects onto
a central subgroup of $B/C=\AGF /\, \Inn_F^\Gamma$. Hence, iv) holds.
\end{proof}

\begin{prf} {\it Proof of Theorem \ref{ftop}.}
We may now apply Proposition \ref{endprop} which asserts that 
there exists a $\BQ$-defined linear algebraic group
\begin{equation*} 
{\cal O}_\Gamma =\cA_{N}=\cA_{\Inn_\Gamma}  
\end{equation*}
which contains an isomorphic copy of the group 
$\Out(\Gamma)=A/N$ as an arithmetic subgroup. 
This already establishes the arithmeticity of $\Out(\Gamma)$.

Consider next the algebraic outer automorphism group  
$$ \Out_a(\bH_\Gamma)=    \Aut_a(\bH_\Gamma)   \big/\, \Inn_{\bH_\Gamma} \, , $$
and let $\pi_\Gamma:  \Out(\Gamma)\ra \Out_a(\bH_\Gamma)$
be the homomorphism induced on $\Out(\Gamma)$.
Since the natural map 
$ \Aut_a(\bH_\Gamma) \ra   
\Out_a(\bH_\Gamma)$
is a $\BQ$-defined homomorphism, it
induces a $\bbQ$-defined homomorphism $\cA \ra \Out_a(\bH_\Gamma)$.
Since $\pi_\Gamma$ contains $N = \Inn_\Gamma$ in its kernel,  
Proposition \ref{endprop} asserts that the
homomorphism $\pi_\Gamma :\Out(\Gamma)\to \Out_a(\bH_\Gamma)$ 
can be extended
to a $\BQ$-defined homomorphism 
$$ \pi_{\cO_\Gamma} : \, {\cal O}_\Gamma\to \Out_a(\bH_\Gamma) \; .$$
This proves the first part of Theorem \ref{ftop}.

To show the statements about the kernel of $\pi_{\Gamma}$, we define:
\begin{equation*}
K=\Inn_\bH\cap\Aut(\Gamma),\qquad K_F=\Inn_\bH\cap \AGF,\qquad
E_F = \Inn_{\bf F}^{\bH_\Gamma}   \cap \AGF   \; . 
\end{equation*}

\begin{lemma}\label{kerra}
With the above notation the following hold: \begin{itemize}

\item[i)] $\Inn^\Gamma_F$ has finite index in $E_F$.

\item[ii)] There is a finite index normal subgroup $T\le \AGF$ 
such that $T\cap E_F\le \Inn^\Gamma_F$.

\item[iii)]  The commutator group $[\AGF, K_F]$ is contained in $E_F$. 
\end{itemize}
\end{lemma}
\begin{proof} Note that
$\Inn_F^\bH$ is Zariski-dense and arithmetic 
in the unipotent group $\Inn_{\bf F}^\bH$.
Since $\AGF$ is an arithmetic subgroup of
its Zariski-closure in $\cA$, $E_F$ is arithmetic
in $\Inn_{\bf F}^\bH$ as well. This implies i). 

Now ii) follows from i) together with the congruence 
subgroup property for $\Inn^\Gamma_F$.
(Compare the proof of Lemma \ref{lemmapropquot}.)

For iii),  let $\Phi$ be the extension of 
$\phi\in \AGF$ to an automorphism of $\bH_\Gamma$.  
Let $h\in\bH$ such that $\Inn_h=\psi \in K_F$. 
Since $\Phi \in  \cA_{\H_\Gamma| \bF}$
$$\Phi\circ \Inn_h\circ\Phi^{-1}=\Inn_{\Phi(h)}=\Inn_h\circ \Inn_{f_h}.$$ 
This in turn gives $\psi^{-1}\circ\phi\circ\psi\circ \phi^{-1}\in
\Inn_{\bf F}^{\bH_\Gamma}  \cap \AGF = E_F$, proving iii).  
\end{proof}

The kernel of $\pi_\Gamma$ is the image of $K= \Inn_{\bH_\Gamma}  \cap \Aut(\Gamma)$
in $\Out(\Gamma)$. Let $\bar K_F$ and $\bar E_F$ be the images of $K_F$, $E_F$ in 
$\AGF\big/\, \Inn_F^\Gamma$. Let $\bar T$ be the corresponding image of the finite index 
subgroup $T\le\AGF$ as in Lemma \ref{kerra} ii).  
Lemma \ref{kerra} i) shows that $\bar E_F$ is finite. By ii) and  iii) of the same 
Lemma, $\bar K_F$ is centralized by the finite index subgroup $\bar{T}
\leq \AGF\big/\, \Inn_F^\Gamma$. In particular, 
$\bar K_F$ is abelian-by-finite. Consider now the commutative diagram
$$
\xymatrix{
\AGF\big/\, \Inn_F^\Gamma \ar[r]^{ } \ar[d]_{ } 
& \Out(\Gamma) = \Aut(\Gamma)\big/\, \Inn_\Gamma \ar[ld]^{ } \\
\Out_a({\bf H}_\Gamma) & 
}
$$
of natural homomorphisms.  By AR3,
every abelian subgroup of an arithmetic group is finitely generated. 
Hence, the image of $\bar K_F$ in $\Out(\Gamma)$ is so. 
We may also infer that $\bar K_F$ is finitely generated.

Since $\Inn_\Gamma \cdot \AGF$ has finite index in $\Aut(\Gamma)$, 
the normal subgroup $\bar K_F$ maps onto a finite index subgroup of the image of $K$ in 
$\Aut(\Gamma)\big/\, \Inn_\Gamma$. This proves that $\ker \pi_{\Gamma}$
is finitely generated, abelian-by-finite and centralized by a finite index 
subgroup of $\Out(\Gamma)$.  If $\Gamma$ is nilpotent-by-finite then $E_F$ is of finite
index in  $K_F$, and hence $\ker \pi_{\Gamma}$ is finite. 
This finishes the proof of Theorem \ref{ftop}. 
\end{prf}


\subsection{The case of a general polycyclic-by-finite group}
\label{generalpofi}


In this subsection, we explain the transfer of our arithmeticity results from 
the case of wfn-groups to general polycyclic-by-finite groups. Thereby, 
we provide the final step in the proofs of Theorem \ref{teoa} and Theorem  \ref{teoc}. 
We also prove Proposition \ref{equivi} and Theorem \ref{teosep}.  

Let $\Gamma$ be a polycyclic-by-finite group. Let $\tau_\Gamma$ denote the maximal finite normal 
subgroup of $\Gamma$.  Note that $\tau_\Gamma$ is
characteristic in $\Gamma$ that is,
it is normalized by every automorphism of $\Gamma$.
The quotient group $$ \tGam := \, \Gamma\big/ \tau_\Gamma$$
is a wfn-group. We let ${\mathsf j}: \Gamma \ra \tGam$ denote the
quotient homomorphism. By Theorem \ref{ftop}, the group  
$\Out(\tGam)$ is an arithmetic group.  We 
shall show that $\Out(\Gamma)$ has the same property. 

Let $\Gamma_{0} \leq \Gamma$ be a characteristic
finite index subgroup with 
$\Gamma_{0} \cap \tau_\Gamma = \{ 1\}$. We may suppose that 
the image $\tGam_0 \leq \tGam$ of $\Gamma_0$ in
$\tGam$ is also characteristic. (To obtain such a subgroup,  let $n\in\BN$ be the
index of a torsion-free subgroup of finite index in $\Gamma$. Take $\Gamma_0$ to
be the subgroup generated by all $\gamma^{n!}$, $\gamma\in\Gamma$). 
 
Let us put  $$ \mu: = \, \Gamma/ \Gamma_{0} \; . $$
The quotient homomorphism ${\mathsf j}$ and the projection to $\mu$
induce injective homomorphisms  
$$  {\mathsf j}_\mu: \, \Gamma \rightarrow \tGam \times \mu  \, , \, \; \text{  and } \; \;
 {\mathsf k}_\mu: \, \Aut(\Gamma)  \rightarrow \Aut(\tGam) \times\Aut( \mu) \; . $$
 
We define finite index subgroups of $\Aut(\Gamma)$ and $\Aut(\tGam)$,
respectively:
\begin{equation*}
A_{0}: =  \;  \{\, \phi \in \Aut(\Gamma) \mid 
\phi_{\Gamma/ \Gamma_{0}} = \id_{\Gamma/ \Gamma_{0}}\,  \} \;  \leq \Aut(\Gamma) \; , 
\end{equation*}
\begin{equation*}
\tilde{A}_{0}: = \; \{\, \phi \in \Aut(\tGam) \mid 
\phi_{\tGam/ \tGam_{0}} = \id_{\tGam/ \tGam_{0}}\, \} \; \leq \Aut(\tGam) \; . 
\end{equation*}

\begin{lemma} 
\label{lemma:findex} 
With the above notation the following hold: 
\begin{itemize} 
\item[i)]  The group ${\mathsf j}_\mu(\Gamma)$ is of finite index in $\tGam \times \mu$.
\item[ii)]  Let $F = \Fitt(\Gamma)$. Then ${\mathsf j}(F)$  is of finite index in $\Fitt(\tGam)$.
\item[iii)]   The induced homomorphism ${\mathsf k}: \Aut(\Gamma) \ra \Aut(\tGam)$ maps the
group $A_0$ isomorphically onto $\tilde{A}_{0}$.  In particular, ${\mathsf k}(\Aut(\Gamma))$
is of finite index in $\Aut(\tGam)$.
\item[iv)]   The subgroup $\AGF \leq \Aut(\Gamma)$ is mapped by ${\mathsf k}$ onto
a finite index subgroup of $\tAGF$.  
\end{itemize} 
\end{lemma}
\begin{proof} Part i) is clear. For ii), note first that ${\mathsf j}(F)$
is a nilpotent ideal in $\tGam$, and hence ${\mathsf j}(F) \leq \Fitt(\tGam)$. 
Then $F_0 = {\mathsf j}^{-1}(\Fitt(\tGam) \cap \Gamma_0)$ is a nilpotent 
normal subgroup of $\Gamma$ which is of finite
index in the preimage ${\mathsf j}^{-1}(\Fitt(\tGam))$. Moreover, 
$F_0 \leq F$. Hence $F$ is of finite index in this preimage.
This implies ii). 

To prove iii), we show that 
${\mathsf k}_\mu(A_{0}) =  \tilde{A}_{0} \times \{ 1 \} $.
Clearly,  ${\mathsf k}_\mu(A_{0})$ is contained in $\tilde{A}_{0} \times \{ 1 \}$. 
Let $\psi \in  \tilde{A}_{0} $. We show that there exists
$\phi \in A_{0}$ such that  $\psi =\mathsf{k}({\phi})$ is induced by $\phi$.  
If $\gamma \in \Gamma$, we let $\tilde{\gamma} = \mathsf{j}(\gamma)$
denote its projection into $\tilde{\Gamma}$.
Since the projection ${\mathsf j}$ maps $\Gamma_0$ isomorphically
onto the invariant subgroup $\tGam_0$, the automorphism 
$\psi$ defines $\psi_0 \in \Aut(\Gamma_0)$ uniquely with the property that 
\[   \mathsf{j} (\phi_0( \gamma) )   =   \psi( \tilde{\gamma})   \;  \; \;  (\gamma \in\Gamma_0 ) \; . \]
Let 
$\Gamma\big / \Gamma_0 = \bigcup_i \gamma_i \Gamma_0$ be the
coset decomposition.
Now we can write 
$\psi(\tilde{\gamma}_i) =  \tilde{\gamma}_i  \epsilon_i$, $\epsilon_i \in \tGam_0$.
There exist unique $ \delta_i \in \Gamma_0$ such that  $\epsilon_i = \tilde{\delta}_i$. 
We declare now 
$$  \phi(\gamma_i s) = \tilde{\gamma}_i \delta_i \phi_0(s)    \;  \; \;  (s \in\Gamma_0 )  \; . $$
It is easy to verify that  this actually defines 
an automorphism $\phi \in \Aut(\Gamma)$.  This
$\phi$ is clearly a lift of  $\psi$.

For iv) consider the quotient homomorphism $\Gamma/ F \ra \tGam/\tF$. 
By ii), this homomorphism has finite kernel. Since it is surjective, the group
$\AGF$ is mapped into $\tAGF$ by ${\mathsf k}$. Let $A_1$ be the
preimage of $\tAGF$ in $\Aut(\Gamma)$. Now applying the
reasoning of iii) to the above quotient homomorphism with finite kernel, 
we deduce that finite index subgroup of $A_1$ acts as the identity on $\tGam/\tF$.
That is, $\AGF$ has finite index in $A_1$. This shows iv). 
\end{proof}

Part iii) of the above lemma immediately implies Proposition \ref{equivi}:
 
\begin{prf} {\it Proof of Proposition \ref{equivi}:}
The groups $\Aut(\Gamma)$ and $\Aut(\tGam)$ are 
commensurable, since they have isomorphic finite index
subgroups $A_0$ and $\tilde{A}_{0}$.  
If  $\Aut(\tGam)$ is arithmetic, then the product  
$\Aut(\tGam) \times\Aut( \mu)$ is arithmetic. 
Since $\Aut(\Gamma)$ embeds as a subgroup of 
finite index in the latter product,  
$\Aut(\Gamma)$  is an arithmetic group as well.  
Conversely, if $\Aut(\Gamma)$ is arithmetic, 
the finite index subgroup  $A_0 \leq \Aut(\Gamma)$ is arithmetic too. 
Therefore, the subgroup $\tilde{A}_0 \leq \Aut(\tGam)$ is arithmetic.
\end{prf}

A related result is: 
\begin{proposition} The subgroup $\AGF$ of $\Aut(\Gamma)$ is arithmetic. \label{AGFisa}
\end{proposition} 
\begin{proof} By iv) of  Lemma \ref{lemma:findex}, the injection ${\mathsf j}_\mu$ 
maps $\AGF$ onto a finite index subgroup of $\tAGF \times \Aut(\mu)$. 
Since $\tAGF$ is arithmetic, by Corallary \ref{AGFisa1}, we can infer that 
$\AGF$ is arithmetic.  
\end{proof}

\begin{prf} {\it Proof of Theorem \ref{teoc} in the general case.}
As remarked above, the projection $\mathsf k$ maps $\AGF$ onto a finite
index subgroup of $\tAGF$.  Let $B \leq \Gamma \cap \Gamma_0$ be a 
nilpotent subgroup such that $\tAGF \cdot {\mathsf j}(B)$ is of finite
index in $\Aut(\tGam)$ (see Section \ref{secautoauto}). Then $\AGF \cdot B$ is of finite index 
in $\Aut(\Gamma)$.  Together with Proposition \ref{AGFisa}, this
proves the required decomposition of $\Aut(\Gamma)$. 
\end{prf}

As another consequence of Lemma \ref{lemma:findex}, we 
infer that $\Out(\Gamma)$ and $\Out(\tGam)$ are S-commensurable: 

\begin{prf}{Proof of Proposition \ref{finiin}.}
Since ${\mathsf k}(\Inn_{\Gamma_0}) \leq \Inn_\tGam$ is
of finite index in $\Inn_\tGam$ and ${\mathsf k}: A_0 \ra \tilde{A}_0$ is an isomorphism, 
$\A_0 \cap \Inn_\Gamma$ is a subgroup of finite index in  ${\mathsf k}^{-1}(\tA_{0} \cap \Inn_\tGam)$. 
We consider the map on quotients 
$$   A_{0} \big/ A_{0} \cap \Inn_\Gamma 
\xrightarrow{{\mathsf k}^*} \tA_{0 }\big/ \tA_{0} \cap \Inn_\tGam$$
which is induced by $\mathsf k$.
The above implies that ${\mathsf k}^*$ has finite kernel. 
Since the left hand side 
is a finite index subgroup of $\Out(\Gamma)$, the corollary follows.
\end{prf} 

\begin{prf} {\it Proof of Theorem \ref{teoa} in the general case.}
Let $\Gamma$ be a polycyclic-by-finite group.  Then we know that $\Out(\Gamma)$
is residually finite (by \cite{Wehrfritz}). By Corollary
\ref{finiin},  it projects with finite kernel onto a finite index subgroup of
the arithmetic group $\Out(\tGam)$. Thus Proposition \ref{finker} 
implies that $\Out(\Gamma)$ is arithmetic.  
\end{prf}
   
\begin{prf} {\it Proof of Theorem \ref{teosep}.}
By Proposition \ref{equivi} we reduce to the case of $\tilde\Gamma$.
Using Proposition \ref{InnGAGFfi} 
we infer that  
$\A_{\tilde\Gamma|F}$ is of finite index in $\Aut(\tilde\Gamma)$.
We finally use Theorem \ref{teoagf}.
\end{prf}


\section{Polycyclic groups with non-arithmetic automorphism groups}


We present examples of polycyclic groups 
whose automorphism groups are not isomorphic 
to any arithmetic group. In particular, we shall prove Theorem \ref{teob}.


\subsection{Automorphism groups of semi-direct products}\label{splitext}


Here are some remarks concerning the automorphism group of groups $\Gamma$
which are  semi-direct products $F\rfish D$ where 
$D$ is a group and $F$ is a (commutative) $D$-module. 

We write the group product in $F$ additively, and 
for $h \in D$, we write $f \mapsto h \cdot f$, $f \in F$, 
to denote the action of the element $h$ on $F$.
Let $\Theta$  be a subgroup of $\Gamma$. We write ${\rm Inn}_\Theta$ for the
subgroup of $\Aut(\Gamma)$ consisting of the inner automorphisms defined by
the elements of $\Theta$.
Similarly as before,  we put
$$\AGF:=\{\, \phi \in \Aut(\Gamma)\Mid \phi(F)=F,\, 
\phi |_{\Gamma /F}=  \id_{\Gamma /F} \,\}.$$
There are two constructions for auto\-morphisms in $\AGF$. For the first, let 
${\rm Der}(D,F) = \{ d: D \to F \mid d(h_1 h_2) = z(h_1)+ h_1 \cdot z(h_2) \}$ 
be the group of derivations from $D$ into $F$. The group of
derivations naturally obtains the structure of a $D$-module by setting
$$g*d\; (h):=g\cdot d(g^{-1}hg)\qquad (g,\, h\in D,\, d\in {\rm Der}(D,F)).$$
A derivation $d\in {\rm Der}(D,F)$ gives rise to 
an automorphism $\phi_d :\Gamma\to \Gamma$ by 
$$\phi_d((m,g)):=(m+d(g),g)\qquad\qquad (m\in F,\, g\in D) \; .$$
We write $\Aut^{\rm d}(\Gamma)$ for the (abelian) subgroup of $\AGF$
consisting of these automorphisms. We remark that the homomorphism 
$${\rm Der}(D,F)\to \AGF,\qquad d\mapsto \phi_d \quad (d\in {\rm
  Der}(D,F))$$
is $D$-equivariant with respect to the above $D$-action on ${\rm Der}(D,F)$ 
and conjugation by elements of  $\Inn_D$ on $\AGF$.

Let us define $\Aut_D(F)$ to be the group of $D$-equivariant automorphisms of $F$.
Given a $D$-equivariant automorphism $\rho :F\to F$, we define
an automorphism $\phi_\rho :\Gamma\to \Gamma$ by
$$\phi_\rho((m,g)):=(\rho(m),g)\qquad\qquad (m\in F,\, g\in D) \; .$$
We write $\Aut^{\rm a}(\Gamma)$ for the subgroup of $\AGF$
consisting of these automorphisms.

\begin{proposition}\label{splita} Let $\Gamma=F\rfish D$ be a semi-direct
product of a abelian $D$-module $F$ by the group $D$. We then have:
\begin{itemize}
\item[i)]  $\AGF=\Aut^{\rm d}(\Gamma)\cdot \Aut^{\rm a}(\Gamma)$.

\item[ii)] ${\rm Inn}_\Gamma\cdot\AGF=
\left(\Aut^{\rm d}(\Gamma)\cdot \Aut^{\rm a}(\Gamma)\right)\cdot {\rm Inn}_D$.

\item[iii)] $\Aut^{\rm a}(\Gamma)$ centralizes $\Inn_D$.

\item[iv)] $\Aut^{\rm d}(\Gamma)\cap \Aut^{\rm a}(\Gamma)=\{1\}$.

\item[v)] $\Aut^{\rm d}(\Gamma)$ is an abelian normal subgroup in 
${\rm Inn}_\Gamma\cdot \AGF$.  
\end{itemize}
\end{proposition}

The proof of this proposition is straightforward, we skip it.


\subsection{Examples}\label{exnona}


In order to show that certain groups are not arithmetic we use the following
simple criterion.

\begin{proposition}\label{splitprop}
For a matrix $A\in\GL(n,\BZ)$, let $\Gamma(A)=\BZ^n\rfish \langle
A \rangle$ be the split extension of $\BZ^n$ by the cyclic group generated by
$A$. If $\Gamma(A)$ is an arithmetic group then either $A$ is of finite order
or a power of $A$ is unipotent or $A$ is semisimple.
\end{proposition}

In the first two cases of Proposition \ref{splitprop} that is if  
either $A$ is of finite order or a power of $A$ is unipotent the group 
$\Gamma(A)$ is arithmetic. In case $A$ is semisimple $\Gamma(A)$ can be
arithmetic but examples in \cite{GP4} show that it need not have this
property. 

\begin{proof}
Suppose that $\Gamma(A)$ is an arithmetic group and that $A$ is 
not of finite order nor a power of $A$ is unipotent. In this case, we have
$\Fitt(\Gamma(A))=\BZ^n$. Assume further that $\Gamma(A)$ is an arithmetic group.
We can find (compare \cite[Theorem 3.4]{GP1})
a solvable $\BQ$-defined linear
algebraic group ${\bf H}$, having a strong unipotent radical 
so that there is an isomorphism 
$\psi: \Gamma(A)\to \Gamma$
where $\Gamma$ is a Zariski-dense arithmetic
subgroup of ${\bf H}(\BQ)$. Let $\bU$ be the unipotent radical of ${\bf H}$ and
let $\lu$ denote the Lie-algebra of $\bU$. 
The exponential map
$ \exp:  \lu \rightarrow  \bU $ 
is a $\bbQ$-defined isomorphism of varieties.
The adjoint representation leads to a $\BQ$-defined rational
representation $\alpha_{\bf H}: {\bf H}\to \Aut({\lu})$ which is defined by 
$$\alpha_{\bf H}(g)\; (x)=\; \exp^{-1}(g\exp(x)g^{-1})\quad (g\in {\bf H},\, 
x\in {\lu}).$$
Since $\alpha_\bH$
is $\bbQ$-defined,  the image $\alpha_{\bf H}(\Gamma)$ is a Zariski-dense and
arithmetic (by AR1) subgroup of $\alpha_{\bf H}({\bf H})$.
The kernel of $\alpha_{\bf H}$ is equal to $\bU$. 
Taking the image under $\exp^{-1}$ of the standard basis in $\BZ^n$ we obtain a
$\BQ$-basis of $\lu(\BQ)$. Expressed in this basis $\psi(A)$ acts
by the matrix $A$ on the subspace $\lu_0$ spanned by these elements. 
Also ${\bf H}_0$ stabilizes $\lu_0$.
Hence the cyclic subgroup generated by $A$ is arithmetic and Zariski-dense in
$\bH_1 = \alpha_{\bf H}(\bH)$.

Let $A=SJ$ be the Jordan-decomposition of $A$, that is $J\in {\bf H}_1(\BQ)$ is
unipotent and $S\in {\bf H}_1(\BQ)$ is semisimple and $JS=SJ$ holds. 
There is a $n\in\BN$
such that $J^n\in {\bf H}_1(\BZ)$,  and hence a $m\in\BN$
such that $J^m\in \langle A\rangle$. This implies 
$S^m\in \langle A\rangle$. We infer that $J=1$. 
\end{proof}
   
We shall discuss now the example from the introduction. 
That is, we choose $d\in\BN$ not a square, 
 set $\omega=\sqrt d$ and
let $K=\BQ(\omega)$ be the
corresponding real quadratic number field. We write $x\mapsto \bar x$ for the
non-trivial element of the Galois group  of $K$ over $\BQ$. 
We consider the subring ${\cal O}=\BZ+\BZ\omega\subset K$ and choose 
a unit $\epsilon=a+b\omega$ of ${\cal O}$ which is of infinite
order and satisfies $\epsilon \bar\epsilon=1$. 

Let  $D_\infty$ be the infinite dihedral group as in (\ref{dihe}).
We further take $F={\cal O}\times \BZ$ with the $D_\infty$-module structure
defined as in (\ref{modu}).
As done in the introduction, we put 
$$\Gamma(\epsilon):=F\rfish D_\infty.$$

We describe four derivations $d_1,\ldots, d_4$
in  ${\rm Der}(D_\infty,F)$ by specifying their values on 
the generators $A,\, \tau$ of $D_\infty$. 
We define $l$ to be the greatest common factor of $a+1$ and $bd$.
Now put:
$$d_1(A)=(0,1),\ d_1(\tau)=(0,0)\, ; \quad d_2(A)=(0,0),\ d_2(\tau)=(0,1) \, ;$$
$$d_3(A)=(\omega,0),\ d_3(\tau)=(\omega,0) \, ;
\quad d_4(A)=\left(\frac{(\epsilon+1)\omega}{l},0\right),\ d_4(\tau)=(0,0).$$
Each of the above pairs of values defines a
derivation by extension. 

We also define
\begin{equation}\label{mamat}
\hat A:=\left(\begin{array}{cccc}
1 & -2 & 0 & 0 \\
0 & 1  & 0 & 0 \\
0 & 0 & -1 & \frac{-2(a+1)}{l} \\
0 & 0 &  g & 2a+1 
\end{array} \right).
\end{equation}
The structure of $\Aut\big(\Gamma(\epsilon)\big)$ is described in the following
proposition. 

\begin{proposition}\label{specialgr} The following hold in
$\Aut(\Gamma(\epsilon))$:  \begin{itemize}
\item[i)] $\Aut_{D_\infty}(F)$ is finite.

\item[ii)] The derivations $d_1,\ldots, d_4$ are a $\BZ$-basis of 
${\rm Der}(D_\infty,F)$.

\item[iii)] The action of\/ $\Inn_A$ on ${\rm Der}(D_\infty,F)$ 
expressed relative  to the basis  $d_1,\ldots, d_4$ is given by the matrix $\hat A$.

\item[iv)] $\Aut\big(\Gamma(\epsilon)\big)$ contains a subgroup of finite index which is
   isomorphic to $\Gamma(\hat A)$.
\end{itemize}   
\end{proposition}

\begin{proof} Items i), ii), iii) are proved by some 
straightforward computations
which we skip. Remark that $F$ is the Fitting-subgroup
of $\Gamma(\epsilon)$. Setting $\Gamma=\Gamma(\epsilon)$ 
we know from Proposition \ref{InnGAGFfi} 
that ${\Inn}_\Gamma\cdot \AGF$ has finite index  
in $\Aut\big(\Gamma(\epsilon)\big)$. The rest follows from Proposition \ref{splita}.
\end{proof}

We are now ready for the proof of Theorem \ref{teob}.

\begin{prf} {\it Proof of Theorem \ref{teob}} Suppose that
$\Aut\big(\Gamma(\epsilon)\big)$ contains a subgroup of finite index which is an 
arithmetic group. We infer from iv) in Proposition \ref{specialgr} that    
$\Gamma(\hat A)$ is an 
arithmetic group, where $\hat A$ is as in \eqref{mamat}. 
We finish by the remark that $\hat A$ 
does not satisfy the necessary properties 
in Proposition \ref{splitprop}.
\end{prf}

Building on the above method it is possible to construct many more
examples of polycyclic groups $\Gamma$ with an
automorphism group $\Aut(\Gamma)$ which does not contain an arithmetic
subgroup of finite index. For example, as a slight variation of the above
groups $\Gamma(\epsilon)$, we may replace the dihedral group
$D_\infty$ by the non-trivial semi-direct product $D_1$ of $\bbZ$ with
itself, and let $D_1$ act on $F= \cO \times \bbZ$ 
via its natural homomorphism to $D_\infty$. 
We thus obtain a torsion-free, arithmetic polycyclic group $\Gamma_1(\epsilon)$
of rank five with non-arithmetic automorphism group. 
Another interesting class of examples may be constructed 
by starting with the (non-arithmetic) polycyclic
groups constructed in  \cite{GP4}. For these examples
the failure of arithmeticity is of rather different nature  
than in the groups  $\Gamma(\epsilon)$.


\section{Cohomology representations of $\Out(\Gamma)$}\label{topram1}


In this section we study the representation of $\Aut(\Gamma)$, $\Gamma$ a 
torsion-free polycyclic-by-finite group, on the cohomology groups 
$\H^*(\Gamma,R)$, where $R = \bbZ, \bbQ, \bbC$.
Since inner automorphisms act trivially on the cohomology of $\Gamma$,  
the outer automorphism group 
$\Out(\Gamma)$ is represented on the cohomology 
ring $\H^*(\Gamma,R)$. 
Considering the special case $R=\BC$, we find that the 
complex vector space $\H^*(\Gamma,\BC)$ comes with a natural 
$\BZ$-structure which is given by the image of the base change homomorphism 
$\H^*(\Gamma,\bbZ) \ra  \H^*(\Gamma, \bbC)$. Recall that this image is a
finitely generated subgroup containing a basis of $\H^*(\Gamma,\BC)$.
We fix here this $\BZ$-structure and its resulting $\BQ$-structure on 
$\H^*(\Gamma,\BC)$.
The representation of  $\Out(\Gamma)$ 
is an  \emph{integral representation} on $\H^*(\Gamma, \bbC)$, that is, 
$\Out(\Gamma)$ normalizes the
$\BZ$-lattice in $\H^*(\Gamma, \bbC)$ just described.
The $\BQ$-structure on  
$\H^*(\Gamma, \bbC)$ allows us to identify 
the group of invertible linear maps $\GLHC$ 
with a $\BQ$-defined linear algebraic group.
The Zariski-closure of the image 
of $\Out(\Gamma)$ in $\GLHC$ is a $\bbQ$-closed 
subgroup. We will show
that the representation of $\Out(\Gamma)$ on 
$\H^*(\Gamma, \bbC)$ 
is  an \emph{arithmetic representation}, 
that is, the image of $\Out(\Gamma)$ in $\GLHC$
is an arithmetic subgroup in its Zariski-closure.
In particular, this establishes 
our main results of Section \ref{topram}.

To carry over the information from the embedding of $\Out(\Gamma)$
into a linear algebraic group to 
topology and to the study of the cohomology $\H^*(\Gamma,R)$, 
we apply geometric methods originating from \cite{Baues2}.


\subsection{Automorphisms of Lie algebra cohomology}\label{liecohom}


An important special case in our theory is that of a finitely generated
torsion-free nilpotent group $\Gamma$. In this case, the cohomology of 
$\Gamma$ is intimately related to the Lie algebra cohomology of the Lie
algebra of the Malcev completion of $\Gamma$, see Section \ref{geometrc}.
We add here some well known facts about Lie algebra cohomology.  

Let $\lg$ denote a Lie algebra. 
The Lie product
of $\lg$ is expressed by a map $\varphi: \lg \wedge \lg \ra \lg$ which
satisfies the Jacobi-identity. The cohomology 
ring $\H(\lg)$ of $\lg$ is defined as the cohomology 
of the \emph{Koszul-complex} $\cK$ of $\lg$, cf.\
\cite{Koszul}. The complex $\cK$ has the structure 
of a differential graded algebra. 
As a graded algebra $\cK = \bigwedge \lg^*$ is the
exterior algebra of  the dual of $\lg$. The 
differential $d$ of $\cK$ is determined in degree one,
where $d:\lg^* \ra \bigwedge^2 \lg^*$ is defined as the 
dual of the Lie product $\varphi$. In particular, 
the cohomology of $\lg$
in degree one is computed as $\H^1(\lg) =
Z^1(\lg) =  [\lg, \lg]^\perp$. Note furthermore
that, via the duality, the automorphism group  
of the differential graded algebra $\cK$ 
identifies with the group of Lie algebra
automorphisms $\Aut(\lg)$. The automorphism 
group $\Aut(\lg)$ acts on the cohomology 
$H(\lg)$ with the inner automorphisms, generated 
by the exponentials of inner derivations  
of $\lg$ acting trivially. 
 
Assume now that $\lg$ is nilpotent. We consider the
descending central series of $\lg$ which is defined by $\lg^0 = \lg$,
$\lg^{i+1} = [\lg, \lg^i]$. Since $\lg$ is nilpotent,
$\lg^{k} = \{ 0\}$, for some (minimal) $k \in \bbN$.
Dualizing the descending central series,
we obtain a filtration $\lg_0 = \{0 \} \subset \lg_1 \ldots 
\subset \lg_k = \lg^*$, where $\lg_i = \left(\lg^i\right)^\perp$,
and $d \lg_{i} \subset \bigwedge^2 \lg_{i-1}$.
 
\begin{lemma} Let $\Phi$ be a semi-simple automorphism of the 
Koszul-complex of the nilpotent Lie algebra $\lg$. 
If $\Phi$ induces the identity on $\H^1(\lg)$ then 
$\Phi = \id$.
\end{lemma}  
\begin{proof} Since $\Phi$ is the identity on $\H^1$, it
is the identity on $\lg_1$. We prove by induction that
$\Phi$ is the identity on the subalgebra $\cK_j$, generated 
by $\lg_j$, $j>1$. Now $\lg_j$ is obtained 
from $\lg_{j-1}$ by 
adding finitely many generators $x \in \lg_j$. 
Since $\Phi$ is semisimple, $x$ may be chosen in a $\Phi$-invariant 
complement $W$ of $\lg_{j-1}$ in $\lg_j$. 
Since $d x \in  \cK_{i-1}$,
$d \Phi x =  \Phi dx = dx$ and
$d \left(\Phi x -x\right) = 0$. Since $W$ 
has no intersection with $\ker d^1 = \lg_1$,
this implies $\Phi x = x$. Therefore,
$\Phi$ is the identity on $\lg_{j}$, and hence 
on $\cK_j$.  
\end{proof}

We thus obtain the following result:  

\begin{proposition}\label{cohom1} Let $\lg$ be a nilpotent Lie algebra. 
Then the kernel of the natural representation of $\Aut(\lg)$ 
on the cohomology\/  $\H(\lg)$ is unipotent. In particular,
any reductive subgroup of $\Aut(\lg)$ acts faithfully
on $\H(\lg)$, even on $\H^1(\lg)$. 
\end{proposition}


\subsection{Computation of $\H^*(\Gamma,\bbC)$ via geometry
and $\Aut(\Gamma)$-actions}\label{geometrc} 


Before treating the general case,
we start by recalling some known facts which allow us to
compute the complex cohomology of a finitely generated torsion-free 
nilpotent group in terms of Lie algebra cohomology.
Let $\Theta$ be a finitely generated torsion-free nilpotent group, and
$\bU$ the complex Malcev-completion of $\Theta$,  
$\lu$ the Lie algebra of $\bU$. Thus $\Theta \leq \bU(\bbQ)$ and 
$M_{\Theta}= \Theta \backslash \bU({\bbR})$ is a smooth manifold  
which is an Eilenberg-Mac Lane
space of type $K(\Theta,1)$. In particular, there 
is a natural identification
of $\H^*(\Theta,\bbZ)$ with the singular cohomology group 
$\H^*(M_{\Theta},\bbZ)$, see the discussion later in this section. 
By de Rham's theorem,
the singular cohomology ring $\H^*(M_{\Theta},\bbC)$
of the smooth manifold $M_{\Theta}$ is isomorphic to the
cohomology $\H^*_{\rm DR}(M_{\Theta},\bbC)$ of
complex valued
$C^{\infty}$-differential forms on $M_{\Theta}$. 
In this situation, Nomizu \cite{Nomizu} proved that 
the natural map from $\lu$ into the differential forms 
on $\bU(\bbR)$  induces an  isomorphism of cohomology rings 
\begin{equation}\label{enili101}
{\mathsf n}: \,  \H^*(\lu, \bbC) \xrightarrow{\cong}  
\, \H^*_{\rm DR}(M_{\Theta},\BC) \; .
\end{equation}
Composing this map with 
the natural isomorphisms 
$$ \H^*_{\rm DR}(M_{\Theta},\BC)\to \H^*(M_{\Theta},\BC)\to
\H^*({\Theta},\bbC) $$ 
gives thus a linear isomorphism
\begin{equation}\label{enili102} 
{\it n_{\Theta}}: \H^*(\lu, \bbC) \xrightarrow{\cong}
\H^*({\Theta},\bbC) \; .
\end{equation}
Since $\Aut(\Theta)$ acts on 
$\bU(\bbR)$ by algebraic automorphisms, it also acts 
through smooth maps on $M_\Theta$.  Moreover, the
isomorphisms ${\it n}$ and ${\it n}_{\Theta}$ are compatible with 
the induced cohomology actions of $\Aut(\Theta)$ on  $ \H^*(\lu, \bbC) $, 
$\H^*_{\rm DR}(M_{\Theta},\BC)$ and $\H^*(\Theta,\bbC)$. 

A similar picture carries over to
our general situation where we start with a
torsion-free polycyclic-by-finite group $\Gamma$. 
We explain now some geometric constructions which 
extend the above picture from the 
case of torsion-free nilpotent groups 
to the more general situation. These constructions are closely 
connected with the algebraic setup discussed
so far in this paper. 
 
Let $\bH_{\Gamma}$ be the algebraic hull of 
$\Gamma$, $\bS$ a maximal $\BQ$-closed d-subgroup, 
and $\bU$ the unipotent radical of $\bH_\Gamma$. We have 
$\bH_{\Gamma}=\bU\cdot \bS$. 
We report from \cite{Baues2} the construction of the standard
$\Gamma$-manifold $M_\Gamma$. To construct this manifold we write a 
$\gamma\in\Gamma$ (uniquely) as 
$\gamma=u s$ with $u\in \bU(\BQ),\, s\in \bS(\BQ)$ 
and set
\begin{equation}\label{op12}
\gamma * x:= u s x s^{-1}=\gamma x s^{-1}\qquad ( x\in \bU(\BR) ).
\end{equation}
As noted in  \cite{Baues2}, this establishes a fixed-point-free, differentiable 
and properly discontinuous action of $\Gamma$ on $\bU(\BR)$. Moreover, 
the quotient space  
$$M_\Gamma=\Gamma \big \backslash \bU(\BR)$$
is a compact 
${\rm C}^\infty$-manifold and 
an Eilenberg-Mac Lane space of type $\KGone$.

We now explain how to calculate the complex 
cohomology of $M_\Gamma$, and hence the cohomology 
of $\Gamma$. 
Let $\lu$ denote the Lie algebra of $\bU$, and let 
$\cK_{\lu}$ be the Koszul-complex of $\lu$ (see Section \ref{liecohom}).
We let $\bS$ act by conjugation on $\bU$ and by the adjoint action on $\lu$
and $\cK_{\lu}$.
Let $\cK_{\lu}^\bS \subset \cK_{\lu}$ denote the
differential subcomplex of invariants for $\bS$. 
Now, as is proved  \cite[\S 3]{Baues2},  the obvious map 
\begin{equation}\label{epili101}
{\mathsf  n}:   \H^*(\cK_{\lu},\bbC)^\bS   \ra  \H^*_{\rm DR}(M_{\Gamma},\BC).
\end{equation}  
is an isomorphism of cohomology rings.

We explain next how $\Aut(\Gamma)$ acts on the above 
cohomology spaces. 
Let $\phi\in \Aut(\Gamma)$ and $\Phi$ be its extension to an algebraic
automorphism of $\bH_{\Gamma}$. 
We choose a $w_\phi\in\bU(\BQ)$ with the
property that $\Phi(\bS)= w_\phi\bS w_\phi^{-1}$ and set
\begin{equation}\label{op13}
X_\phi(x) =\Phi(x)w_\phi\qquad (x\in \bU(\BR)) \; .    
\end{equation}
This defines  a ${\rm C}^\infty$-map $X_\phi :\bU(\BR)\to \bU(\BR)$. A
straighforward computation yields that
\begin{equation}\label{op14} 
X_\phi(\gamma * x)=\phi(\gamma) * X_\phi(x).
\end{equation}
This shows that the map $X_\phi$ descends to a map 
$\bar X_\phi:M_\Gamma\to M_\Gamma$. 

Let
$${\rm Z}_\bU(\bS):=\{\, v\in \bU\Mid vsv^{-1}=s\ {\rm for \ all}\ s\in
\bS\,\}$$
be the centralizer of $\bS$ in $\bU$. This is a $\BQ$-closed subgroup of
$\bU$. Note that the similarly defined normalizer 
of $\bS$ in $\bU$ is in fact equal to ${\rm Z}_\bU(\bS)$. Multiplication from
the right defines an action of ${\rm Z}_\bU(\bS)(\BR)$ on $\bU(\BR)$ which 
commutes with the action of $\Gamma$ on $\bU(\BR)$ defined in (\ref{op12}). 
Since ${\rm Z}_\bU(\bS)(\BR)$ is connected, this action is homotopically
trivial, and so is the induced action on $M_\Gamma$.

Let $\phi,\, \psi\in \Aut(\Gamma)$ be automorphisms. A straighforward
computation shows that $X_\phi\circ X_\psi$ differs from    
$X_{\phi\circ\psi}$ by an element of ${\rm Z}_\bU(\bS)(\BR)$ acting on 
$\bU(\BR)$. This shows that $X_\phi\circ X_\psi$ and $X_{\phi\circ\psi}$
are homotopic maps, 
as well as the maps $\bar X_\phi\circ\bar X_\psi$ and $\bar X_{\phi\circ\psi}$. 
In particular, this implies that,   via the maps $\bar X_\phi$, 
$\phi\in \Aut(\Gamma)$, we obtain 
an action of the group $\Aut(\Gamma)$ on the cohomology spaces  
$\H^*_{\rm DR}(M_{\Gamma},\BC)$ and $\H^*(M_{\Gamma},\BC)$. 
The de Rham isomorphism 
\begin{equation}\label{deRa}
{\mathsf I}^*: \; \H^*_{\rm DR}(M_{\Gamma},\BC)\to\H^*(M_{\Gamma},\BC)
\end{equation} 
is obviously equivariant. Since $M_\Gamma$ is a $\KGone$, 
there is an isomorphism (compare
 \cite[Theorem 11.5]{Macl}), 
\begin{equation}\label{Macl} 
{\mathsf l} :\;  \H^*(M_{\Gamma},\BC)\to \H^*(\Gamma,\BC) \; .
\end{equation}
The isomorphism ${\mathsf l}$ is natural with respect to the pairs
$(\bar X_\phi,\phi)$. (For this, property (\ref{op14}) is essential,
as is explained in \cite{Macl}.) In particular, ${\mathsf l}$
is $\Aut(\Gamma)$-equivariant.

The group $\Aut_{a}(\bH_{\Gamma})$ of algebraic automorphisms of 
$\bH_{\Gamma}$ stabilizes the unipotent radical $\bU$, and hence it 
acts on $\lu$ and $\cK_{\lu}$. Since the group of inner automorphisms 
$\Inn_\bU$ acts
trivially on $\H^*(\cK_{\lu},\bbC)$,  we obtain an action of 
$\Aut_{a}(\bH_{\Gamma})= \Inn_\bU \cdot \Aut_a(\bH_{\Gamma})_\bS$ on 
 $\H^*(\cK_{\lu},\bbC)^\bS$. Here $\Aut_a(\bH_{\Gamma})_\bS$ stands for the
stabilizer of $\bS$ in $\Aut_{a}(\bH_{\Gamma})$.  
In particular, identifying $\Aut(\Gamma)$ as usually
with a subgroup of $\Aut_a(\bH_{\Gamma})$, 
this constructs a representation 
of $\Aut(\Gamma)$ on $\H^*(\cK_{\lu},\bbC)^\bS$. 
The isomorphism \eqref{epili101} 
can then be easily  seen to be 
$\Aut(\Gamma)$-equivariant. 

Let us define
\begin{equation}\label{epili102}
\mathsf{n}_{\Gamma} = {\mathsf l}  \circ {\mathsf I}^*  \circ {\mathsf n}: \;   
\H^*(\cK_{\lu},\bbC)^\bS   \, \ra  \, \H^*({\Gamma}, \bbC) \; .
\end{equation}
We have proved:
\begin{proposition}\label{equivi12}
Let $\Gamma$ be a torsion-free polycyclic-by-finite group. The isomorphism
${\mathsf n}_\Gamma: \H^*(\cK_{\lu},\bbC)^\bS  \rightarrow  \H^*({\Gamma}, \bbC)$ is equivariant 
with respect to the action of $\Aut(\Gamma)$ 
on $\H^*(\cK_{\lu},\bbC)^\bS$ (as defined above) and 
the natural action on $\H^*({\Gamma},\bbC)$. 
\end{proposition}  
 
\subsection{Rational action of $\Out_a(\bH_\Gamma)$ on $\H^*({\Gamma},
\bbC)$} \label{geometrc3}

As remarked before, the cohomology $\H^*({\Gamma},
\bbC)$ carries a natural $\bbQ$-structure 
induced by the coefficient homomorphism 
$\H^*(\Gamma,\bbQ) \ra \H^*(\Gamma,\bbC)$.
Thus, in particular, the group $\GLHC$ 
attains the natural structure of a $\bbQ$-defined group.
We discuss now the naturally defined $\bbQ$-structure on $\H^*(\cK_{\lu},\bbC)^\bS$.

Note first that the Lie algebra $\lu$ is defined over $\bbQ$. This means, there
exists a $\bbQ$-Lie algebra $\lu_{\bbQ}$ such that 
$\lu = \lu_{\bbQ} \otimes \BC$ is 
the scalar extension of $\lu_{\bbQ}$. 
The $\bbQ$-subalgebra $\lu_{\bbQ}$ is called a $\bbQ$-structure on $\lu$.
It is  induced by the $\bbQ$-strucure on $\bH_{\Gamma}$ (or,  equivalently,
by the unipotent shadow $\Theta$ of $\Gamma$) on $\lu$. (For related details
concerning $\bbQ$-structures on nilpotent Lie algebras and unipotent
groups one may consult  \cite{Segal,GO}.)   
Since $\lu$ is defined over $\BQ$,  we obtain a rational structure for the
vector space $\cK_{\lu}$ of the Koszul-complex of $\lu$, and 
all differentials are defined over $\BQ$.
Since $\bS$ is $\bbQ$-closed in $\bH_\Gamma$, 
it follows that the complex $\cK_{\lu}^\bS$ and its 
cohomology vector spaces $\H^*(\cK_{\lu}^\bS,\bbC)$
inherit a natural $\bbQ$-structure from  $\cK_{\lu}$, 
representing $\GLHC$  as a $\bbQ$-defined 
linear algebraic group.

Recall the construction of the $\bbQ$-defined
algebraic structure of $\Aut_{a}(\bH_{\Gamma})$ which
is discussed in Section \ref{subsect:algebraicstructure}. 
It is obtained by taking the natural
quotient $\Inn_\bU \rtimes \Aut_a(\bH_{\Gamma})_\bS 
\ra \Aut_{a}(\bH_{\Gamma})$. 
Since, by definition of its algebraic structure, the natural 
representation of   
$\Aut_a(\bH_{\Gamma})_\bS$ on $\lu$ is defined over $\bbQ$,  the natural 
representation of $\Aut_a(\bH_{\Gamma})$ on  
$\H^*(\cK_{\lu},\bbC)^\bS$ is $\bbQ$-defined as well. 
This also implies that the representation of $\Aut_{a}(\bH_{\Gamma})$
factors via a $\bbQ$-defined representation $\eta$ 
of $\Out_a(\bH_\Gamma)$. 
In particular, the induced representation 
of $\Aut(\Gamma)$ on $\H^*(\cK_{\lu},\bbC)^\bS$ factors through     
$\Out(\Gamma)$, and,  taking a basis with respect to the above 
constructed $\BQ$-structure, 
every element of  
$\Out(\Gamma)$ acts by a matrix with rational entries on  
$\H^*(\cK_{\lu},\bbC)^\bS$.

\begin{proposition}   \label{prop:cequivariance}
The natural 
representation of\/ $\Out(\Gamma)$ on $\H^*(\Gamma,\bbC)$
is arithmetic. 
\end{proposition} 
\begin{proof} 
By construction, $\Out(\Gamma)$ is contained in 
the $\bbQ$-defined linear algebraic group $\cO_\Gamma$
as a Zariski-dense arithmetic subgroup. Via the 
homomorphism $\pi_{\cO_\Gamma}: \cO_\Gamma \ra \Out_a(\bH_\Gamma)$
the representation of $\Out(\Gamma)$
on $\H^*(\cK_{\lu},\bbC)^\bS$ is induced by a $\bbQ$-defined
representation of $\cO_\Gamma$. Let $\rho_{\cO_\Gamma}$ denote the
representation of  $\cO_\Gamma$ on $\H^*(\Gamma,\bbC)$
which is obtained by conjugating with the isomorphism 
$\mathsf{n}_{\Gamma}$ defined in \eqref{epili102}. Let $\rho$
denote the natural representation of $\Out(\Gamma)$ on 
$\H^*(\Gamma, \bbC)$. Since $\mathsf{n}_{\Gamma}$ is $\Aut(\Gamma)$-equivariant,
$\rho_{\cO_\Gamma}\big(\Out(\Gamma)\big) = \rho\big(\Out(\Gamma)\big)$ consists of 
integral matrices in $\GLHC$ (with respect to a basis of 
$\H^*(\Gamma,\bbC)$ taken in the image
of $\H^*(\Gamma,\BZ)$).

We show that $\rho_{\cO_\Gamma}$ is a $\bbQ$-defined representation 
for the natural $\bbQ$-structure on $\H^*(\Gamma,\bbC)$.
As $\Out(\Gamma)$ is Zariski-dense in $\cO_\Gamma$ and consists of
rational points, it follows that the homomorphism $\rho_{\cO_\Gamma}$  maps
a Zariski-dense subset of rational points of  $\cO_\Gamma$
to rational points of  $\GLHC$. 
The Galois-criterion for rationality applied to the case of the extension 
$\bbC$ over $\bbQ$ implies that $\rho_{\cO_\Gamma}$ is defined over $\bbQ$.
By AR1, we infer that the image $\rho\big(\Out(\Gamma)\big)$ in $\GLHC$
is arithmetic in $\rho_{\cO_\Gamma}(\cO_\Gamma)$. 
\end{proof} 

Proposition  \ref{prop:cequivariance} proves 
Theorem \ref{topolo1} of the introduction. Remark that 
the kernel of $\rho$ is a finitely
generated group since it is an arithmetic subgroup of the kernel of
$\rho_{\cO_\Gamma}$.  

\begin{remark}
The proof of Proposition \ref{prop:cequivariance} gives no information
 about the rationality of the isomorphism ${\mathsf n}_{\Gamma}$.
In fact, the representation of $\cO_\Gamma$ might be trivial.
\end{remark}

For $\Gamma=\Theta$ a nilpotent group, some considerations on
rationality questions for the isomorphism ${\mathsf n}_\Theta$ 
may be found  in \cite{LP}. We state here:

\begin{proposition}   The isomorphism ${\mathsf n}_\Gamma$ is defined over $\bbQ$.
\end{proposition}
\begin{prf}{Sketch of proof.} We  show that the map 
${\mathsf  n}:   \H^*(\cK_{\lu},\bbC)^\bS   \ra  \H^*_{\rm DR}(M_{\Gamma},\BC)$ is defined over $\bbQ$. This can be seen as follows. 
In exponential coordinates for $\bU(\bbR)$, the forms in $\cK_{\lu}(\bbQ)$ define polynomial differential 
forms with rational coefficients. The de Rham isomorphism  ${\mathsf I}^*$ 
factorizes over the
cohomology of the piecewise differential forms, see \cite[VIII]{GM}.
The natural map from differential forms to piecewise-forms, maps the forms 
corresponding to elements of $\cK_{\lu}(\bbQ)$ into piecewise linear forms with rational coefficients.  
By Sullivan's  \emph{P.L.\ de Rham 
theorem} (see \cite{Sullivan, GM}) these are mapped into $\H^*(M_\Gamma, \bbQ)$. 
\end{prf}


\end{document}